\documentclass{article}

\usepackage{PRIMEarxiv}

\usepackage[utf8]{inputenc} 
\usepackage[T1]{fontenc}    
\usepackage{hyperref}       
\usepackage{doi}            
\usepackage{url}            
\usepackage{booktabs}       
\usepackage{amsfonts}       
\usepackage{amsthm}
\usepackage{nicefrac}       
\usepackage{microtype}      
\usepackage{lipsum}
\usepackage{fancyhdr}       
\usepackage{graphicx}       
\graphicspath{{media/}}     
\usepackage{xcolor}
\usepackage{mathtools}
\usepackage{bm}
\usepackage{subcaption}
\usepackage[square,numbers]{natbib}
\usepackage{soul}
\usepackage{tikz}
\usetikzlibrary{positioning, fit, calc, decorations.pathreplacing, shapes.geometric, arrows.meta, backgrounds}

\pgfdeclarelayer{background}
\pgfdeclarelayer{foreground}
\pgfsetlayers{background,main,foreground}

\usepackage{circuitikz}
\usepackage{xurl}
\usepackage{amssymb}
\usepackage{multirow}

\usepackage{tablefootnote}

\usepackage{pifont}
\usepackage{xcolor}

\usepackage[acronym]{glossaries}
\newacronym{r3bp}{R3BP}{Restricted Three-Body Problem}
\newacronym{r3bps}{R3BPs}{Restricted Three-Body Problems}
\newacronym{cr3bp}{CR3BP}{Circular Restricted Three-Body Problem}
\newacronym{er3bp}{ER3BP}{Elliptic Restricted Three-Body Problem}
\newacronym{bcr4bp}{BCR4BP}{Bi-Circular Restricted Four-Body Problem}
\newacronym{hr3bp}{HR3BP}{Hill Restricted Three-Body Problem}
\newacronym{hr4bp}{HR4BP}{Hill Restricted Four-Body Problem}
\newacronym{qbcp}{QBCP}{Quasi Bi-Circular Problem}
\newacronym{qhr4bp}{QHR4BP}{Quasi-Hill Restricted Four-Body Problem}
\newacronym{qhr4bps}{QHR4BPs}{Quasi-Hill Restricted Four-Body Problems}

\newacronym{iqhr4bp}{I-QHR4BP}{In-plane Quasi-Hill Restricted Four-Body Problem}
\newacronym{oqhr4bp}{O-QHR4BP}{Out-of-plane Quasi-Hill Restricted Four-Body Problem}

\newacronym{hfem}{HFEM}{Higher-Fidelity Ephemeris Model}

\newacronym{nrho}{NRHO}{Near Rectilinear Halo Orbit}

\newacronym{qpo}{QPO}{Quasi-Periodic Orbit}
\newacronym{qpos}{QPOs}{Quasi-Periodic Orbits}

\newacronym{myt}{MYT}{Multi-Year Trajectory}
\newacronym{myts}{MYTs}{Multi-Year Trajectories}

\newacronym{gmos}{GMOS}{G{\'o}mez Mondelo Olikara Scheeres}

\newacronym{po}{PO}{Periodic Orbit}
\newacronym{pos}{POs}{Periodic Orbits}

\newacronym{dft}{DFT}{Discrete Fourier Transform}
\newacronym{fft}{FFT}{Fast Fourier Transform}
\newacronym{cft}{CFT}{Continuous Fourier Transform}

\newacronym{dro}{DRO}{Distant Retrograde Orbit}
\newacronym{lpo}{LPO}{Low Prograde Orbit}
\newacronym{dpo}{DPO}{Distant Prograde Orbit}
\newacronym{stm}{STM}{State Transition Matrix}
\newacronym{vsr}{VSR}{Vertical Self-Resonant}

\newacronym{qpts}{QPTs}{Quasi-Periodic Trajectories}

\newacronym{elfos}{ELFOs}{Elliptical Lunar Frozen Orbits}
\newacronym{elfo}{ELFO}{Elliptical Lunar Frozen Orbit}
\newacronym{qelfo}{Q-ELFO}{Quasi-Elliptical Lunar Frozen Orbit}
\newacronym{qelfos}{Q-ELFOs}{Quasi-Elliptical Lunar Frozen Orbit}
\newacronym{dadm}{DADM}{Doubly Averaged Dynamical Model}

\newacronym{laskar}{L-NAFF}{Laskar-Numerical Analysis of Fundamental Frequency}
\newacronym{gomez}{GMS-C}{G{\'o}mez-Mondelo-Sim{\'o}-Collocation}

\newacronym{fli}{FLI}{Fast Lyapunov Indicator}

\newacronym{fddc}{FDDC}{Frequency-Domain Differential Corrector}

\newacronym{brf}{BRF}{Barycentric Rotating Frame}
\newacronym{eof}{PIF}{$\mathcal{P}_2$ Inertial Frame}
\newacronym{mci}{JIF}{J2000 Inertial Frame}
\newacronym{hrf}{HRF}{Hill Rotating Frame}

\newacronym{raan}{RAAN}{Right Ascension of the Ascending Node}
\glsdisablehyper
\definecolor{threeCat}{rgb}{1.00,0.00,0.00} 
\definecolor{nuOneThree}{rgb}{1,0,1}
\definecolor{nuTwoThree}{rgb}{0,1,1}

\definecolor{cf}{rgb}{1,0,1}
\definecolor{ccu}{rgb}{0,1,1}
\definecolor{ccs}{rgb}{1,0,0}

\newcommand{\mb}{\begin{bmatrix}}
\newcommand{\me}{\end{bmatrix}}

\newcommand{\da}[1]{\ensuremath{#1}}

\newcommand{\po}{\ensuremath{\mathcal{P}_1}}
\newcommand{\pt}{\ensuremath{\mathcal{P}_2}}

\newcommand{\TT}{\mathbb{T}}

\definecolor{goals}{rgb}{0.00,0.00,1.00}

\definecolor{southern}{rgb}{0.00, 1.00, 1.00}
\definecolor{northern}{rgb}{0.00, 0.00, 1.00}
\definecolor{ascending}{rgb}{1.00, 0.00, 0.00}
\definecolor{descending}{rgb}{1.00, 0.00, 1.00}

\tikzset{
    markerbase/.style={
        regular polygon,
        regular polygon sides=3,
        draw=black,
        line width=0.4pt,
        inner sep=0pt,
        minimum size=1.0em, 
        baseline=-0.3ex
    }
}

\newcommand{\xp}[1]{%
  \tikz\node[markerbase, fill=#1, shape border rotate=-90] {};%
}

\newcommand{\xm}[1]{%
  \tikz\node[markerbase, fill=#1, shape border rotate=90] {};%
}

\newcommand{\yp}[1]{%
  \tikz\node[markerbase, fill=#1, shape border rotate=0] {};%
}

\newcommand{\ym}[1]{%
  \tikz\node[markerbase, fill=#1, shape border rotate=180] {};%
}

\newcommand{\ID}[2]{\ensuremath{\mathrm{#1}_{#2}}}

\usepackage{colortbl}
\definecolor{prune_light}{gray}{0.90} 
\definecolor{prune_dark}{gray}{0.70}  

\usepackage{enumitem}
  
\title{\Large{Linking Averaged and Unaveraged Three-Body Dynamics\\ Near Smaller Primaries: Symmetric Periodic Orbits}}

\author{
  Beom Park\\
  Apollo 11 Postdoctoral Fellow \\
  School of Aeronautics and Astronautics \\
  Purdue University \\
  West Lafayette, IN, USA, 47907 \\
  \texttt{park1103@purdue.edu} \\
  \And
  Kathleen C. Howell \\
  Hsu Lo Distinguished Professor \\ 
  School of Aeronautics and Astronautics \\
  Purdue University \\
  West Lafayette, IN, USA, 47907 \\
  \texttt{howell@purdue.edu}
}
\begin{document}
\maketitle


\begin{abstract}
Within a three-body system comprised of two celestial bodies and a spacecraft, the dynamical environment near a smaller primary is significantly perturbed, motivating a balance between global insight and model fidelity. While averaged dynamics offer an integrable model to classify solution landscapes, they inherently lack the accuracy of the unaveraged dynamics, such as the \acrfull{hr3bp} and \acrfull{cr3bp}.
 This work establishes a systematic bridge between the averaged and unaveraged regimes by explicitly linking averaged equilibria to symmetric periodic orbits in the unaveraged three-body systems. A unified frequency framework is introduced to characterize the mapping of invariant tori across the dynamical models. Leveraging the parity of the resonance ratio, an initialization scheme is developed to identify admissible apse configurations, enabling the \textit{a priori} prediction of solution multiplicity and symmetry types. Furthermore, the global evolution of families derived from averaged equilibria is traced via bifurcation and frequency analysis. These findings are synthesized into archetypical bifurcation diagrams, providing a comprehensive atlas of the symmetric periodic orbit web within the \acrshort{hr3bp} and \acrshort{cr3bp}. The resulting framework not only clarifies the topological origins of complex periodic orbit families but also offers a versatile tool for trajectory design in cislunar and multi-body environments.
\end{abstract}


\section{Introduction}

Many smaller primaries within various planet-moon systems are important targets in space exploration. For example, Earth’s Moon remains the focal point of near-term exploration, envisioned as a platform for science, a logistics hub for deep-space missions, and a potential outpost for both robotic and human operations \cite{fuller2022gateway}. Icy satellites of Jupiter and Saturn may host subsurface oceans and thus carry exceptional scientific value for understanding the origins of life and the possibility of life beyond Earth \cite{grasset2013jupiter, spencer2013enceladus}.

To support mission concepts involving smaller primaries, ``orbiters,'' i.e., ballistic trajectories that remain bounded near a smaller primary, are especially valuable. Yet around a smaller body, a purely Keplerian model centered on that body is often inadequate; third-body gravity from the larger primary can be comparable to, or exceed, the influence from the central body. A natural next step is the \acrfull{r3bp} \cite{szebehely2012theory} that incorporates the third body explicitly. There are two principal approaches to treat the third-body perturbation within this context. The first is a Keplerian-based perturbative approach that averages the third-body potential first over the satellite’s orbital period with respect to the small primary and then over the perturber’s period, yielding the doubly averaged model \cite{von1910application, lidov1961evolution, kozai1962secular}\footnote{Originally introduced by von Zeipel \cite{von1910application} in 1910, the method was independently re-discovered by Lidov \cite{lidov1961evolution} in 1961 and Kozai \cite{kozai1962secular} in 1962. For a thorough historical background and summary of the doubly-averaged model, refer to Ito and Ohtsuka \cite{ito2019lidov}. Modern developments of the doubly averaging approach appear in \cite{scheeres2001stability,de2003third,broucke2003long}.}. The averaging process renders the reduced dynamics integrable, enabling a clean classification of equilibria and phase-space structures, notably the librating and circulating structures in terms of the argument of periapsis \cite{de2003third} \cite{broucke2003long}. These orbits guide numerous mission designs targeting smaller primaries, for instance, high-latitude ``frozen'' solutions in supplying persistent polar coverage at the Moon \cite{ely2005stable,folta2006lunar,ely2006constellations} and at Europa \cite{lara2007computation}, as well as catalogs and use cases for circulating orbits \cite{russell2009circulating,mcardle2021circulating,franz2022database}.

The complementary route retains the unaveraged, non-integrable dynamics of the three-body formulation, rendering the \acrfull{hr3bp} and \acrfull{cr3bp}. While these unaveraged models are non-integrable, the insights and computational strategies for the Hamiltonian system are highly effective; Lagrange equilibria, (quasi-)periodic orbits, and their invariant manifolds organize the phase space. In particular, periodic orbits form one-parameter families in the \acrshort{r3bps}, enabling routine numerical continuation along branches; monodromy-based stability and bifurcation analyses then reveal a web of inter-family connections. Examples include numerical strategies for continuation and Earth-Moon applications \cite{howell1984almost,zimovan2020near}, global searches for periodic orbits \cite{russell2006global,restrepo2018database} and, more recently, symplectic methods tailored to the unaveraged dynamics that sharpen bifurcation mapping across families \cite{moreno2024bifurcation,aydin2025exploration,aydin2025studying}.

These two strategies are complementary. The averaged approach renders the reduced system integrable, but its accuracy degrades as the effective semi-major axis grows; once the third-body gravity becomes comparable to, or larger than, the central gravity, a Keplerian baseline around the small primary is no longer adequate, and double averaging rapidly loses fidelity. Moreover, integrability removes the natural transport mechanisms, i.e., stable/unstable manifolds, that organize global motion within a realistic model. The unaveraged models address precisely these gaps: the third-body perturbation is represented explicitly, and although the system is non-integrable (so many invariant tori of the averaged model are destroyed), their remnants may persist as invariant tori, whose families and invariant manifolds supply the backbone of the solution space.

The current effort seeks a unified characterization of third-body environments by explicitly linking the averaged and unaveraged dynamics. The idea is not new; multiple prior efforts successfully generate periodic orbits in the unaveraged dynamics by seeding from resonance-based orbital elements derived from averaged models \cite{howell2007design,russell2007long,mcardle2021circulating, lara2005dynamic}. However, these approaches are typically pragmatic rather than systematic; they enforce commensurability and use numerical correction schemes to locate specific solutions, but often lack a rigorous formalization of the mapping in terms of symmetries and frequency constraints. Consequently, fundamental questions remain regarding the completeness of the solution space. Specifically, the total number of expected branches for a given resonance and the topological origins of their bifurcations are often left unexplained. In a parallel thread, extensive continuations of periodic orbits in the \acrshort{r3bps} reveal branches resembling structures that originate from the averaged dynamics (e.g., \cite{lara2007classification,koblick2025novel,peng2025analog, zhang2025time}). Yet, most references stop short of explicitly grounding these branches in the underlying averaged invariant structures. Without this explicit link, the identification of the specific ``generating structure'' for a given family and the determination of physically realizable symmetries remain elusive.

The present contribution addresses this gap by focusing on the \emph{symmetric periodic orbits} within the unaveraged dynamics that originate from the averaged equilibria. A rigorous taxonomy regarding the configurations, solution multiplicity (number of branches), and symmetry hierarchy is established. Methodologically, this framework yields reliable, reproducible initialization schemes for symmetric periodic orbit families within the \acrshort{r3bps}. To achieve these objectives, the investigation is structured around the following primary goals:
\vspace{-\parskip}
\begin{enumerate}[label= \color{goals}(Goal \arabic*), leftmargin=3mm, itemindent=12mm, labelsep=!, labelwidth=2mm, labelindent=12mm]
    \item\label{goal1} \ul{Establish a common frequency framework:} Place averaged and unaveraged dynamics on a unified frequency lattice to consistently define the resonance conditions and the corresponding resonant periodic orbits.
    \item\label{goal2} \ul{Identify admissible symmetric configurations and solution multiplicity:} Recognize that only specific phase configurations within the averaged dynamics yield symmetric periodic orbits in the \acrshort{r3bps}. By leveraging the parity of the resonance ratio, this framework determines the physically admissible apse alignments, thereby unambiguously prescribing the exact number of solution branches and their symmetry types \textit{a priori}. This analysis spans both the \acrshort{hr3bp} and \acrshort{cr3bp} to ensure a unified classification of the symmetry evolution.
    \item\label{goal3} \ul{Trace the global evolution of symmetric periodic orbit families:} Map the averaged circular and frozen equilibria to their unaveraged counterparts, establishing explicit connections to well-known structures such as halo orbits. Through stability and frequency analysis, this objective quantifies the dynamical connectivity between these structures, elucidating the manifestation of ideal bifurcations (e.g., pitchfork) predicted by the averaged model as either continuous transitions or disconnected, ``broken'' branches in the \acrshort{r3bps}. Together, this analysis aids in predicting the number of distinct family branches originating from a given averaged equilibrium. 
    \item\label{goal4} \ul{Construct archetypical bifurcation diagrams within the unaveraged dynamics:} Synthesize a comprehensive atlas of symmetric families by grounding them in the underlying averaged equilibria. This framework explicitly illustrates the global connectivity and the exact number of family branches of different symmetry types, effectively organizing the complex topology of the \acrshort{r3bps} based on the structure of the integrable limit.
\end{enumerate}
\vspace{-\parskip}
Overall, a symmetry-aware correspondence is established between averaged equilibria and the resonant periodic orbit skeleton of the \acrshort{cr3bp}, turning integrable structures into a concrete bifurcation atlas in the unaveraged setting. This unified framework not only enhances the fundamental understanding of the dynamical environment near smaller primaries but also offers versatile utility for rapid trajectory design. The remainder of this paper is organized as follows: Section~\ref{sec:dynamics} details the dynamical formulations and coordinate transformations utilized throughout the work. The subsequent sections, Section~\ref{sec:dadm} through Section~\ref{sec:archetype}, are structured to independently address \textcolor{blue}{(Goal 1)} through \textcolor{blue}{(Goal 4)}, respectively, systematically building the connections from the averaged to the unaveraged dynamics.



\section{\label{sec:dynamics}Frames and Dynamical Models}

This section introduces the foundational concepts essential for the analysis in subsequent sections. Different reference frames and their associated state transformations serve as the basis, along with dynamical models reflecting varying fidelity levels. A common notation is first established. The larger and smaller primary bodies are denoted as $\po$ and $\pt$, respectively. Their dimensional standard gravitation parameters are $\tilde{\mu}_1$ and $\tilde{\mu}_2$. The mass ratio is defined as $\mu = \tilde{\mu}_2/(\tilde{\mu}_1 +\tilde{\mu}_2)$. All dynamical models under investigation assume a constant distance $l_*$ between the primaries. The characteristic time is defined as $t_* = \sqrt{l_*^3/(\tilde{\mu}_1 + \tilde{\mu}_2)}$. These two quantities serve as characteristic units that normalize the distance and time variables. The nondimensional (nd) time variable is denoted as $t$, serving as the common independent variable throughout the current work. Differentiation with respect to $t$ is denoted by an overdot, i.e., $\dot{(\cdot)} = d(\cdot)/dt$. 

\subsection{Frames and Transformations}

\begin{figure}[htpb]
    \centering
    \includegraphics[width=0.30\linewidth]{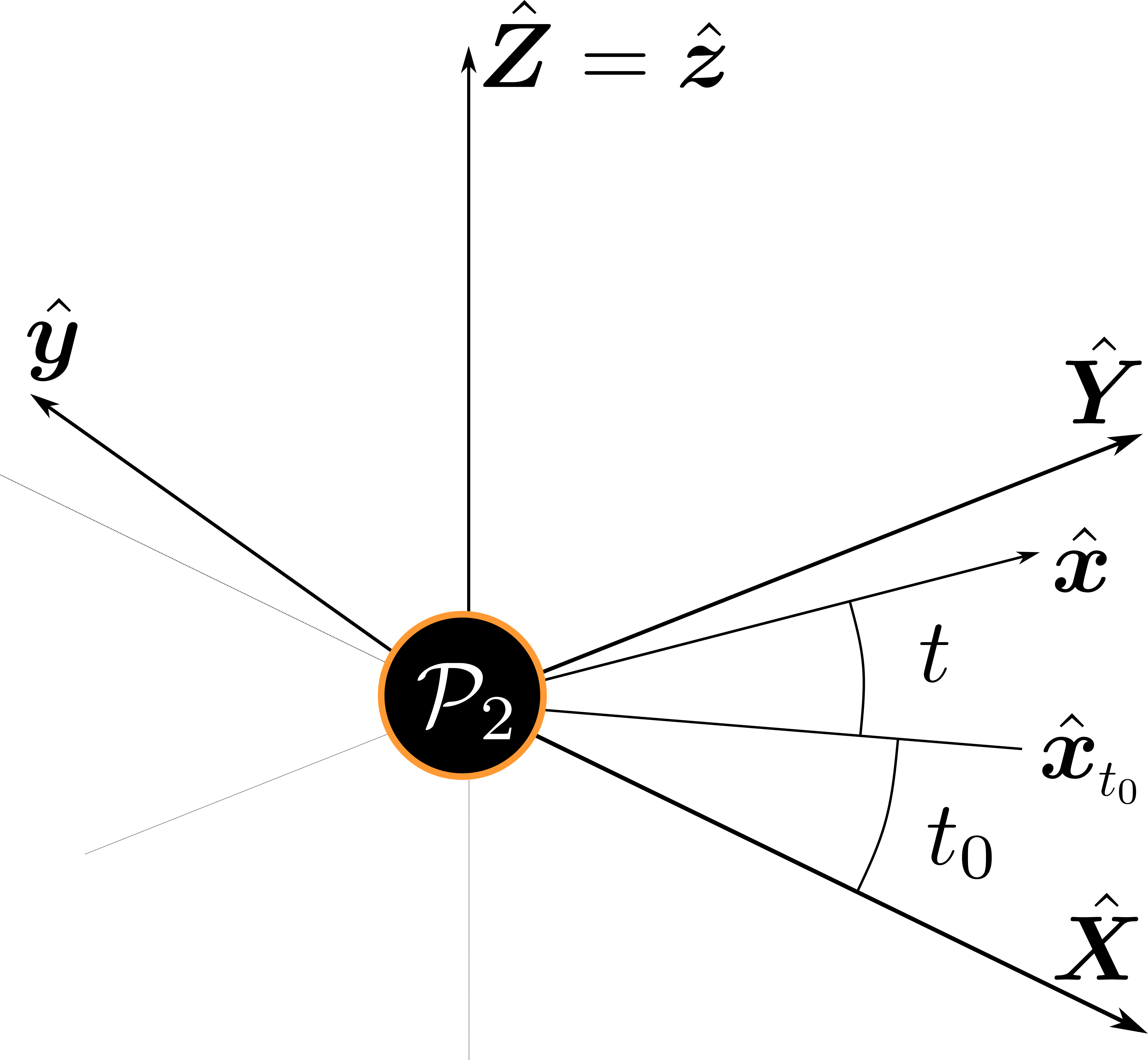}
    \caption{Geometry of the basis vectors for the inertial and rotating frames.}
    \label{fig:frames}
\end{figure}

Three different frames are employed in the current work. These various frames facilitate analyses of third-body gravity across models as follows:
\begin{itemize}[topsep=0pt, partopsep=0pt]
    \item \ul{Barycentric Rotating Frame} (\acrshort{brf}): This frame adopts an orthogonal basis $\hat{\bm{x}}-\hat{\bm{y}}-\hat{\bm{z}}$ defined as: (1) $\hat{\bm{x}}$ directed from $\po$ to $\pt$, (2) $\hat{\bm{z}}$ coincides with the angular momentum vector of $\pt$ with respect to $\po$, and (3) $\hat{\bm{y}}$ completes the dextral triad. The origin is located at the $\po$-$\pt$ barycenter. With the $\po$-$\pt$ distance normalized to unity, $\po$ and $\pt$ are fixed on the $\hat{\bm{x}}$ axis at $-\mu$ and $(1-\mu)$, i.e., $\bm{r}_1 = -\mu\hat{\bm{x}}$ and $\bm{r}_2 = (1-\mu)\hat{\bm{x}}$, respectively. The nd spacecraft position vector is $\bm{r} = x\hat{\bm{x}} + y\hat{\bm{y}} + z\hat{\bm{z}}$, and the nd velocity vector is $\bm{v} = d\bm{r}/dt = \dot{x}\hat{\bm{x}} + \dot{y}\hat{\bm{y}} + \dot{z}\hat{\bm{z}}$.
\item \ul{$\pt$-Centered Hill Rotating Frame} (\acrshort{hrf}): This frame shares the orthogonal basis $\hat{\bm{x}}-\hat{\bm{y}}-\hat{\bm{z}}$ with the \acrshort{brf} but is centered at $\pt$. The units are renormalized by $\mu^{1/3}$; thus, the position vector is $\bm{r}_H = x_H\hat{\bm{x}} + y_H\hat{\bm{y}} + z_H\hat{\bm{z}}$, where $x_H = (x-(1-\mu))/\mu^{1/3}$, $y_H = y/\mu^{1/3}$, and $z_H = z/\mu^{1/3}$. This scaling facilitates consistent analysis across systems with varying $\mu$. The velocity vector follows such that $\bm{v}_H = d\bm{r}_H/dt = \bm{v}/\mu^{1/3}$. The term ``rotating frames'' denotes both the \acrshort{brf} and \acrshort{hrf}. The \acrshort{hrf} is mainly adopted for visualizing trajectories in the current analysis.
\item \ul{$\pt$-Centered Inertial Frame} (\acrshort{eof}): The orthonormal basis $\hat{\bm{X}}-\hat{\bm{Y}}-\hat{\bm{Z}}$ is defined such that $\hat{\bm{Z}}$ aligns with $\hat{\bm{z}}$, while $\hat{\bm{X}}-\hat{\bm{Y}}$ are inertial. The rotating frame ($\hat{\bm{x}}-\hat{\bm{y}}$) rotates at a constant unit rate relative to this inertial frame. The origin is located at $\pt$. The spacecraft position vector in the \acrshort{eof} is constructed as:
    \begin{equation}
        \label{eq:eof_brf} 
        \bm{R} = \bm{C}_{\hat{\bm{z}}} (\bm{r}-\bm{r}_2) = X\hat{\bm{X}} + Y\hat{\bm{Y}} + Z\hat{\bm{Z}}, 
        \quad 
        \bm{C}_{\hat{\bm{z}}} = 
        \begin{bmatrix} 
        \cos (t + t_{0}) & -\sin (t + t_{0}) & 0 \\ 
        \sin (t + t_{0}) & \cos (t + t_{0}) & 0 \\ 
        0 & 0 & 1 
        \end{bmatrix}.
    \end{equation}
    Figure \ref{fig:frames} illustrates the unit vectors for the \acrshort{eof} and rotating frames; vector $\hat{\bm{x}}_{t_0}$ denotes the direction of the $\po$-$\pt$ line within the \acrfull{eof} at the initial epoch, $t_{0}$. This frame definition is consistent with prior studies \cite{broucke2003long, ely2005stable, folta2006lunar}, assuming the motion of $\po$ is confined to the $XY$-plane. The time derivative of Eq. \eqref{eq:eof_brf} yields the inertial velocity $\bm{V} = d\bm{R}/dt$. 
\end{itemize}
The states within the inertial frame may be represented via classical Keplerian elements: $a$ (semi-major axis), $e$ (eccentricity), $i$ (inclination), $\omega$ (argument of periapsis), $\Omega$ (right ascension of the ascending node, RAAN), and $M$ or $\theta$ (mean or true anomaly). Standard transformations (e.g., \citet{vallado2001fundamentals}) convert these elements to inertial position and velocity ($\bm{R}, \bm{V}$) that are then mapped to the rotating frames ($\bm{r}, \bm{v}$) via Eq. \eqref{eq:eof_brf} and its derivative.

\subsection{Dynamical Models}

Three dynamical models are investigated. The first model, termed the \acrfull{dadm}, represents the averaged dynamical model. The other two represent unaveraged dynamics in the form of the \acrshort{r3bp}, namely, the \acrshort{hr3bp} and \acrshort{cr3bp}. These models are defined in order.

\subsubsection{\acrfull{dadm}}

Within the $\pt$-dominated dynamical regime, several assumptions are introduced to supply the \acrshort{dadm} employed in the current investigation. Firstly, $\po$ is the only perturbing force. The point-mass bodies $\po$ and $\pt$ are orbiting around their mutual barycenter in circular orbits. While it is possible to extend the \acrshort{dadm} by relaxing these assumptions (e.g., including perturbations from the nonspherical gravity of $\pt$), such extensions remain out of scope for the current analysis. The perturbing potential of $\po$ is truncated leveraging the Legendre expansion as detailed by Longuski et al. \cite{longuski2022introduction}. Subsequently, the perturbation from $\po$ is averaged over (1) the period of $\po$ with respect to $\pt$ (medium-period) within the \acrshort{eof} and (2) the Keplerian orbital period of the spacecraft with respect to $\pt$ (short-period) within the \acrshort{eof}. This process yields the doubly-averaged, truncated $\po$ perturbing potential\footnote{In some work, an approximation $1-\mu \approx 1$ \cite{folta2006lunar} is leveraged. In the current analysis, this approximation is not introduced, and the factor $1-\mu$ is retained in all terms, following derivations by Broucke \cite{broucke2003long}.}, 
\begin{align}
    \label{eq:perturbing_potential} U_{\po} = \frac{1}{32}(1-\mu)n_1^2\da{a}^2 \left[ (1+3\cos 2 \da{i} )(2 + 3\da{e}^2) + 30\da{e}^2 \sin^2 \da{i} \cos 2 \da{\omega} \right],
\end{align}
where $n_1 = 1/t_*$ is the dimensional mean angular rate of $\po$ relative to $\pt$. The orbital elements in Eq. \eqref{eq:perturbing_potential} are defined within the \acrshort{eof}. From Lagrange's planetary equations, these elements evolve at the following rates with respect to the nd time $t$,
\begin{align}
    \label{eq:dadt} & \dot{a} = 0 \\
    \label{eq:dedt}& \dot{e} = \dot{e}(\da{e}, \da{i}, \da{\omega}) = \frac{15}{8}\frac{(1-\mu)n_1}{\da{n}}\da{e}(1-\da{e}^2)^{1/2}\sin^2 \da{i}\sin 2 \da{\omega}  \\
    \label{eq:didt}& \dot{i}= \dot{i}(\da{e}, \da{i}, \da{\omega}) = - \frac{15}{16}\frac{(1-\mu)n_1}{\da{n}}\frac{\da{e}^2}{(1-\da{e}^2)^{1/2}}\sin 2\da{i} \sin 2\omega \\
    \label{eq:domegadt}& \dot{\omega}= \dot{\omega}(\da{e}, \da{i}, \da{\omega}) = \frac{3}{16}\frac{(1-\mu)n_1}{\da{n}}\frac{1}{(1-\da{e}^2)^{1/2}} \left[ (3 + 2\da{e}^2 + 5\cos 2\da{i}) + 5(1-2\da{e}^2 - \cos 2\da{i} ) \cos 2\da{\omega} \right] \\
   \label{eq:dOmegadt}   & \dot{\Omega}= \frac{3}{8}\frac{(1-\mu)n_1}{\da{n}}\frac{1}{(1-\da{e}^2)^{1/2}} (5\da{e}^2 \cos 2\da{\omega} - 3\da{e}^2 - 2) \cos \da{i} \\
   \label{eq:dthetadt}  & \dot{M} = \da{n}t_* = \sqrt{\frac{\tilde{\mu}_2}{\da{a}^3}}t_* = \sqrt{\frac{l_*^3}{\da{a}^3}\mu} = \sqrt{\frac{1}{a_H^{3}}},
\end{align}
where $\da{n}$ is the dimensional mean motion for the satellite relative to $\pt$ supplied with $\da{a}$, the dimensional semi-major axis. For additional insights, the nd semi-major axis is also introduced as,
\begin{align}
   \label{eq:nd_a} \da{a}_H = \frac{\da{a}}{l_* \mu^{1/3}},
\end{align}
scaled to be consistent with the unit as defined in the \acrshort{hrf}. 

\subsubsection{\acrfull{hr3bp}}

Within the \acrshort{hr3bp}, the perturbing potential for $\po$ corresponds to the second-order term of the Legendre expansion; however, position and velocity components within the \acrshort{hrf} are directly leveraged without the averaging process. A defining characteristic of the \acrshort{hr3bp} is the limiting assumption that $\mu \to 0$, representing the dynamics infinitely close to the secondary primary with a negligible mass. Despite this assumption, the \acrshort{hr3bp} serves as a suitable intermediate model between the \acrshort{dadm} and higher-fidelity dynamical models \cite{mcardle2021circulating}. The equations of motion are rendered as,
\begin{align}
    \label{eq:hr3bp} 
    \dot{\bm{v}}_H = -2\hat{\bm{z}}  \times {\bm{v}_H}+\nabla U_H,
\end{align}
where the vector cross product is denoted as $\times$, and the potential function is $U_H = \frac{1}{2}(3x_H^2 - z_H^2 ) + 1/r_H$ with $r_H = |\bm{r}_H|$.

\subsubsection{\acrfull{cr3bp}}

The \acrshort{cr3bp} dynamics are formulated within the \acrshort{brf}. The equations of motion are given by,
\begin{align}
    \label{eq:cr3bp} 
    \dot{\bm{v}} = -2\hat{\bm{z}} \times \bm{v} + \nabla U.
\end{align}
Here, the pseudo-potential function $U$ is defined as $U = \frac{1}{2}(x^2 + y^2) + \frac{1-\mu}{d} + \frac{\mu}{r}$, where $d = |\bm{r} - \bm{r}_1|$ and $r = |\bm{r} - \bm{r}_2|$. In contrast to the \acrshort{hr3bp}, the \acrshort{cr3bp} incorporates the full gravitational influence of $\po$ without truncation. As $\mu \to 0$, the dynamics in the vicinity of $\pt$ asymptotically approach the \acrshort{hr3bp}. While the model is applicable to general three-body systems, this work specifically adopts $\mu \approx 0.01215$ that represents the Earth-Moon system.

\section{\label{sec:dadm}Invariant Tori and Frequency Mapping in the Averaged Dynamics}

Initially, the solution space within the averaged model, i.e., \acrshort{dadm}, is reviewed. Inspection of Eqs. \eqref{eq:dadt}-\eqref{eq:dthetadt} reveals three independent integrals of motion \cite{lidov1961evolution, broucke2003long},
\begin{align}
    \label{eq:C_Z} C_Z &= (1-e^2)\cos^2 i\\
    \label{eq:C_H} C_H &= e^2(2/5 - \sin^2 i \sin^2 \omega) \\
    \label{eq:C_a} C_a &= a,
\end{align}
related to the out-of-plane angular momentum, total energy, and $\da{a}$ with respect to $\mathcal{P}_2$, respectively. Given these three integrals for the three-degree-of-freedom Hamiltonian dynamics (see Nie and Gurfil \cite{nie2018lunar}), the \acrshort{dadm} is completely integrable in the sense of Liouville-Arnold \cite{arnol2013mathematical}. Consequently, any motion within the \acrshort{dadm} is confined to invariant tori, ensuring long-term stability \cite{celletti2023infinite}.

\subsection{Equilibria in the Reduced Space}

Analysis proceeds in a reduced space because the dynamics for $e, i, \omega$ (Eqs. \eqref{eq:dedt}-\eqref{eq:domegadt}) are separable from those of $a, \Omega, M$. Since $a$ is constant, the first two integrals ($C_Z, C_H$) depend solely on $e, i, \omega$. Hence, the triad $(e, i, \omega)$ traces either fixed points (equilibria) or periodic loops within the reduced phase space. To avoid the general singularity for $\omega$ at $e = 0$, the nonsingular Laplace-Runge-Lenz vector components $(k, h) = (e \cos{\omega}, e\sin\omega)$ are adopted. In these coordinates, Eqs. \eqref{eq:dedt} and \eqref{eq:domegadt} are recast as,
\begin{align}
    \label{eq:dkdt} \frac{dk}{dt} &= \frac{3}{16}K \left[ \frac{20 h k^2 }{e^2 }  \sin^2 i \sqrt{1-e^2} - \frac{h}{1-e^2} \Phi(k,h,i) \right] \\
    \label{eq:dhdt} \frac{dh}{dt} &= \frac{3}{16}K \left[ \frac{20 h^2 k }{e^2 }  \sin^2 i \sqrt{1-e^2} + \frac{k}{1-e^2} \Phi(k,h,i) \right],
\end{align}
where $e^2 = k^2 + h^2$ and $\Phi(k,h,i) = (3+2e^2 + 5\cos 2i) + 5(1-2e^2 -\cos 2i) \frac{k^2 -h^2}{e^2} $. The constant coefficient is defined as $K = (1-\mu)n_1/n$. Two types of equilibria for these elements are reviewed below.

\subsubsection{Circular Equilibria ($\da{e} = 0$)}

 The first type of equilibria, denoted the circular equilibria, occurs at $\da{e} = 0$; as such, $k = h = 0$ and $\da{i}$ is constant since $d\da{i}/dt = 0$ from Eq. \eqref{eq:didt}. Linearizing Eqs. \eqref{eq:dkdt}-\eqref{eq:dhdt} at $k=h=0$ results in the following,
 \begin{align}
     \label{eq:linearized_k_h}\mb \dot{k} \\ \dot{h} \me \approx \mb 0 & -\frac{3}{8}K(5 \cos 2\da{i} -1) \\ \frac{3}{2}K & 0  \me \mb k \\ h \me.
 \end{align}
 Linear stability for the circular equilibria depends solely on $\da{i}$, where,
\begin{align}
    \sin^2 \da{i} < 2/5 \; & : \; \text{Elliptic (stable) } \\  
    \sin^2 \da{i} > 2/5 \; & : \; \text{Hyperbolic (unstable)} \\ 
    \sin^2 \da{i} = 2/5 \; & : \; \text{Critical value (bifurcation)}.
\end{align}
The next subsection examines the noncircular equilibria that emerge from the bifurcation at the critical value of $\da{i}$. 

\subsubsection{Frozen Equilibria ($e > 0$)}

At the critical inclination ($\sin^2 i_{\mathrm{crit}} = 2/5$), the linearized system in Eq. \eqref{eq:linearized_k_h} undergoes a bifurcation that yields ``frozen'' equilibria. The eigenvector associated with the zero eigenvalue is $\begin{bmatrix} 0 & 1 \end{bmatrix}^\intercal$, aligned with the $h$-axis, indicates the direction of the emerging branch. The new equilibria satisfy
\begin{align}
    & k = 0 \\
    \label{eq:h_frozen} & h = \pm \sqrt{\frac{5\sin^2 i - 2}{3}}\quad \left(\text{for } \sin^2 i > \frac{2}{5}\right),
\end{align}
corresponding to $(e, \omega) = ( \sqrt{\frac{5\sin^2 i - 2}{3}}, \pm 90^\circ)$. These equilibria exist only in the inclination regime where circular equilibria are unstable. At each frozen equilibrium, the elements $e, i, \omega$ remain constant with $e > 0$. Linear stability analysis confirms that these frozen equilibria are elliptic within the reduced space.

\subsection{Global Solution Space in the Reduced Domain: Lidov Diagram}

\begin{figure}[h!]
    \centering
    \begin{subfigure}[b]{0.33\textwidth}
        \centering
        \includegraphics[width = 0.99\textwidth]{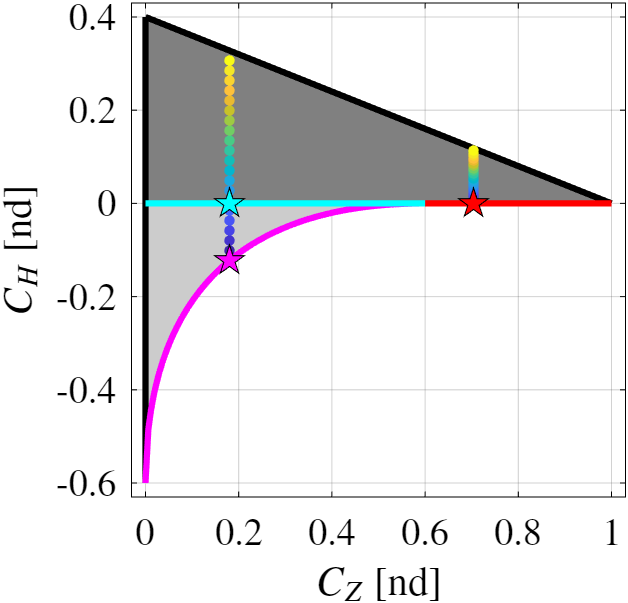}   \caption{\label{fig:2d_surface_constants_averaged_stable_unstable}Lidov diagram.}
    \end{subfigure}   
    \begin{subfigure}[b]{0.33\textwidth}
        \centering
        \includegraphics[width = 0.99\textwidth]{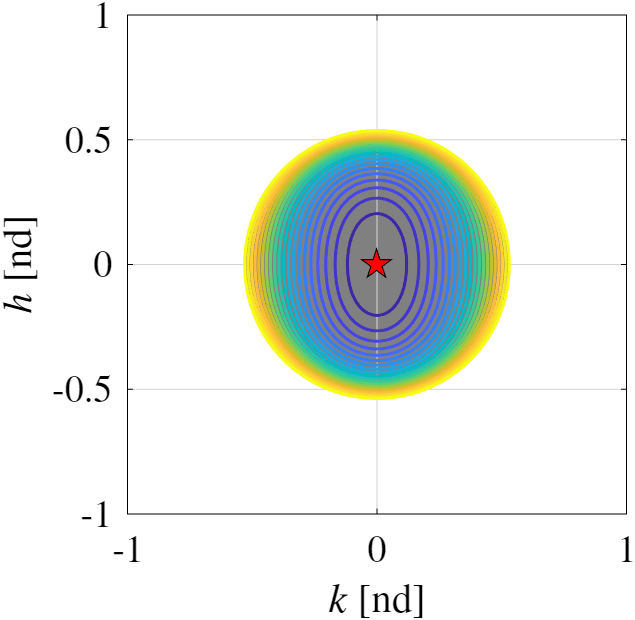}        \caption{\label{fig:periodic_orbits_in_e_w_stable}$k-h$ space ($C_Z \approx 0.705$).}
    \end{subfigure} 
    \begin{subfigure}[b]{0.33\textwidth}
        \centering
        \includegraphics[width = 0.99\textwidth]{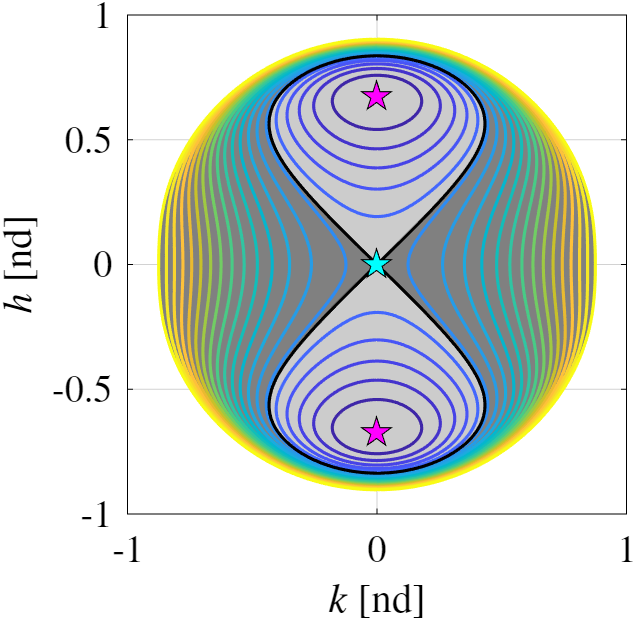}        \caption{\label{fig:periodic_orbits_in_e_w_unstable}$k-h$ space ($C_Z \approx 0.180$).}
    \end{subfigure}  
    \caption{\label{fig:dadm_solution_configuration}\acrshort{dadm} solution space within the reduced domain, $(k, h) = (\da{e} \cos{\da{\omega}}, \da{e}\sin\da{\omega})$.}
\end{figure}

The separable dynamics in $(e,i,\omega)$ facilitate a global visualization of the phase space via the ``Lidov diagram'' \cite{lidov1961evolution, broucke2003long}. The first two integrals, $C_Z$ and $C_H$, define the axes in Fig. \ref{fig:2d_surface_constants_averaged_stable_unstable}. Admissible motions are confined within the boundaries. Physical interpretations of the regions are detailed in Broucke \cite{broucke2003long}; key features relevant to this work are summarized here. Circular equilibria satisfy $C_H = 0$ with $C_Z = C_Z(i)$ for an arbitrary inclination. In the Lidov diagram (Fig. \ref{fig:2d_surface_constants_averaged_stable_unstable}), the red and cyan segments along the horizontal axis denote stable ($\sin^2 i < 2/5$) and unstable ($\sin^2 i > 2/5$) circular equilibria, respectively. Frozen equilibria emerge along the magenta curve. Two representative slices at constant $C_Z$ values illustrate the phase portraits in the $(k,h)$ plane. For a higher $C_Z \approx 0.705$ (Fig. \ref{fig:periodic_orbits_in_e_w_stable}), a single elliptic circular equilibrium (\textcolor{ccs}{$\bigstar$}) is surrounded by periodic orbits. Conversely, for a lower $C_Z \approx 0.180$ where $\sin ^2 i > 2/5$ (Fig. \ref{fig:periodic_orbits_in_e_w_unstable}), the circular equilibrium becomes hyperbolic (\textcolor{ccu}{$\bigstar$}), and two elliptic frozen equilibria (\textcolor{cf}{$\bigstar$}) appear. A homoclinic separatrix (black curve) partitions the circulating and librating motions of $\omega$. These distinct dynamical regimes are shaded dark and light grey in Fig. \ref{fig:dadm_solution_configuration}. A note on interpretation is necessary: while the unstable circular fixed point lies on $C_H = 0$ (for $C_Z < 3/5$), the converse is not unique. The level set $C_H = 0$ also contains the homoclinic orbit characterized by $\sin^2 i\sin^2 \omega = 2/5$ (see Eq. \eqref{eq:C_H}) with $e \neq 0$. Furthermore, the definition $C_Z = (1-e^2)\cos^2 i$ does not distinguish between prograde ($i < 90^\circ$) and retrograde ($180^\circ - i$) orbits. These degeneracies require careful identification of equilibria within the diagram.

\subsection{Resonance Conditions and Frequency Mapping in the Full Phase Space}

While the separable flow in the $(e, i, \omega)$ subspace clarifies qualitative behaviors, a rigorous comparison with the unaveraged dynamics requires restoring the full dimensionality of the phase space. Generally, the separation of variables exploited in the \acrshort{dadm} does not persist in the unaveraged models. The \acrshort{dadm} is characterized by three integrals of motion (Eqs. \eqref{eq:C_Z}-\eqref{eq:C_a}) and three associated phase angles. Recovering the previously omitted coordinates ($a, \Omega, M$) reintroduces the missing dimensions. For both the circular and frozen equilibria, three fundamental phase angles are defined as\footnote{Note, for non-equilibria solutions in the reduced space, $\dot{\Omega}$ is not constant from Eq. \eqref{eq:dOmegadt} and, thus, $\phi_m$ is not equivalent to $-\Omega_{\mathrm{R}}$. Rather, $\nu_m$ is defined as the time-averaged rate of change of $-\Omega_{\mathrm{R}}$, evaluated over one revolution of the long-period motion. Such distinction does not apply to the current analysis since the focus is on equilibria in the reduced space, where $\phi_m = -\Omega_{\mathrm{R}}$ holds exactly.},
\begin{align}
    \label{eq:phase_s} \phi_s &:= M \quad \text{(short-period angle, mean anomaly)} \\
    \label{eq:phase_m} \phi_m &:= t + t_0 - \Omega = -\Omega_{\mathrm{R}} \quad \text{(medium-period angle, negative RAAN in the rotating frames)} \\
    \label{eq:phase_l} \phi_l & \quad \text{(long-period angle)}.
\end{align}
The associated (fixed) frequencies, $\nu_{s,m,l} = d\phi_{s,m,l}/dt$, are functions of the integrals of motion. In contrast to $\phi_s$ and $\phi_m$, the long-period angle $\phi_l$ does not possess a direct geometrical analog in the inertial frame; instead, its frequency $\nu_l$ corresponds to the periodicity of the bounded motion within the reduced $(k, h)$ subspace (i.e., the closed loops in Figs. \ref{fig:periodic_orbits_in_e_w_stable} and \ref{fig:periodic_orbits_in_e_w_unstable}). Consequently, at an equilibrium where the reduced state is stationary, $\nu_l$ vanishes, and the phase $\phi_l$ becomes degenerate. Note that $\phi_m$ is defined using the negative of the RAAN to ensure a positive frequency $\nu_m > 0$ in the rotating frame.

Then, the dimensionality of the \acrshort{dadm} structures within the full phase space are examined. Denote an $n$-dimensional invariant torus by $\TT^n$. In the separated \acrshort{dadm}, solutions are either $\TT^0$ (equilibria) or $\TT^1$ (periodic orbits, PO); these structures correspond to the first column of Table~\ref{table:tori}. Restoring the omitted angles $(\phi_s, \phi_m)$ generically adds two degrees of quasi-periodicity, whereas integer relations among the phase rates can reduce the resulting dimension as \begin{align}
D_{\mathrm{sep}} &\in \{0,1\},\\
\boldsymbol{\nu} &= (\nu_s,\nu_m,\nu_l),\\
\mathrm{r} &= \#\ \text{of independent resonances among }(\nu_s,\nu_m,\nu_l).
\end{align}
Then the dimensional lift from the separated to the full space is supplied as,
\begin{align}
\TT^{D_{\mathrm{sep}}}
\;&\xrightarrow[\ \text{restore }(\phi_s,\phi_m)\ ]{\ \text{resonance count }\mathrm{r}\ }\;
\TT^{D_{\mathrm{full}}}, \label{eq:dim-lift-map}\qquad \,D_{\mathrm{full}} = D_{\mathrm{sep}} + 2 - \mathrm{r}\,,\qquad
0 \le \mathrm{r} \le D_{\mathrm{sep}}+1.
\end{align}
Resonances are expressed by an integer relation,
\begin{align}
\exists\,\mathbf m \in \mathbb Z^{3}\setminus\{\mathbf 0\}
\ \ \text{s.t.}\ \ 
\mathbf m \cdot \boldsymbol{\nu} &= 0.
\label{eq:resonance}
\end{align}
Note that for equilibria ($D_{\mathrm{sep}} = 0$), the resonance vector takes the form $\mathbf{m}=(m_s, m_m, 0)$ as the long-period frequency vanishes, implying $\mathrm{r} \leq 1$. Conversely, when the separated solution is periodic ($\TT^1$), up to two independent relations are possible ($\mathrm{r} \leq 2$). In the absence of resonances ($\mathrm{r}=0$), $\TT^0$ and $\TT^1$ from the reduced domain lift to $\TT^2$ (2D quasi-periodic orbit, 2D QPO) and $\TT^3$ (3D QPO), respectively. Such scenarios are included in the second column of Table~\ref{table:tori}. When a single resonance exists ($\mathrm{r} = 1$), $\TT^0$ and $\TT^1$ lift to $\TT^1$ (PO) and $\TT^2$ (2D QPO), respectively (the third column in Table~\ref{table:tori}). Finally, the presence of two independent resonances ($\mathrm{r} = 2$) only applies to $\TT^1$ in the reduced space that lifts back to $\TT^1$ (PO) in the full space (the last column in Table~\ref{table:tori}). In summary, the resonance conditions dictate the dimensionality of the lifted structures in the full phase space, with $\mathrm{r}$ serving as a critical parameter in this mapping; possible dimensions of the \acrshort{dadm} tori in the full phase space are compiled in Table \ref{table:tori}.

The current analysis focuses on $\TT^1$ (PO) that lift from $\TT^0$ with $\mathrm{r} = 1$ (highlighted in red in Table~\ref{table:tori})\footnote{This class of periodic orbits draws a contrast to the blue entry in Table~\ref{table:tori}, i.e., $\TT^1$ originating from $\TT^1$ via $\mathrm{r} = 2$; the latter typically carry longer periods and more complex geometries, and are not the primary focus of this work.}. Going forward, the terms \textbf{``circular PO''} and \textbf{``frozen PO''} denote \textcolor{red}{$\TT^1$ (PO)} that lift from the circular and frozen equilibria ($\TT^0$), respectively. For the assumed, single resonance case ($\mathrm{r} = 1$), denote the commensurability as $p\nu_m = q\nu_s$ with coprime $p,q\in\mathbb Z_{>0}$. The resonance ratio is defined as,
\begin{align}
\label{eq:ratio}\eta \;:=\; \frac{\nu_s}{\nu_m} \;=\; \frac{p}{q}.
\end{align}
Physically, the spacecraft completes $p$ cycles in $\phi_s$ (mean anomaly) while completing $q$ cycles in $\phi_m$ (negative of the RAAN within the rotating frames) with $p > q$. The short-period frequency depends solely on the semi-major axis ($\nu_s = n/n_1$), while the medium-period frequency $\nu_m$ varies with both $a$ and $i$,
\begin{align}
    \label{eq:nu_m_circular} & \text{Circular:} \quad \nu_m = 1+\frac{3}{4}(1-\mu)\frac{n_1}{n} \cos i \quad (0^\circ < i < 180^\circ)  \\
    \label{eq:nu_m_frozen}& \text{Frozen:} \quad  \nu_m = 1+ \frac{1}{4}(1-\mu) \frac{n_1}{n} \sqrt{\frac{3}{5}} (20 \sin^2 i -5) \frac{\cos i }{|\cos i|} \quad (i \neq 90^\circ).
\end{align}
Note that the leading ``1'' in $\nu_m$ in both Eqs. \eqref{eq:nu_m_circular} and \eqref{eq:nu_m_frozen} reflects the rotation rate of the rotating frames, and $n_1/n = \sqrt{a_H^3}$ via the Hill scaling. A key distinction exists: the circular locus varies smoothly through $i = 90^\circ$, whereas the frozen locus jumps due to the $\text{sgn}(\cos i)$ term.

\begin{figure}[h!]
    \centering
    \begin{subfigure}[b]{0.45\textwidth}
        \centering
        \includegraphics[width = 0.99\textwidth]{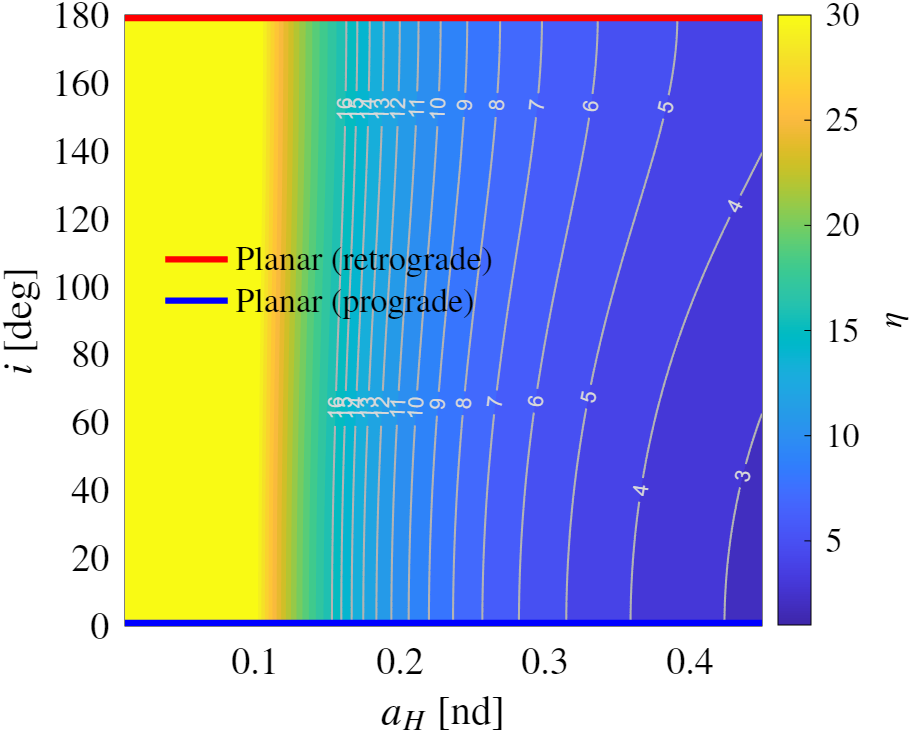}   \caption{\label{fig:a_i_eta_C_hr3bp}Circular periodic orbits ($\mu = 0$).}
    \end{subfigure}   
    \begin{subfigure}[b]{0.45\textwidth}
        \centering
        \includegraphics[width = 0.99\textwidth]{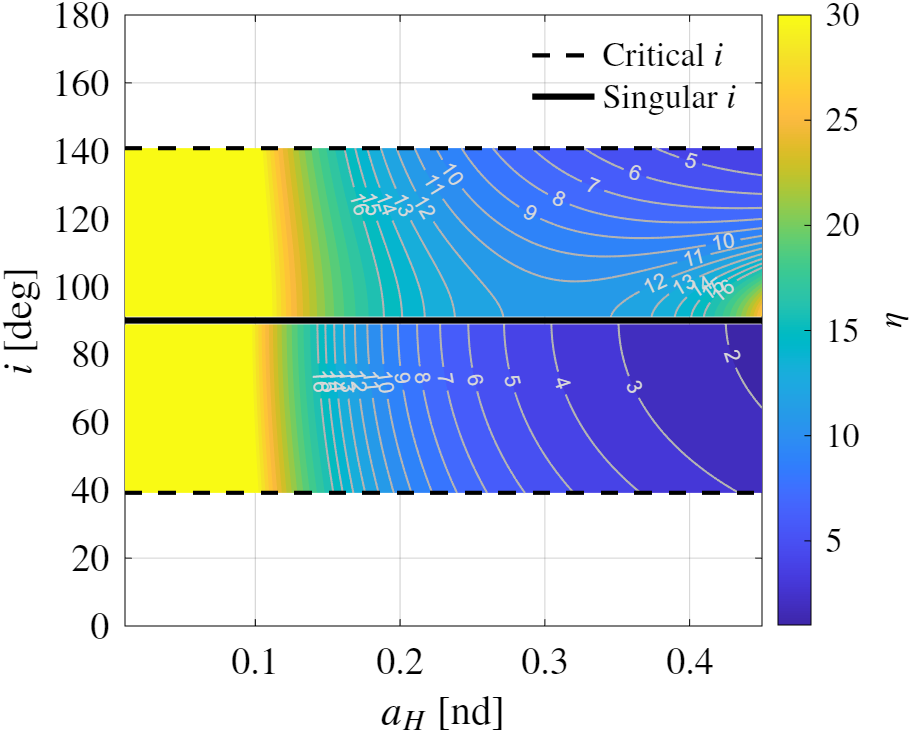}        \caption{\label{fig:a_i_eta_F_hr3bp}Frozen periodic orbits ($\mu = 0$).}
    \end{subfigure} 
    \begin{subfigure}[b]{0.45\textwidth}
        \centering
        \includegraphics[width = 0.99\textwidth]{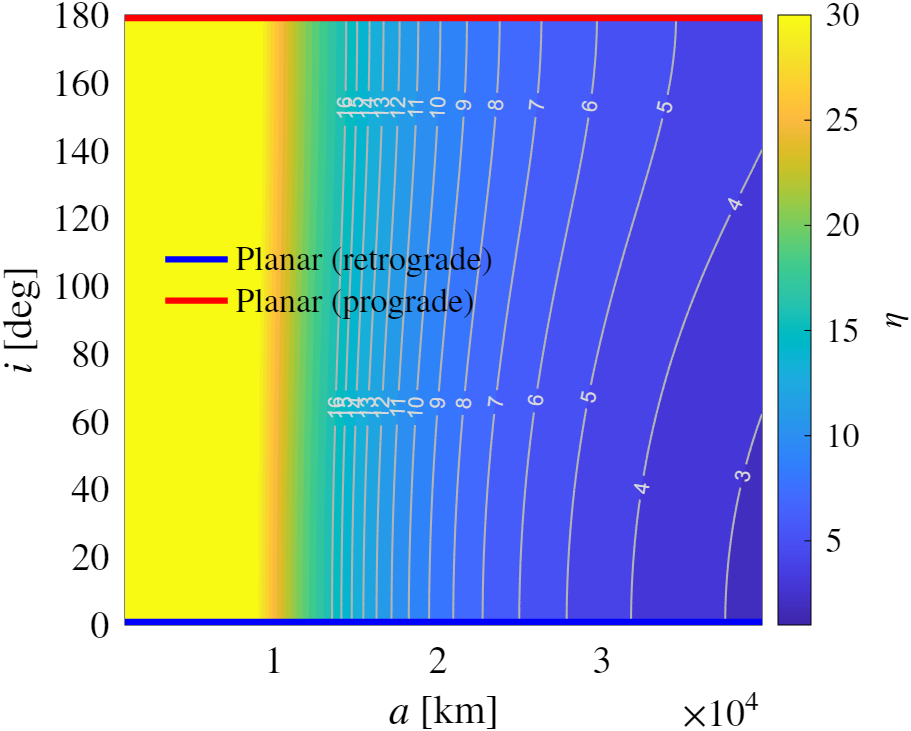}   \caption{\label{fig:a_i_eta_C_cr3bp}Circular periodic orbits ($\mu \approx 0.01215$).}
    \end{subfigure}   
    \begin{subfigure}[b]{0.45\textwidth}
        \centering
        \includegraphics[width = 0.99\textwidth]{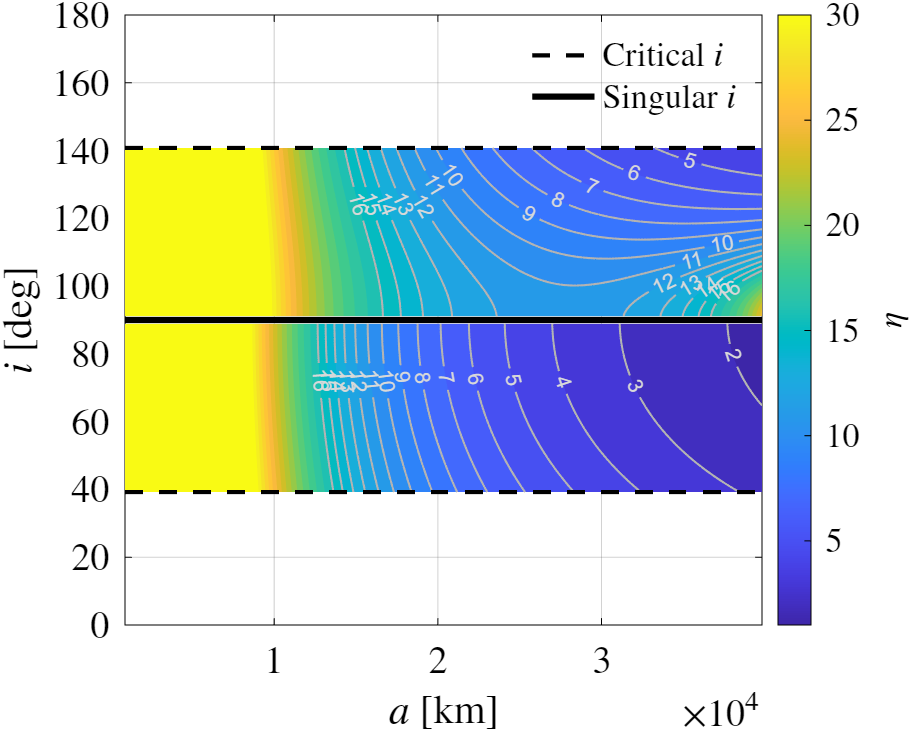}        \caption{\label{fig:a_i_eta_F_cr3bp}Frozen periodic orbits ($\mu \approx 0.01215$).}
    \end{subfigure} 
    \caption{\label{fig:a_i_eta} $\eta = \nu_s/\nu_m$ for the \acrshort{dadm} constructed for $\mu = 0$ and $\mu \approx 0.01215$ (Earth-Moon system).}
\end{figure}

The resonance ratio $\eta$ (Eq. \eqref{eq:ratio}) allows for mapping resonant bands over the parameter space \cite{lara2005dynamic}. Figure \ref{fig:a_i_eta} presents contour plots of $\eta$ within the $(a_H, i)$ domain. The top row (Figs. \ref{fig:a_i_eta_C_hr3bp}, \ref{fig:a_i_eta_F_hr3bp}) displays the limiting case $\mu=0$ (\acrshort{hr3bp} limit), while the bottom row (Figs. \ref{fig:a_i_eta_C_cr3bp}, \ref{fig:a_i_eta_F_cr3bp}) depicts the Earth-Moon system ($\mu \approx 0.01215$). The parameter range is restricted to $a_H \leq 0.45$ [nd], corresponding to a short-period frequency $\nu_s \gtrapprox 3.3$ [radians/nd]. This upper bound encompasses the validity regime of the \acrshort{dadm}; since the averaging approximation relies on a clear timescale separation ($\nu_s \gg 1$), extending beyond this limit diminishes the model's fidelity. The nearly identical contour patterns observed for $\mu = 0$ and $\mu \approx 0.01215$ confirm that these maps serve as a robust guide for locating $\TT^1$ families with $r = 1$ within the \acrshort{dadm} regardless of the mass ratio $\mu$. As indicated in Fig. \ref{fig:a_i_eta_C_hr3bp}, the circular POs evolve smoothly over $0^\circ < i < 180^\circ$. The lines at $i = 0^\circ$ and $180^\circ$ represent planar circular orbits (prograde and retrograde, respectively), further examined in subsequent sections. In contrast, the frozen POs exist only where $\sin^2 i > 2/5$ ($i_\mathrm{crit}^-<i<i_\mathrm{crit}^+$ with $i_\mathrm{crit}^- \approx 39.2^\circ$ and $i_\mathrm{crit}^+ \approx 140.8^\circ$) and exhibit a discontinuity at $i= 90^\circ$ due to the singularity in Eq. \eqref{eq:nu_m_frozen}. Finally, while prograde frozen POs ($i < 90^\circ$) generally display monotonic behavior with respect to $a$, retrograde frozen POs ($i > 90^\circ$) reveal more intricate structural variations.

\begin{table}[htbp]
\centering
\caption{\label{table:tori}Dimensional lift from the reduced (separated) DADM phase space to the DADM full phase space. Here, $\mathrm{r}$ denotes the number of independent resonance relations among $(\nu_{s,m,l})$.}
\label{tab:dim-lift}
\begin{tabular}{@{}lccc@{}}
\toprule
\multirow{2}{*}{\textbf{DADM tori (separated)}} 
  & \multicolumn{3}{c}{\textbf{DADM tori (full phase space): resonance count} $\mathrm{r}$}  \\
\cmidrule(l){2-4}
& \boldmath$0$ & \boldmath$1$ & \boldmath$2$ \\
\midrule
\textcolor{red}{$\TT^0$ (Equilibria)}       & $\TT^2$ (2D QPO)  & \textcolor{red}{$\TT^1$ (PO)} & N/A \\
$\TT^1$ (PO) & $\TT^3$ (3D QPO)  & $\TT^2$ (2D QPO)   & \textcolor{blue}{$\TT^1$ (PO)} \\
\bottomrule
\end{tabular}
\end{table}

The tori from the \acrshort{dadm} (full phase space) may persist as tori with the same dimension in the unaveraged dynamics, a case that remains central to the current work. Thus, the current analysis focuses on the evolution of \textcolor{red}{$\TT^0$ (Equilibria)} from the separated \acrshort{dadm} to analog \textcolor{red}{$\TT^1$ (PO)} within the unaveraged dynamics (\acrshort{hr3bp}, \acrshort{cr3bp}) with $\mathrm{r} = 1$. As the exact frequency values may shift across dynamical models, the integer resonance relation $\eta$ itself serves as the robust matching criterion. It is noted that torus destruction is expected outside the regime of validity, for example, at large $a$ where the perturbation from $\po$ is significantly large. Such breakdowns remain out-of-scope for the current work. Rather, the analysis focuses on the ``nominal'' scenarios where circular and frozen periodic orbits are well-defined within the unaveraged dynamics.

\section{\label{sec:r3bp}Symmetric Periodic Orbits within Unaveraged Dynamics Lifted from Averaged Equilibria}

While the frequency analysis in Section \ref{sec:dadm} identifies the orbital element combinations that theoretically support resonant periodic orbits ($\TT^1$) within the \acrshort{dadm}, key discrepancies remain in realizing these solutions within the unaveraged dynamics. These gaps include: (1) the \emph{kinematic transformation} from the orbital elements in the inertial frame to the states within the rotating frames, and (2) the \emph{dynamical disparity} between the integrable averaged model and the non-integrable \acrshort{r3bps}. To bridge these gaps, \emph{symmetry} serves as a critical criterion; configurations that satisfy specific symmetric properties are known to persist robustly within the \acrshort{r3bps} as symmetric periodic orbits, shaping the backbone of the solution space. Accordingly, this section investigates the geometric correspondence between the averaged and unaveraged models through the lens of symmetry. The central inquiry is twofold: (1) \emph{``When do the \acrshort{dadm} equilibria satisfy the symmetric configuration conditions within the rotating frames}, and (2) \emph{what is the resulting multiplicity of the solution geometries?''} The answer depends on a combinatorial interplay of the equilibrium type (circular vs. frozen), the parity of the resonance integers $p:q$, and the specific symmetry group of the target unaveraged model. This section derives a systematic initialization strategy based on modular arithmetic to locate these symmetric configurations and determine their count.

The analysis proceeds in three stages. First, Section \ref{sec:cr3bp_symmetry} reviews the intrinsic symmetries of the \acrshort{hr3bp} and \acrshort{cr3bp} and classifies the requisite symmetric apse conditions. Second, Section \ref{sec:cr3bp_init} analyzes the phase evolution of the \acrshort{dadm} equilibria under the kinematic transformation. By leveraging modular arithmetic, this section identifies the admissible symmetric apse configurations for a given resonance parity defined by the integers $p:q$. Finally, Section \ref{sec:initialization_targeting} presents the initialization and targeting results, systematically transforming the \acrshort{dadm} equilibria to the respective periodic orbits within the \acrshort{r3bps}.

\subsection{\label{sec:cr3bp_symmetry}Symmetries within the Unaveraged Dynamics and Symmetric Apse Conditions}

The existence and classification of symmetric periodic orbits are governed by the invariance of the equations of motion under discrete coordinate transformations. The \acrshort{hr3bp}, owing to its idealized assumptions regarding the distant perturbing body ($\mathcal{P}_1$), admits a symmetry group of order 8 \cite{aydin2023linear}, generated by three independent symmetries. An intuitive set of generators includes the time-reversing reflections across two vertical planes ($x_H = 0, y_H = 0$) and the geometric reflection across the horizontal plane ($z_H = 0$). In contrast, the \acrshort{cr3bp} represents a symmetry-broken environment where the reflection across the $x_H = 0$ plane (corresponding to $x - (1-\mu)=0$ in the \acrshort{brf}) is lost. This symmetry breaking reduces the symmetry group to order 4, generated by only two independent symmetries. Consequently, periodic orbits within the \acrshort{hr3bp} may exhibit up to three independent symmetries, whereas those in the \acrshort{cr3bp} are limited to a maximum of two.

This analysis focuses on the time-reversing symmetries ($\rho$). The operator $\rho$ entails a time reversal, $t \rightarrow -t$. The position and velocity transform according to $\bm{r}\rightarrow \bm{S}\bm{r}$ and $\bm{v} \rightarrow -\bm{S} \bm{v}$, where $\bm{S} = \text{diag}(S_x, S_y, S_z)$ is a spatial reflection matrix with components $\pm 1$. The admissible time-reversing symmetries in the \acrshort{hr3bp} and \acrshort{cr3bp} are summarized in Table \ref{tab:symmetry_comparison}. Here, the symmetries OX and OY are termed \textit{axial}, corresponding to $180^\circ$ rotations about the $\hat{\bm{x}}$ and $\hat{\bm{y}}$ axes, respectively. The symmetries XOZ and YOZ are termed \textit{reflectional}, involving reflections across the $xz$- and $yz$-planes. The fixed sets $\Sigma$ of these transformations, i.e., states that remain invariant under $\rho$, define the mirror configurations facilitating the detection and construction of symmetric periodic orbits. Geometrically, these fixed sets correspond to orthogonal apse crossings, satisfied when $\bm{r}\cdot \bm{v} = 0$ and, thus, are also termed ``symmetric apse sets.'' The fixed set for each symmetry is listed in the last column of Table \ref{tab:symmetry_comparison}. Note that while these conditions are expressed in the \acrshort{hrf}, the same constraints may be equivalently represented within the \acrshort{brf} after translation of the origin and re-scaling.

\begin{table}[htbp]
\centering
\caption{\label{tab:symmetry_comparison}Time-reversing symmetries ($\rho$) in the R3BPs and their fixed sets (symmetric apse sets).}
\begin{tabular}{@{}llcc@{}}
\toprule
\textbf{Symmetry Label} & \textbf{Transformation} $(S_x, S_y, S_z)$  & \textbf{Valid Models} & \textbf{Fixed Set} ($\Sigma$) \\
\midrule
\textbf{OX} & \quad $\hat{x}$-axis rotation $(+1, -1, -1)$ & \acrshort{hr3bp} \& \acrshort{cr3bp} & $y_H=z_H=\dot{x}_H=0$ \\
\textbf{XOZ} & \quad $xz$-plane reflection $(+1, -1, +1)$ & \acrshort{hr3bp} \& \acrshort{cr3bp} & $y_H=\dot{x}_H=\dot{z}_H=0$ \\
\textbf{OY} & \quad $\hat{y}$-axis rotation $(-1, +1, -1)$ & \acrshort{hr3bp} only & $x_H=z_H=\dot{y}_H=0$ \\
\textbf{YOZ} & \quad $yz$-plane reflection $(-1, +1, +1)$ & \acrshort{hr3bp} only & $x_H=\dot{y}_H=\dot{z}_H=0$ \\
\bottomrule
\end{tabular}
\end{table}

The structural classification of periodic orbits within the \acrshort{r3bps} is facilitated based on the multiplicity of the time-reversing symmetries they preserve; to maintain a time-reversing symmetry, an orbit must intersect the respective fixed set at two distinct locations within one period. Let $N$ denote the number of distinct types of symmetric apse sets ($\Sigma$) intersected by a spatial ($z_H \neq 0$) periodic orbit. The possible configurations are classified as follows:
\begin{itemize}[topsep=0pt, partopsep=0pt]
    \item \textbf{Singly Symmetric ($N=1$):} The orbit intersects a single fixed set $\Sigma$ twice per period, separated by $t=P/2$, without encountering a distinct fixed set at the quarter period intervals.
    \item \textbf{Doubly Symmetric ($N=2$):} The orbit intersects two distinct fixed sets (e.g., $\Sigma_{\text{OX}}$ and $\Sigma_{\text{XOZ}}$) in an alternating sequence every $t=P/4$, rendering two independent symmetries.
\end{itemize}
With this classification, the task of locating spatial symmetric periodic orbits reduces to evaluating the apse conditions listed in Table \ref{tab:symmetry_comparison} at half and quarter period intervals.

The above classification requires refinement to address the following nuances. First, in the generic spatial domain, triply symmetric orbits do not exist. The imposition of two independent apse conditions ($N=2$) already fixes the quarter-period alternation and ensures closure without encountering a third, distinct apse. Second, geometric symmetries without time reversal also exist in the \acrshort{r3bps}. A notable example is the reflection across the $z=0$ plane (denoted as $\sigma$) that exists in both the \acrshort{hr3bp} and \acrshort{cr3bp}. While the composition of certain two time-reversing symmetries yields this geometric symmetry (i.e., $\rho_{\text{OX}} \circ \rho_{\text{XOZ}} = \rho_{\text{OY}} \circ \rho_{\text{YOZ}} = \sigma$) \cite{aydin2023linear}, it is possible for $\sigma$ to exist independently without any time-reversing symmetries (e.g., vertical Lyapunov orbits for the triangular Lagrange points within the \acrshort{cr3bp}). While such orbits are technically ``singly symmetric'' (possessing $\sigma$), they remain out-of-scope in the current investigation\footnote{These $\sigma$-symmetric structures likely bifurcate from planar periodic orbits that do not encounter the time-reversing symmetric apse sets defined here. As the planar periodic orbits under investigation (see Section \ref{sec:bifurcation}) intersect the symmetric apse sets at least once, the occurrence of purely $\sigma$-symmetric orbits is expected to be rare.}. Thus, for spatial periodic orbits in this work, the terms singly and doubly symmetric strictly refer to $N = 1$ and $N = 2$. Lastly, a degeneracy exists for planar periodic orbits where the $\sigma$-symmetry is trivially satisfied. Consequently, planar periodic orbits are at least singly symmetric ($\sigma$). Encountering the time-reversing apse sets increases the symmetry count from this baseline. Note that in the planar case, the axial and reflectional symmetries become identical (i.e., $\rho_{\text{X}} := \rho_{\text{OX}} \equiv \rho_{\text{XOZ}}$ and $\rho_{\text{Y}} := \rho_{\text{OY}} \equiv \rho_{\text{YOZ}}$). Thus, planar periodic orbits satisfy $\sigma$ (singly), $\sigma$ plus $\rho_{X}$ or $\rho_{Y}$ (doubly), or $\sigma$ plus both $\rho_X$ and $\rho_Y$ (triply). 

\begin{figure}[htbp]
    \centering
    \begin{subfigure}[b]{0.44\textwidth}
        \centering
        \includegraphics[width = 0.99\textwidth]{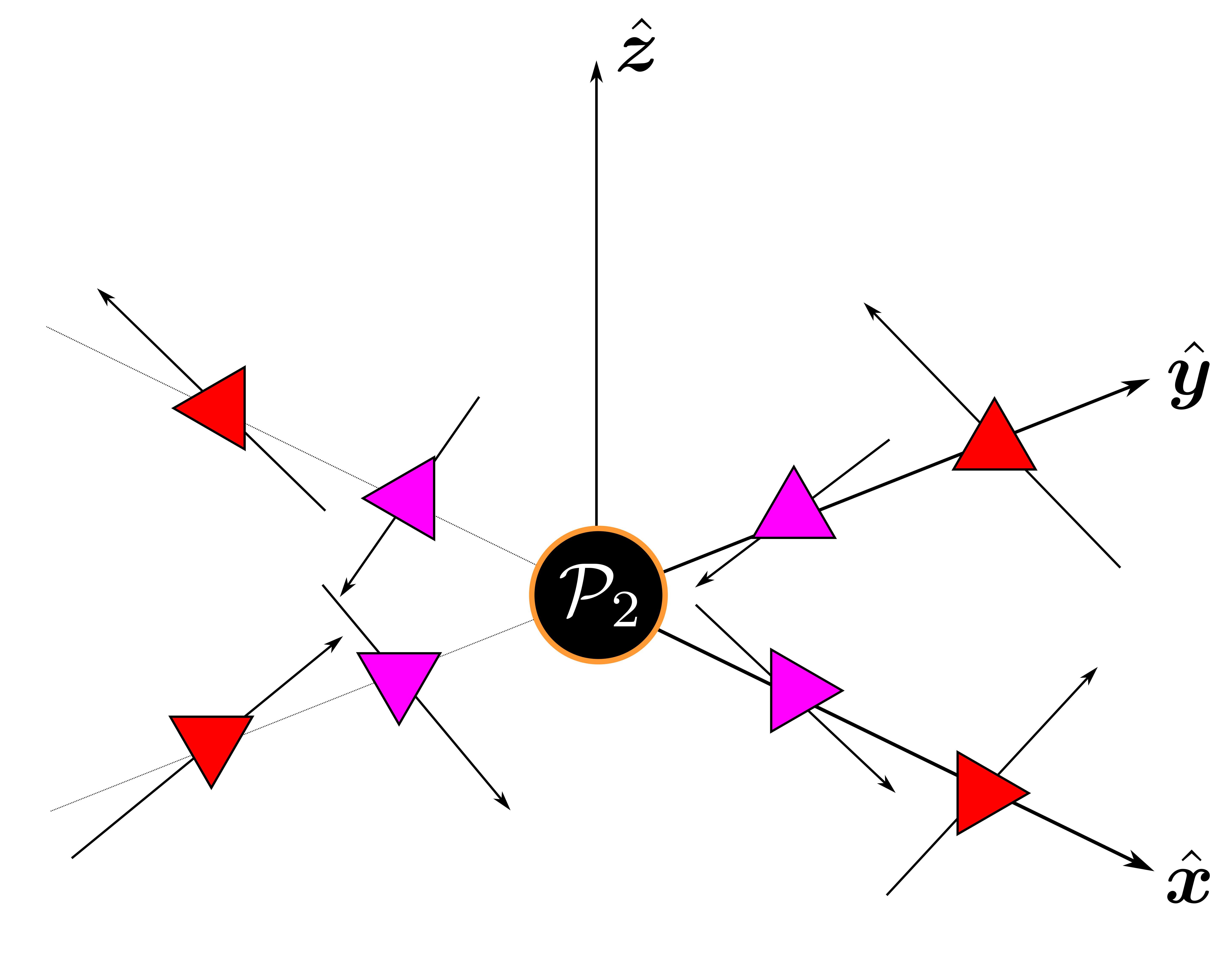}
        \caption{\label{fig:axial_illustration}Axial type.}
    \end{subfigure}
    \begin{subfigure}[b]{0.44\textwidth}
        \centering
        \includegraphics[width = 0.99\textwidth]{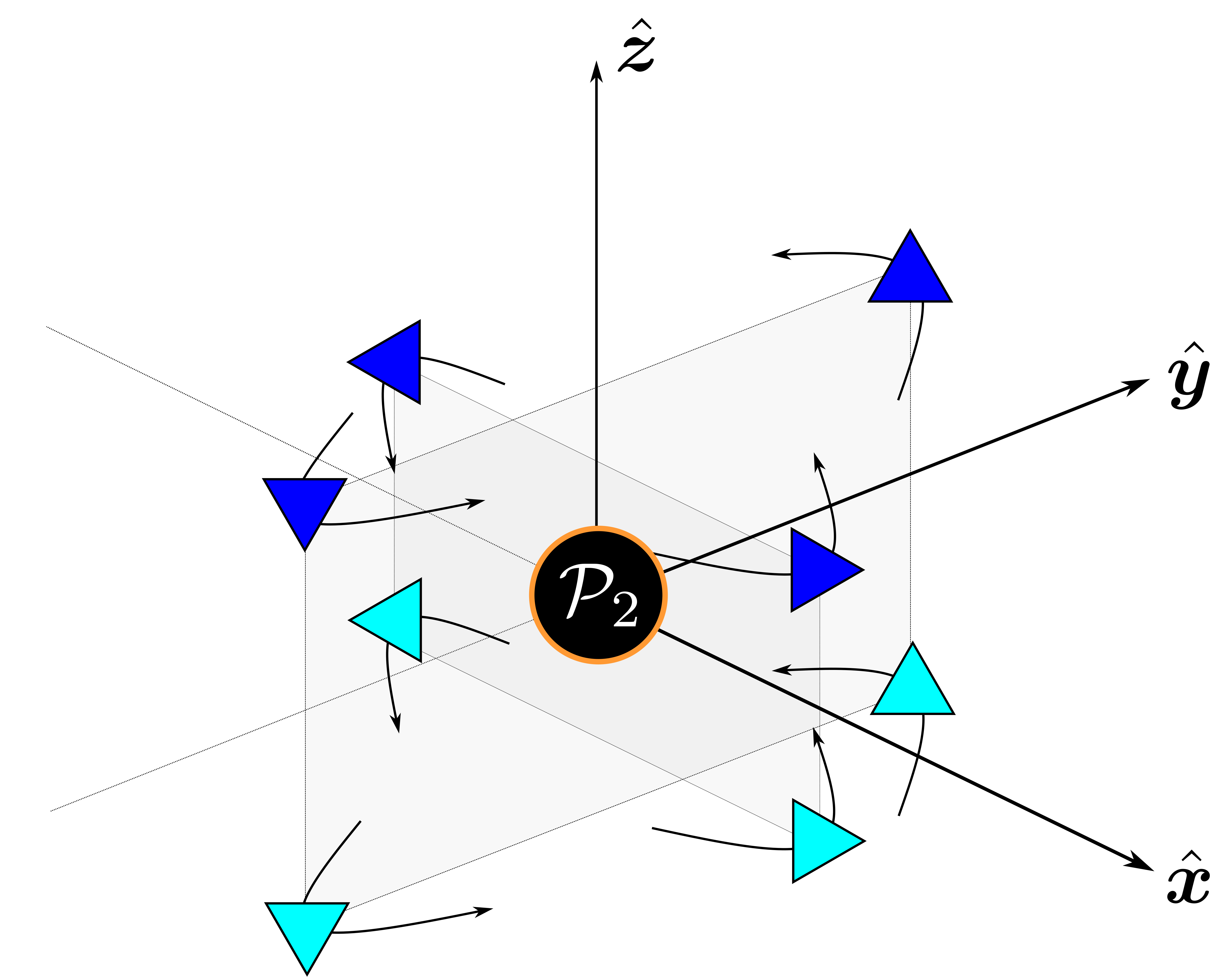}
        \caption{\label{fig:reflectional_illustration}Reflectional type.}
    \end{subfigure}
    \caption{\label{fig:node}Symmetric apse configurations within the \acrshort{hrf} (prograde, $h_z > 0$) with markers defined in Table \ref{tab:symmetry_classification}.}
\end{figure}

To translate the fixed sets into full, six-dimensional states, the geometric ambiguities inherent in the sets must be resolved. Specifically, the symmetric apse conditions in Table \ref{tab:symmetry_comparison} constrain only three state components to zero, leaving the signs of the remaining variables undefined. Systematically addressing these sign ambiguities yields multiple geometric configurations as follows: 
\begin{itemize}[topsep=0pt, partopsep=0pt]
    \item \textbf{Vertical Component Sign:} The axial and reflectional fixed sets enforce $z_H = Z = 0$ and $ \dot{z}_H = \dot{Z} = 0$, respectively. For the axial configuration, $\dot{Z} > 0$ and $\dot{Z} < 0$ result in \emph{ascending} ($\mathrm{A}$) and \emph{descending} ($\mathrm{D}$) nodes. For the reflectional configuration, $Z > 0$ and $ Z< 0$ result in the apex ($\dot{Z} = 0$) at the \emph{northern} ($\mathrm{N}$) and \emph{southern} ($\mathrm{S}$) sides. 
    \item \textbf{Location of the Apse:} While one coordinate is fixed at zero ($y_H=0$ for OX/XOZ; $x_H=0$ for OY/YOZ), the sign of the orthogonal coordinate remains ambiguous. The phase of the apse is defined as $\chi_H = \arctan_2(y_H, x_H)$, locating the apse at the $+\hat{\bm{x}}$ ($0^\circ$), $+\hat{\bm{y}}$ ($90^\circ$), $-\hat{\bm{x}}$ ($180^\circ$), or $-\hat{\bm{y}}$ ($270^\circ$) axes within the \acrshort{hrf}. 
    \item \textbf{Direction of the Orbit:} The remaining velocity component determines the sign of the angular momentum within the \acrshort{hrf} as $h_z = \hat{\bm{z}}\cdot (\bm{r}_H \times \bm{v}_H)$. The orbit is classified as \emph{prograde} ($h_z > 0$) or \emph{retrograde} ($h_z < 0$). While both configurations are admissible, the analysis mainly focuses on prograde orbits for a consistent illustration. 
\end{itemize}
Based on this classification, Table \ref{tab:symmetry_classification} summarizes the possible apse configurations with the sign ambiguities resolved. Tables \ref{tab:init_node} and \ref{tab:init_apex} focus on axial and reflectional symmetries, respectively. The first column of the tables contains the identifiers (ID) for the symmetric apsides; the notation follows the format of $\mathrm{vert}_\mathrm{axis}$ where $\mathrm{vert} = \mathrm{A, N, D, S}$, and $\mathrm{axis}=\pm x, \pm y$. The direction of the orbit $h_z$ is omitted in this identifier for a brief notation, assuming a prograde direction. To facilitate visual identification, a graphical marker system is adopted in the table as well. The marker color distinguishes the vertical geometry: ``warm'' tones (red/magenta) denote axial symmetries (nodal crossings, $z_H=0$), whereas ``cool'' tones (blue/cyan) denote reflectional symmetries (apical crossings, $\dot{z}_H=0$). Furthermore, the marker shape indicates the azimuthal phase ($\chi_H$) of the apse, with the triangle pointing toward the physical location of the apse in the rotating frame (e.g., $\xp{white}$ for $+\hat{\bm{x}}$ and $\yp{white}$ for $+\hat{\bm{y}}$). Table \ref{tab:symmetry_classification} illustrates the combinations between the vertical component sign and the location of the apse, resulting in four distinct cases for each time-reversing fixed set, yielding 16 distinct configurations. These symmetric apse configurations are visualized in Fig. \ref{fig:node}, where Figs. \ref{fig:axial_illustration} and \ref{fig:reflectional_illustration} illustrate the axial and reflectional cases, respectively, following the markers defined in Table \ref{tab:symmetry_classification}.

\begin{table*}[htbp]
\centering
\caption{\label{tab:symmetry_classification}Classification of symmetric apse configurations, visualized in Fig. \ref{fig:node}.}
    \begin{subtable}[t]{0.48\textwidth}
        \centering
        \caption{Axial type ($Z=0$).}
        \label{tab:init_node}
        \begin{tabular}{@{}cccc@{}}
            \toprule
            \textbf{ID} & \textbf{Vert. Velocity} & \textbf{Apse Location} & \textbf{Marker} \\
             & (sgn$(\dot{Z})$) & ($\chi_H$) &  \\
            \midrule
            
            \multicolumn{4}{c}{OX-Symmetric Apse} \\ 
            \midrule
            \ID{A}{+x} & Ascending ($+$)  & $+\hat{\bm{x}}$ ($0^\circ$)   & \xp{ascending} \\
            \ID{D}{+x} & Descending ($-$) & $+\hat{\bm{x}}$ ($0^\circ$)   & \xp{descending} \\
            \ID{A}{-x} & Ascending ($+$)  & $-\hat{\bm{x}}$ ($180^\circ$) & \xm{ascending} \\
            \ID{D}{-x} & Descending ($-$) & $-\hat{\bm{x}}$ ($180^\circ$) & \xm{descending} \\
            \midrule
            
            \multicolumn{4}{c}{OY-Symmetric Apse} \\ 
            \midrule
            \ID{A}{+y} & Ascending ($+$)  & $+\hat{\bm{y}}$ ($90^\circ$)  & \yp{ascending} \\
            \ID{D}{+y} & Descending ($-$) & $+\hat{\bm{y}}$ ($90^\circ$)  & \yp{descending} \\
            \ID{A}{-y} & Ascending ($+$)  & $-\hat{\bm{y}}$ ($270^\circ$) & \ym{ascending} \\
            \ID{D}{-y} & Descending ($-$) & $-\hat{\bm{y}}$ ($270^\circ$) & \ym{descending} \\
            \bottomrule
        \end{tabular}
    \end{subtable}
    \hfill
    \begin{subtable}[t]{0.48\textwidth}
        \centering
        \caption{Reflectional type ($\dot{Z}=0$).}
        \label{tab:init_apex}
        \begin{tabular}{@{}cccc@{}}
            \toprule
            \textbf{ID} & \textbf{Vert. Position} & \textbf{Apse Location} & \textbf{Marker} \\
             & (sgn$(Z)$) & ($\chi_H$) &  \\
            \midrule
            
            \multicolumn{4}{c}{XOZ-Symmetric Apse} \\ 
            \midrule
            \ID{N}{+x} & North ($+$) & $+\hat{\bm{x}}$ ($0^\circ$)   & \xp{northern} \\
            \ID{S}{+x} & South ($-$) & $+\hat{\bm{x}}$ ($0^\circ$)   & \xp{southern} \\
            \ID{N}{-x} & North ($+$) & $-\hat{\bm{x}}$ ($180^\circ$) & \xm{northern} \\
            \ID{S}{-x} & South ($-$) & $-\hat{\bm{x}}$ ($180^\circ$) & \xm{southern} \\
            \midrule
            
            \multicolumn{4}{c}{YOZ-Symmetric Apse} \\ 
            \midrule
            \ID{N}{+y} & North ($+$) & $+\hat{\bm{y}}$ ($90^\circ$)  & \yp{northern} \\
            \ID{S}{+y} & South ($-$) & $+\hat{\bm{y}}$ ($90^\circ$)  & \yp{southern} \\
            \ID{N}{-y} & North ($+$) & $-\hat{\bm{y}}$ ($270^\circ$) & \ym{northern} \\
            \ID{S}{-y} & South ($-$) & $-\hat{\bm{y}}$ ($270^\circ$) & \ym{southern} \\
            \bottomrule
        \end{tabular}
    \end{subtable}
\end{table*}

\subsection{\label{sec:cr3bp_init}Symmetric Configurations from the Averaged Equilibria}

This subsection investigates the existence and multiplicity of symmetric configurations derived from the \acrshort{dadm} equilibria following the frame transformation. Note that the underlying dynamics remain within the context of the \acrshort{dadm}, focusing specifically on the kinematic transformations. Recall that, for a single resonance $\mathrm{r}=1$ with a ratio $\eta=p/q$ (Eq.~\eqref{eq:ratio}), the corresponding equilibria in the \acrshort{dadm} lift to periodic ($\mathbb{T}^1$) solutions in the full phase space with a definitive period determined by the resonance, $P = p(2\pi/\nu_s)$. With this known period, the symmetric apse conditions (Table \ref{tab:symmetry_comparison}) may be evaluated at half and quarter period intervals to identify the conditions for the singly and doubly symmetric orbits. Furthermore, enumerating the distinct initial configurations that re-encounter a symmetric apse set at the half- and quarter-period marks supplies the basis for determining the multiplicity of solution geometries associated with each resonance. 

\begin{figure}
    \centering
    \begin{subfigure}[b]{0.48\textwidth}
        \centering
        \includegraphics[width = 0.90\textwidth]{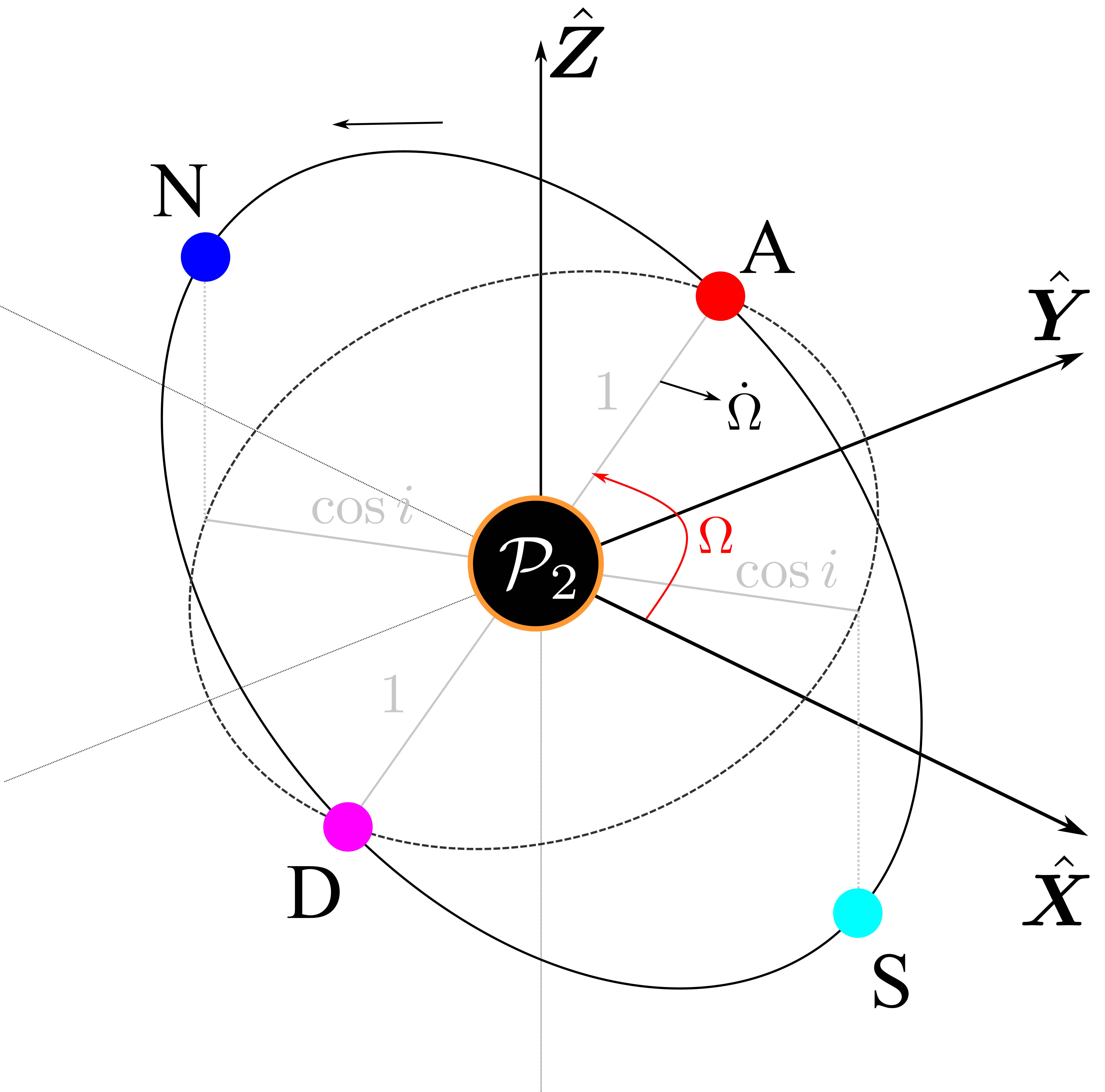}
        \caption{\label{fig:circular_planar}Osculating, circular PO ($h_z > 0$).}
    \end{subfigure}
    \begin{subfigure}[b]{0.48\textwidth}
        \centering
        \includegraphics[width = 0.90\textwidth]{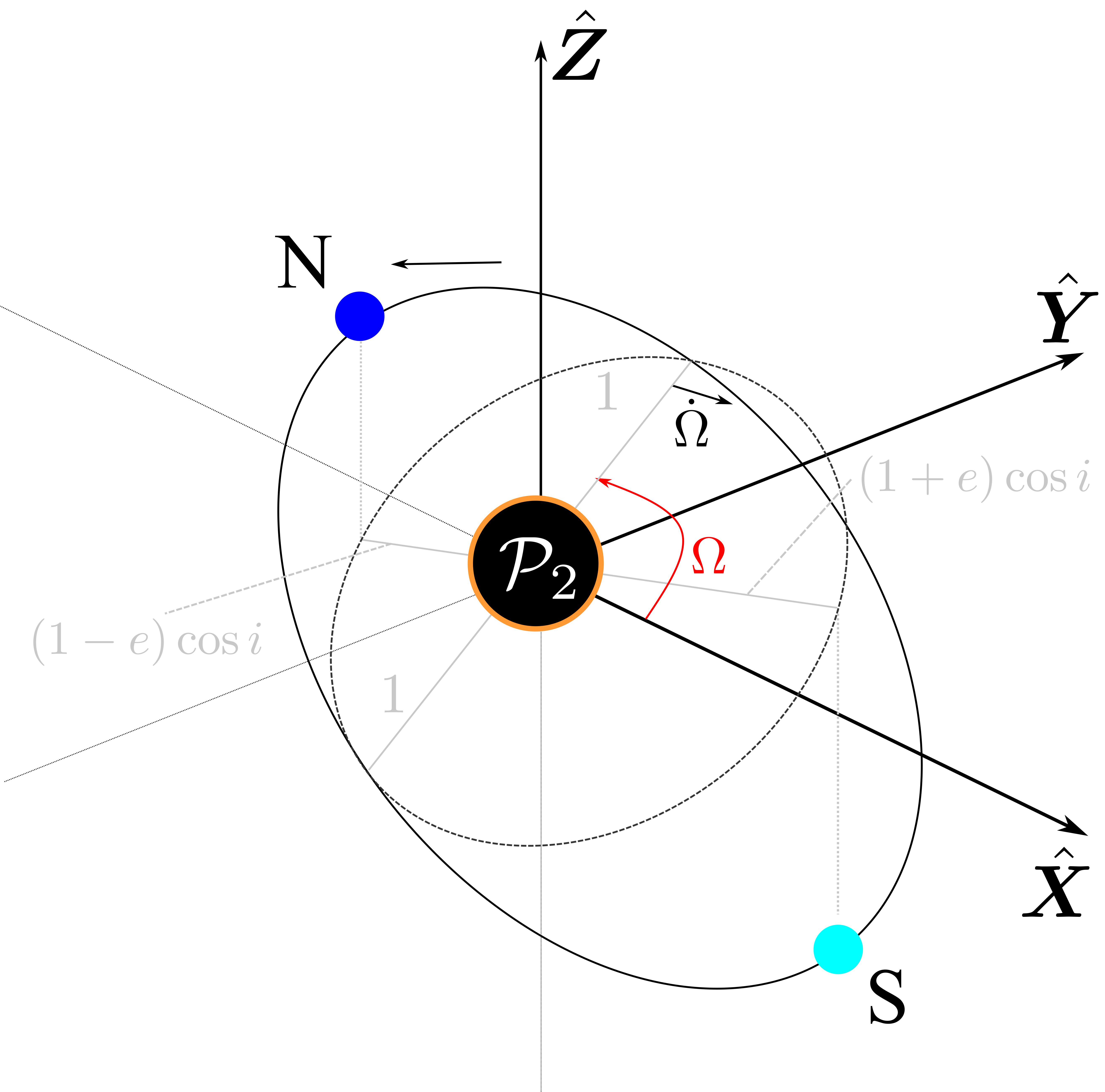}
        \caption{\label{fig:frozen_planar}Osculating, frozen PO ($h_z > 0$, $\omega = 90^\circ$).}
    \end{subfigure}
    \caption{\label{fig:planar_projection}Potential symmetric apse configurations within the \acrshort{eof}. Abbreviations A, D, N, S denote ascending node, descending node, northern apex, and southern apex, respectively. The color scheme follows the same convention as Table \ref{tab:symmetry_classification}.}
\end{figure}

To facilitate this enumeration, a frame-adaptive interpretation of the symmetric apse condition from Table \ref{tab:symmetry_comparison} is first established. Note that as $\hat{\bm{z}} = \hat{\bm{Z}}$, the vertical directional information is preserved. Thus, $Z =0 $ and $\dot{Z} = 0$ is a requirement for supplying the axial and reflectional apse conditions, respectively. Moreover, the symmetric \emph{apse} conditions satisfy $\bm{r}_H\cdot \bm{v}_H = 0$. This geometric property is frame-invariant, i.e., equivalent to $\bm{R}\cdot \bm{V} = 0$ within the \acrshort{eof} as well. This information is graphically analyzed within Fig. \ref{fig:planar_projection}. Figures \ref{fig:circular_planar} and \ref{fig:frozen_planar} depict \emph{osculating} Keplerian orbits for the circular and frozen equilibria solutions within the \acrshort{eof}. The inclined orbit projects onto the $XY$-plane as an ellipse, placing the nodes ($\mathrm{A}, \mathrm{D}$) at in-plane radius unity and the apices ($\mathrm{N}, \mathrm{S}$) at radius $\cos i$ (circular) or $(1\pm e)\cos i$ (frozen); this scaling underlies the scalars $\alpha, \beta$ in the derivations below. From the equilibria conditions, $a, e, i, \omega$ stay constant along the orbits; however, the RAAN ($\Omega$) shifts at a constant rate. For the prograde configurations, $\dot{\Omega} < 0$ from Eq. \eqref{eq:dOmegadt}. Within the figures, the colored dots denote the apsides that \emph{potentially} evolve to symmetric apse conditions within the rotating frames. For the circular equilibria, every location is an apse since $e = 0$, and four locations exist at the nodes ($\mathrm{A}, \mathrm{D}$) and apices ($\mathrm{N}, \mathrm{S}$) that may evolve to a symmetric apse after frame transformation. In contrast, the frozen equilibria only renders apsides at $M = 0^\circ, 180^\circ$; as such, only two locations exist at the apices ($\mathrm{N}, \mathrm{S}$). The evolution of these apsides within the \acrshort{eof} into symmetric apse conditions in the rotating frames is dictated by $t + t_0$ from the rotation matrix in Eq. \eqref{eq:eof_brf} as well as $\Omega$. To test the symmetric configurations, the process proceeds as follows: starting from an initial phase $t_0$ and RAAN $\Omega$ that align an inertial apse in the \acrshort{eof} with a symmetric apse condition in the rotating frame (\acrshort{hrf}), evaluate the state at half and quarter period intervals ($t = P/2, P/4$). The derivation initially considers all symmetric configurations within the \acrshort{hr3bp}; subsequently, configurations incompatible with the \acrshort{cr3bp} are pruned.

\subsubsection{Circular Equilibria to Symmetric Configurations}

This subsection identifies the admissible symmetric apse configurations and their multiplicity for the circular equilibria. Recall from Section~\ref{sec:cr3bp_symmetry} that, for an equilibrium to evolve into a symmetric PO under the kinematic transformation, the lifted orbit must intersect the fixed sets (symmetric apse sets) at two distinct points per period: the same set at $t = P/2$ (singly symmetric) or two distinct sets alternating at $t = P/4$ (doubly symmetric). Accordingly, starting from an arbitrary apse placed on a fixed set, the state is examined at the half- and quarter-period marks within the rotating frame, and the resulting apse encounters govern the number of symmetries. The geometric evolution of circular equilibria within the \acrshort{hrf} yields the following behaviors:
\begin{itemize}[topsep=0pt, partopsep=0pt]
    \item \textbf{Half period examination:} An initial state on any of the symmetric apse sets ($\Sigma$, Table \ref{tab:symmetry_comparison}) returns to the same set at a distinct location after $t = P/2$.
    \item \textbf{Quarter period examination:} An initial state on any of the symmetric apse sets ($\Sigma$, Table \ref{tab:symmetry_comparison}) intersects a different symmetric apse set after $t = P/4$.
\end{itemize}
These behaviors imply that each circular equilibrium lifts to four doubly symmetric PO geometries in the \acrshort{hr3bp} (two per symmetry pairing) that the \acrshort{cr3bp} subsequently prunes; the remainder of this subsection derives and enumerates this count via modular arithmetic. The derivations follow from a straightforward geometric reasoning. Note the rotation between the in-plane position components in the \acrshort{eof} and \acrshort{hrf} as, 
\begin{align}
    \label{eq:half_proof_1}\mb x_H \\ y_H \me & = \bm{C}_{(-t-t_0)} \mb X \\ Y \me, \quad \bm{C}_{\gamma}: =  \mb \cos \gamma & -\sin \gamma \\ \sin \gamma & \cos \gamma \me.
\end{align}
By assumption, at $t = 0$, $x_H$ or $y_H$ is zero, satisfying a symmetric apse condition. Noting that the period is predetermined from the resonance ratio (Eq. \eqref{eq:ratio}) as $P = \frac{2p\pi}{\nu_s}$. After a half period, the in-plane position components within the \acrshort{eof} evolve as,
\begin{align}
    \label{eq:half_proof_2}\mb X \\ Y \me_{t = P/2} & = \bm{C}_{(\Delta \Omega \pm p\pi)} \mb X \\ Y \me_{t = 0},
\end{align} 
following a graphical illustration supplied in Fig. \ref{fig:circular_planar}. Then, from the resonance ratio, $\Delta M = \nu_s \frac{P}{2} = p\pi$; for prograde and retrograde orbits, the in-plane components rotate by $\pm p \pi$, respectively. From this property, nodes/apices remain as nodes/apices after a half period. Also, $\Delta \Omega$ accumulates from Eq. \eqref{eq:dOmegadt} as $\Delta \Omega = \dot{\Omega}P/2$. Note $\Delta \Omega - t = \Delta \Omega_\mathrm{R}$ (change in RAAN within the rotating frames) and from the resonance ratio, $\Delta \Omega_{\mathrm{R}} = -2q\pi/2 = -q\pi$. Combining Eqs. \eqref{eq:half_proof_1}-\eqref{eq:half_proof_2},
\begin{align}
    \label{eq:half_proof_3}\mb x_H \\ y_H\me_{t = P/2} = \bm{C}_{(-t-t_0)} \bm{C}_{(\Delta \Omega \pm p\pi)}\mb X \\ Y \me_{t = 0} = \bm{C}_{(\Delta \Omega  - t \pm p\pi)} \mb x_H \\ y_H \me_{t = 0} = \bm{C}_{(\pm p-q)\pi}\mb x_H \\ y_H \me_{t = 0},
\end{align}
as the rotation matrices commute, i.e., order of the rotation does not matter. Combined, the azimuthal angle ($\chi_H$) of the apse rotates by $(\pm p-q)\pi$, returning to the same axis after $t = P/2$. The state at $t = P/2$ is different from the initial state, as the equivalence requires both $p$ and $q$ to be even numbers, violating the coprime property. Thus, an initial state on any symmetric apse set returns to the same set at a distinct location after $t = P/2$.

A similar process follows for the quarter period examination. From Fig. \ref{fig:circular_planar}, after a quarter period,
 \begin{align}
    \label{eq:quarter_proof_2}\mb X \\ Y \me_{t = P/4} & = \alpha \bm{C}_{(\Delta \Omega \pm \frac{p}{2}\pi)} \mb X \\ Y \me_{t = 0},
\end{align} 
where $\alpha = \cos i$ (node) or $1/\cos i$ (apex) depending on the starting apse at $t = 0$; the sign ambiguity for $\pm p\pi/2$ arises from the prograde ($+$) and retrograde ($-$) orbits. As $\Delta M = p\pi/2$, conditions $Z = 0$ or $\dot{Z} = 0$ is satisfied. Then, 
\begin{align}
    \label{eq:quarter_proof_3}\mb x_H \\ y_H\me_{t = P/4} = \alpha  \bm{C}_{(\Delta\Omega - t \pm \frac{p}{2}\pi)}\mb x_H \\ y_H \me_{t = 0} = \alpha \bm{C}_{\frac{\pm p - q}{2}\pi }\mb x_H \\ y_H \me_{t = 0},
\end{align}
where $\Delta \Omega_\mathrm{R} = \Delta \Omega - t =  -q\pi/2$ from the resonance condition. Equation \eqref{eq:quarter_proof_3} suggests that the azimuthal angle ($\chi_H$) is rotating by $(\pm p-q)\pi/2$ per quarter period. Together, starting with a symmetric apse, after a quarter period, another symmetric apse is located. It belongs to the original set only if $p$ is even and $p-q$ is even, violating the assumption that $p,q$ are coprime integers. Thus, the returned symmetric apse set at $t = P/4$ is distinct from the initial apse set at $t = 0$. Combined, any initial symmetric apse condition generated with the circular equilibria results in a doubly symmetric PO configuration within the \acrshort{dadm}. 

\begin{figure}[htbp]
    \centering
    \begin{subfigure}[b]{0.40\textwidth}
        \centering
        \includegraphics[width = 0.99\textwidth]{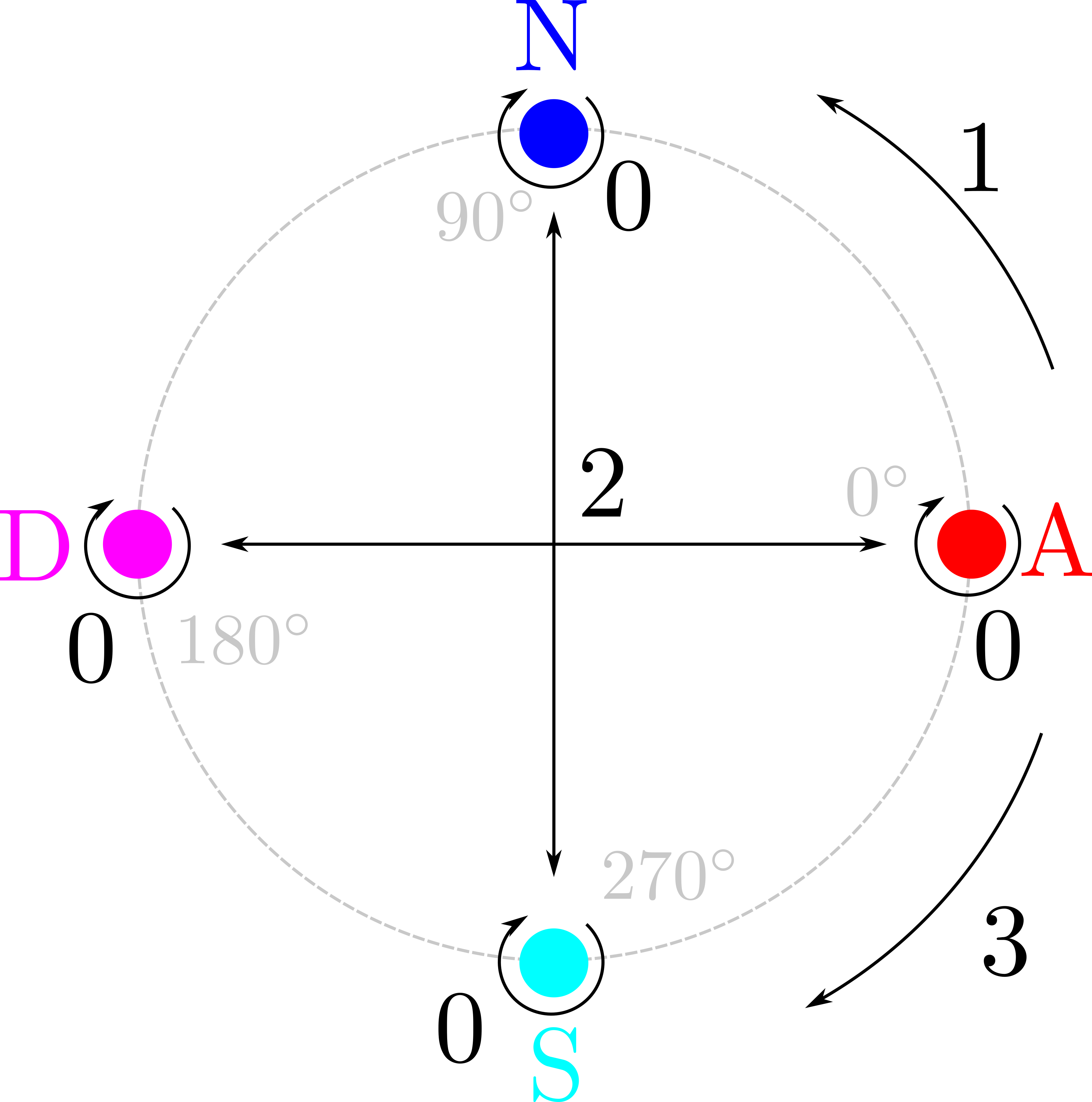}
        \caption{\label{fig:rotation_p} Argument of latitude ($u$) rotation following [$\text{mod}(p,4)$].}
    \end{subfigure}
    \hfill
    \begin{subfigure}[b]{0.40\textwidth}
        \centering
        \includegraphics[width = 0.99\textwidth]{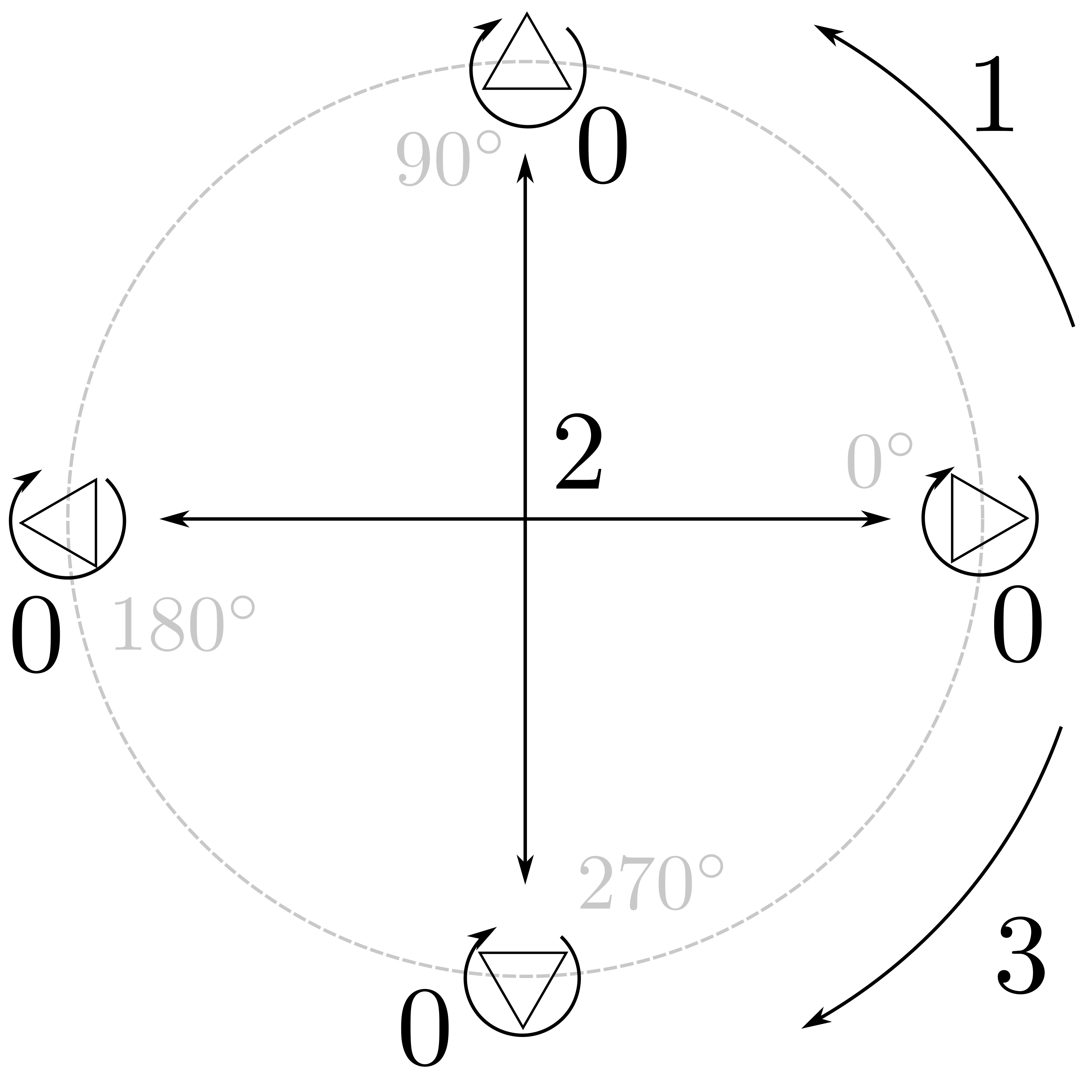}
        \caption{\label{fig:rotation_p_p}Azimuthal angle of the apse ($\chi_H$) rotation following [$\text{mod}(\pm p-q,4)$].}
    \end{subfigure}
    \caption{\label{fig:rotation_modulo}Representation of the modular evolution over a quarter period ($t=0 \to P/4$). The indices dictate the step size of the phase rotation in $u$ and $\chi_H$.}
\end{figure}

The systematic connection between symmetries and the resulting multiplicity of solution geometries is established utilizing modulo-4 arithmetic. As derived in the quarter period evolution (Eqs. \eqref{eq:quarter_proof_2}-\eqref{eq:quarter_proof_3}), the geometric change is governed by the rotation of two distinct phases, as visualized in Fig. \ref{fig:rotation_modulo}. The first index, $\text{mod}(p, 4)$, governs the rotation of the argument of latitude, $u: = \omega + \theta $. For the circular equilibria, $\theta \equiv M$ and $\omega$ may be considered to be zero without losing generality. Since the mean anomaly advances by $p\pi/2$ over a quarter period, the vertical state transitions cyclically through the sequence $\mathrm{A} \to \mathrm{N} \to \mathrm{D} \to \mathrm{S}$, representing a $90^\circ$ phase progression for each integer increment of $\text{mod}(p, 4)$ (Fig. \ref{fig:rotation_p}). Simultaneously, the second index, $\text{mod}(\pm p-q, 4)$, governs the azimuthal rotation of the apse axis ($\chi_H$). This term combines the inertial motion and frame rotation, resulting in a geometric shift of the symmetry axis by $90^\circ$ increments within the rotating frame (Fig. \ref{fig:rotation_p_p}) for each integer increment of $\text{mod}(\pm p-q, 4)$. Again, the sign ambiguity arises from the prograde ($+$) and retrograde ($-$) orbits. 

The sequence of symmetric apsides is determined by the combined evolution of the latitudinal ($u$) and azimuthal ($\chi_H$) angles. Table \ref{tab:modular_classification} compiles the evolution per quarter period depending on $\text{mod}(p, 4)$ and $\text{mod}(q, 4)$. The latitudinal and azimuthal steps are illustrated in Figs. \ref{fig:rotation_p} and \ref{fig:rotation_p_p}, respectively. Notably, the orbital direction (prograde vs. retrograde) dictates the sign of the term $\pm p - q$, thereby directly influencing the azimuthal step $\Delta \chi_H$. The last column of Table \ref{tab:modular_classification} provides sample sequences of symmetric apsides for prograde configurations, utilizing the identifiers and markers defined in Table \ref{tab:symmetry_classification}.

\begin{table*}[htbp]
\centering
\caption{\label{tab:modular_classification}Circular equilibria: evolution of symmetric apsides over a quarter period.}
\begin{tabular}{@{}ccccccl@{}}
\toprule
\textbf{Parity} & \multicolumn{3}{c}{$\mathbf{mod}(\cdot, 4)$} & \textbf{Latitudinal Step} & \textbf{Azimuthal Step} & \multicolumn{1}{c}{\textbf{Sample Sequence ($h_z > 0$)}} \\
\cmidrule(lr){2-4}
($p:q$) & $p$ & $q$ & $\pm p-q$ & ($\Delta u$) & ($\Delta \chi_H$: Pro / Ret) & \multicolumn{1}{c}{($t = 0, P/4, 2P/4, 3P/4$)} \\
\midrule

\multirow{4}{*}{\textbf{odd : odd}} 
 & 1 & 1 & $0\,/\,2$ & $+90^\circ$ & $0^\circ\,/\,180^\circ$ & (\ID{A}{+x}, \ID{N}{+x}, \ID{D}{+x}, \ID{S}{+x}) or  (\xp{ascending}, \xp{northern}, \xp{descending}, \xp{southern}) \\
 & 1 & 3 & $2\,/\,0$ & $+90^\circ$ & $180^\circ\,/\,0^\circ$ & (\ID{A}{+x}, \ID{N}{-x}, \ID{D}{+x}, \ID{S}{-x}) or (\xp{ascending}, \xm{northern}, \xp{descending}, \xm{southern}) \\
 & 3 & 1 & $2\,/\,0$ & $-90^\circ$ & $180^\circ\,/\,0^\circ$ & (\ID{A}{+x}, \ID{S}{-x}, \ID{D}{+x}, \ID{N}{-x}) or (\xp{ascending}, \xm{southern}, \xp{descending}, \xm{northern}) \\
 & 3 & 3 & $0\,/\,2$ & $-90^\circ$ & $0^\circ\,/\,180^\circ$ & (\ID{A}{+x}, \ID{S}{+x}, \ID{D}{+x}, \ID{N}{+x}) or (\xp{ascending}, \xp{southern}, \xp{descending}, \xp{northern}) \\
\midrule

\multirow{4}{*}{\textbf{odd : even}} 
 & 1 & 0 & $1\,/\,3$ & $+90^\circ$ & $+90^\circ\,/\,-90^\circ$ & (\ID{A}{+x}, \ID{N}{+y}, \ID{D}{-x}, \ID{S}{-y}) or (\xp{ascending}, \yp{northern}, \xm{descending}, \ym{southern}) \\
 & 1 & 2 & $3\,/\,1$ & $+90^\circ$ & $-90^\circ\,/\,+90^\circ$ & (\ID{A}{+x}, \ID{N}{-y}, \ID{D}{-x}, \ID{S}{+y}) or (\xp{ascending}, \ym{northern}, \xm{descending}, \yp{southern}) \\
 & 3 & 0 & $3\,/\,1$ & $-90^\circ$ & $-90^\circ\,/\,+90^\circ$ & (\ID{A}{+x}, \ID{S}{-y}, \ID{D}{-x}, \ID{N}{+y}) or (\xp{ascending}, \ym{southern}, \xm{descending}, \yp{northern}) \\
 & 3 & 2 & $1\,/\,3$ & $-90^\circ$ & $+90^\circ\,/\,-90^\circ$ & (\ID{A}{+x}, \ID{S}{+y}, \ID{D}{-x}, \ID{N}{-y}) or (\xp{ascending}, \yp{southern}, \xm{descending}, \ym{northern}) \\
\midrule

\multirow{4}{*}{\textbf{even : odd}} 
 & 0 & 1 & $3\,/\,3$ & $0^\circ$ & $-90^\circ\,/\,-90^\circ$ & (\ID{N}{+x}, \ID{N}{-y}, \ID{N}{-x}, \ID{N}{+y}) or (\xp{northern}, \ym{northern}, \xm{northern}, \yp{northern}) \\
 & 0 & 3 & $1\,/\,1$ & $0^\circ$ & $+90^\circ\,/\,+90^\circ$ & (\ID{N}{+x}, \ID{N}{+y}, \ID{N}{-x}, \ID{N}{-y}) or  (\xp{northern}, \yp{northern}, \xm{northern}, \ym{northern}) \\
 & 2 & 1 & $1\,/\,1$ & $180^\circ$ & $+90^\circ\,/\,+90^\circ$ & (\ID{N}{+x}, \ID{S}{+y}, \ID{N}{-x}, \ID{S}{-y}) or  (\xp{northern}, \yp{southern}, \xm{northern}, \ym{southern}) \\
 & 2 & 3 & $3\,/\,3$ & $180^\circ$ & $-90^\circ\,/\,-90^\circ$ & (\ID{N}{+x}, \ID{S}{-y}, \ID{N}{-x}, \ID{S}{+y}) or  (\xp{northern}, \ym{southern}, \xm{northern}, \yp{southern}) \\
\bottomrule
\end{tabular}
\end{table*}

The interplay between the two rotations of $u$ and $\chi_H$ determines the admissible symmetry pairings for periodic orbits as well. The pairings are governed by the parity in $p:q$ integers and are examined via the evolution per quarter period. For odd-odd parity, the apsides shift their type (node $\leftrightarrow$ apex) while the azimuthal axis rotates by $0^\circ$ or $180^\circ$ (since $\pm p-q$ is even), ensuring the apsides remain aligned with the same axis. Hence, pairings of like-symmetries (OX/XOZ, OY/YOZ) occur. Conversely, for odd-even parity, the azimuthal axis rotates by $\pm 90^\circ$ (since $\pm p-q$ is odd), forcing a transition between apsidal axes and yielding cross-pairings (OX/YOZ, OY/XOZ). Finally, for even-odd parity, the apsides preserve their type (node $\leftrightarrow$ node) while the apsidal axis rotates, connecting the respective axial and reflectional symmetries (OX/OY, XOZ/YOZ). This classification is summarized in the first two columns of Table \ref{tab:symmetry_evolution}. 

The classification in Table \ref{tab:modular_classification} directly determines the multiplicity of distinct solution geometries for each symmetry pairing. Consider a specific resonance $\eta$ with a fixed inclination $0^\circ < i < 90^\circ$ that uniquely determines the semi-major axis $a$. Under these conditions, the 16 symmetric apse configurations in Table \ref{tab:symmetry_classification} each represents a unique state. Since a single periodic orbit traverses a sequence of four unique apsides (as illustrated in Table \ref{tab:modular_classification}), it follows that there are $16/4 = 4$ distinct periodic orbit configurations for the fixed parameters. Enumerating these sequences reveals that the configurations are evenly distributed between the two distinct pairings in Table \ref{tab:symmetry_evolution}. For instance, in the odd-odd parity case, two configurations correspond to the OX/XOZ pairing and two to the OY/YOZ pairing; this breakdown is provided in the third column of Table \ref{tab:symmetry_evolution}. While the \acrshort{hr3bp} admits all four symmetric configurations, the \acrshort{cr3bp} does not, as symmetries involving the $\hat{\bm{y}}$-axis (OY, YOZ) are absent (Table \ref{tab:symmetry_comparison}). This physical constraint leads to a pruning of the solution space: the odd-odd case retains two doubly symmetric geometries (OX/XOZ), whereas the odd-even and even-odd cases are reduced to two \emph{singly symmetric} geometries each (retaining only OX or XOZ), as summarized in the fourth column of Table \ref{tab:symmetry_evolution}.

\begin{table*}[htbp]
\centering
\caption{\label{tab:symmetry_evolution}Evolution of symmetric PO configurations across $p:q$ parity, equilibria type (circular/frozen) as well as unaveraged dynamics (\acrshort{hr3bp}/\acrshort{cr3bp}). The numbers in parentheses indicate the number of distinct symmetric configurations. The acronyms stand for [\emph{DS}: Doubly Symmetric, \emph{SS}: Singly Symmetric, \emph{None}: No symmetric solutions exist.]}
\setlength{\tabcolsep}{5pt}
\begin{tabular}{@{}ccllll@{}}
\toprule
\textbf{Parity} & \textbf{Symmetry Pairing} & \multicolumn{2}{c}{\textbf{Circular Periodic Orbit}} & \multicolumn{2}{c}{\textbf{Frozen Periodic Orbit}} \\
\cmidrule(lr){3-4} \cmidrule(lr){5-6}
($p:q$) & (Quarter-Period) & \multicolumn{1}{c}{\textbf{HR3BP}} & \multicolumn{1}{c}{\textbf{CR3BP}} & \multicolumn{1}{c}{\textbf{HR3BP}} & \multicolumn{1}{c}{\textbf{CR3BP}} \\
\midrule

\multirow{2}{*}[-0.5ex]{\textbf{odd : odd}} 
 & OX--XOZ 
 & \emph{DS}: OX/XOZ (2) & \emph{DS}: OX/XOZ (2) & \emph{SS}: XOZ (4) & \emph{SS}: XOZ (4) \\
 \addlinespace[3pt]
 & OY--YOZ 
 & \emph{DS}: OY/YOZ (2) & \textit{None} & \emph{SS}: YOZ (4) & \textit{None} \\
\midrule

\multirow{2}{*}[-0.5ex]{\textbf{odd : even}} 
 & OX--YOZ 
 & \emph{DS}: OX/YOZ (2) & \emph{SS}: OX (2) & \emph{SS}: YOZ (4) & \textit{None} \\
 \addlinespace[3pt]
 & OY--XOZ 
 & \emph{DS}: OY/XOZ (2) & \emph{SS}: XOZ (2) & \emph{SS}: XOZ (4) & \emph{SS}: XOZ (4) \\
\midrule

\multirow{2}{*}[-0.5ex]{\textbf{even : odd}} 
 & OX--OY 
 & \emph{DS}: OX/OY (2) & \emph{SS}: OX (2) & \textit{None} & \textit{None} \\
 \addlinespace[3pt]
 & XOZ--YOZ 
 & \emph{DS}: XOZ/YOZ (2) & \emph{SS}: XOZ (2) & \emph{DS}: XOZ/YOZ (4) & \emph{SS}: XOZ (4) \\
\bottomrule
\end{tabular}
\end{table*}

\subsubsection{Frozen Equilibria to Analog Periodic Orbits}

The same procedure utilized for the circular equilibria now identifies the admissible symmetric apse configurations and their multiplicity for the frozen equilibria. The geometric evolution of frozen equilibria within the \acrshort{hrf} is evaluated at half and quarter period intervals. The following behaviors are observed:
\begin{itemize}[topsep=0pt, partopsep=0pt]
    \item \textbf{Axially symmetric configurations cannot be generated from frozen equilibria:} The frozen equilibria only satisfy the apse condition ($\bm{R}\cdot\bm{V} = 0$) at $\theta = 0^\circ, 180^\circ$.  The frozen condition also enforces $\sin^2 i > 2/5$, $\omega = \pm 90^\circ$. As such, at the apse, $z_H = Z$ cannot be zero and, thus, the symmetric apse configurations only exist as the reflectional type (Fig. \ref{fig:frozen_planar}). 
    \item \textbf{Half period examination:} An initial state on the reflectional symmetric apse sets ($\Sigma_{\text{XOZ}}$ or $\Sigma_{\text{YOZ}}$, Table \ref{tab:symmetry_comparison}) returns to the same set at a distinct location after $t = P/2$.
    \item \textbf{Quarter period examination:} An initial state on the reflectional symmetric apse sets ($\Sigma_{\text{XOZ}}$ or $\Sigma_{\text{YOZ}}$, Table \ref{tab:symmetry_comparison}) intersects the other reflectional symmetric apse set after $t = P/4$ only when $p$ is even and $q$ is odd. 
\end{itemize}
Derivations similar to the circular case are supplied as follows. Along frozen equilibria, after a half period,   
\begin{align}
    \label{eq:frozen_proof_2}\mb X \\ Y \me_{t = P/2} & = \beta  \bm{C}_{(\Delta \Omega \pm p\pi)} \mb X \\ Y \me_{t = 0},
 \end{align}   
where a scalar $\beta$ is (1) $1$ for an even $p$, (2) $(1+e)/(1-e)$ for an odd $p$ and $f = 0^\circ$ at $t = 0$, and (3) $(1-e)/(1+e)$ for an odd $p$ and $f = 180^\circ$ at $t = 0$ from the schematic in Fig. \ref{fig:frozen_planar}. Then, a similar structure follows from Eq. \eqref{eq:half_proof_3} as, 
\begin{align}
    \mb x_H \\ y_H \me_{t = P/2} = \beta \bm{C}_{(\pm p-q)\pi} \mb x_H \\ y_H \me_{t = 0},
\end{align}
where the scalar $\beta$ does not impact the zero values and, thus, the apsidal rotations. For the quarter period examination, observe that $p$ needs to be an even number to remain at an apse $(\bm{R}\cdot \bm{V} = 0)$ after $t = P/4$. Then, an equation similar to Eq. \eqref{eq:quarter_proof_3} follows as,
\begin{align}
    \mb x_H \\ y_H\me_{t = P/4} = \beta' \bm{C}_{\frac{\pm p - q}{2}\pi }\mb x_H \\ y_H \me_{t = 0},
\end{align}
where $\beta'$ is adjusted for the quarter period but defined similarly as $\beta$. To alternate between different symmetric apse sets, $\frac{\pm p - q}{2}\pi = \pm \frac{\pi}{2}$. Thus, $q$ is an odd number. 

Based on these properties, the classification of the symmetric configurations is determined from the resonance parity, similar to the circular case. For resonances where $p$ is odd (for both odd-odd and odd-even), the orbit does not reach an apse condition at the quarter period mark ($t=P/4$), rendering singly symmetric geometries. Conversely, for the even-odd case where $p$ is even, the orbit satisfies the symmetric apse condition at the quarter period, resulting in doubly symmetric geometries connecting the XOZ and YOZ symmetries. The modulo-4 arithmetic as detailed in Table \ref{tab:modular_classification} generated for the circular equilibria remains informative to the frozen equilibria as well; half period apsidal sequences are available via applying the quarter period evolution twice. However, a critical distinction arises within the \acrshort{cr3bp}, where the YOZ symmetry is non-existent. As a result, these potentially doubly symmetric geometries are effectively ``pruned'' to singly symmetric (XOZ) configurations. This information is detailed in the last two columns of Table \ref{tab:symmetry_evolution}.

The frozen equilibria exhibit a higher multiplicity than the circular equilibria in terms of distinct configurations. In circular equilibria, the geometry is invariant with respect to the argument of periapsis. In contrast, for frozen orbits, the periapsis location relative to the secondary body ($\mathcal{P}_2$) creates two distinct configurations at $\omega = \pm 90^\circ$. This breakage doubles the number of distinct geometrical configurations, as summarized in Table \ref{tab:symmetry_evolution}. For instance, in the odd-odd parity, while the circular equilibria yield two doubly symmetric geometries (OX/XOZ), the frozen equilibria render four singly symmetric configurations (XOZ), reflecting the gain of $\omega = \pm 90^\circ$ multiplicity.

\subsection{\label{sec:initialization_targeting}Transitioning to Unaveraged Dynamics}

\begin{figure}[htbp] 
    \centering
\begin{tikzpicture}[
    >=Stealth,
    font=\footnotesize,
    block/.style={
        draw, rectangle, rounded corners, 
        align=center, minimum width=2.4cm, minimum height=0.9cm,
        inner sep=3pt, fill=white
    },
    slimblock/.style={
        draw, rectangle, rounded corners, 
        align=center, minimum width=2.4cm, minimum height=0.5cm,
        inner sep=2pt, fill=white
    },
    container/.style={
        draw, dashed, inner sep=0.15cm, rounded corners, fill=gray!10
    },
    arrowlabel/.style={
        midway, above, align=center, font=\scriptsize, color=black,
        text width=2.0cm
    }
]

    \def\coldist{4.2cm}  
    \def\vshift{1.5cm}   

    \coordinate (c1) at (0, 0);              
    \coordinate (c2) at (\coldist, 0);       
    \coordinate (c3) at (2*\coldist, 0);     
    \coordinate (c4) at (3*\coldist, 0);     

    
    \node[slimblock] (box2) at (0, \vshift) {Fixed $a$ or $i$};
    \node[slimblock, above=0.15cm of box2] (box1) {$\eta = p/q$};
    \node[block, below=0.15cm of box2] (box3) {Eq. Type\\(Circular/Frozen)};

    \node[slimblock] (box4) at (0, -\vshift + 0.35cm) {Symmetry Type};
    \node[slimblock, below=0.15cm of box4] (box5) {Initial Apse};


    \node[block, minimum height=0.6cm] (box6) at (\coldist, 0.7cm) {$a, e, i, \omega$};
    
    \node[block, minimum height=0.6cm] (box7) at (\coldist, -0.7cm) {$\Omega, M, t_0$};

    \node[block] (box8) at (c3) {Initial Guess\\($\bm{r}_H, \bm{v}_H$) or ($\bm{r}, \bm{v}$)};

    \node[block] (box9) at (c4) {Converged \\ ($\bm{r}_H, \bm{v}_H$) or ($\bm{r}, \bm{v}$)};

    \begin{scope}[on background layer]
        \node[container, fit=(box1) (box2) (box3), label={[font=\tiny, text=gray]above:Input Parameters}] (group1) {};
        
        \node[container, fit=(box4) (box5), label={[font=\tiny, text=gray]below:Symmetry Inputs}] (group2) {};
        
        \node[container, fit=(box6) (box7), label={[anchor=south east, font=\tiny, text=gray]south east:}] (box10) {};
    \end{scope}


    \draw[->] (group1.east) -- (box6.west) 
        node[arrowlabel] {Eqs. \eqref{eq:nu_m_circular}--\eqref{eq:nu_m_frozen}};

    \draw[->] (group2.east) -- (box7.west) 
        node[arrowlabel] {Tables \ref{tab:modular_classification}--\ref{tab:initialization_matrix}};

    \draw[->] (box10.east) -- (box8.west) 
        node[arrowlabel] {Kinematic\\Rotation};

    \draw[->] (box8.east) -- (box9.west) 
        node[arrowlabel] {Differential\\Corrector};
        
\end{tikzpicture}
\caption{\label{fig:schematic_init}Process for transitioning \acrshort{dadm} equilibria to analog symmetric POs within the unaveraged dynamics.}
\end{figure}
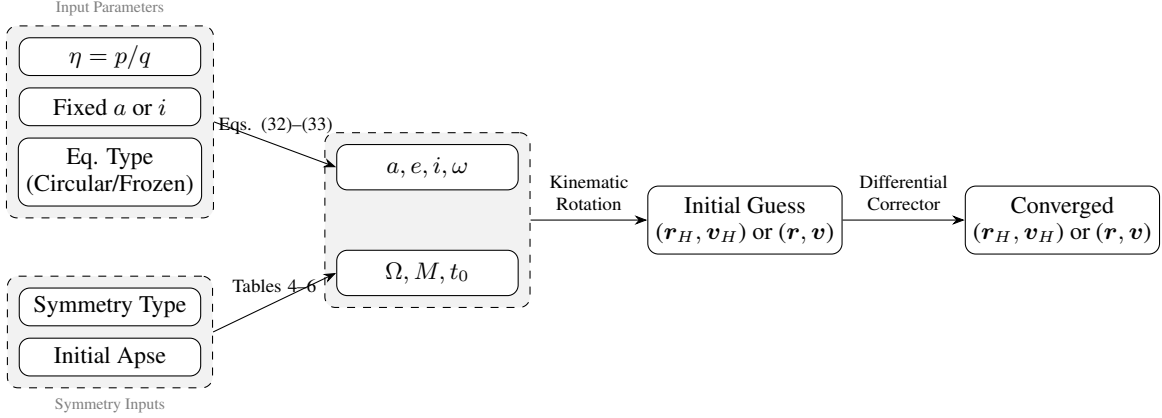

Building on the preceding analysis, the transition of \acrshort{dadm} equilibria into analog symmetric periodic orbits within the \acrshort{cr3bp} and \acrshort{hr3bp} is accomplished. The overall process is summarized schematically in Fig. \ref{fig:schematic_init}. Initialization requires specific user inputs. For a given resonance ratio $\eta = p/q$, fixing either the semi-major axis ($a$) or inclination ($i$) uniquely determines the remaining parameter via the frequency relations in Eqs. \eqref{eq:nu_m_circular}--\eqref{eq:nu_m_frozen}. Subsequently, the equilibrium type dictates the eccentricity ($e$) and argument of periapsis ($\omega$). For circular equilibria, $e = 0$, and the argument of periapsis is set to $\omega = 0^\circ$ without loss of generality. For frozen equilibria, the eccentricity is defined as $e = \sqrt{(5\sin^2 i - 2)/3}$ with $\omega = \pm 90^\circ$. Additional necessary inputs include the desired symmetry pairing (e.g., OX/XOZ) and the specific initial symmetric configuration at $t_0$. The rotational logic depicted in Fig. \ref{fig:rotation_modulo} then yields the remaining orbital elements, $\Omega$ and $M$. Specifically, the initial argument of latitude $u = \omega + \theta $ is selected from $\{0^\circ, 90^\circ, 180^\circ, 270^\circ\}$ to correspond with the apsidal states $\{ \mathrm{A, N, D, S} \}$, respectively. The azimuthal angle of the apse, $\chi_H$, is determined by $\Omega$ and the initial phase $t_0$. To enforce $\Omega \equiv \chi_H$, the initial phase ($t_0$) is defined as $t_0 = \pm u$, where the sign corresponds to prograde ($+$) or retrograde ($-$) motion. This phase offset ensures the configuration aligns with the $+\hat{\bm{x}}$ direction when $\Omega = 0^\circ$. Consequently, the set $\Omega \in \{0^\circ, 90^\circ, 180^\circ, 270^\circ\}$ maps one-to-one to the apse locations at $\{ +\hat{\bm{x}}, +\hat{\bm{y}}, -\hat{\bm{x}}, -\hat{\bm{y}} \}$ relative to the \acrshort{hrf}.

Table \ref{tab:initialization_matrix} summarizes the required initialization parameters ($u, \Omega$) for each target configuration. Since OY and YOZ symmetric configurations do not exist within the \acrshort{cr3bp}, the azimuthal phases $\Omega = \chi_H = 90^\circ$ and $270^\circ$ are excluded (indicated by grey columns). Similarly, $u = 0^\circ$ and $180^\circ$ are inapplicable to frozen equilibria, as axial configurations are invalid within this context (indicated by grey rows). The resulting six-dimensional state in the \acrshort{eof} is transformed to provide an initial guess for the unaveraged dynamics, i.e., $(\bm{r}_H, \bm{v}_H)$ for the \acrshort{hr3bp} and $(\bm{r}, \bm{v})$ for the \acrshort{cr3bp}\footnote{In transforming the orbital elements to position and velocity, this work assumes the doubly averaged elements (governing the \acrshort{dadm} dynamics in Eqs. \eqref{eq:dadt}-\eqref{eq:dthetadt}) are equivalent to the osculating elements. While higher-fidelity transition strategies exist, e.g., from \citet{nie2018lunar}, they are omitted here as differential correctors sufficiently bridge the dynamical gap.}. Differential correctors are then employed to target the symmetric apse conditions at $t = P/2$ (for singly symmetric orbits) or $t = P/4$ (for doubly symmetric orbits), following standard procedures, e.g., by \citet{robin1980numerical}. For orbits with substantially long periods, the initial trajectory derived from the \acrshort{dadm} may be discretized into multiple segments to ensure robust numerical convergence. Upon convergence, the analog symmetric periodic orbits are supplied within the respective unaveraged dynamics.

\begin{table}[htbp]
\centering
\caption{\label{tab:initialization_matrix}Initialization lookup table for $u = \omega + \theta$ and $\Omega$ in transitioning \acrshort{dadm} equilibria to the \acrshort{hr3bp} and \acrshort{cr3bp}. The grey rows are pruned for frozen equilibria, and the grey columns are pruned for the \acrshort{cr3bp} dynamics. Initial phase angle $t_0 = \pm u$ ($+$, $-$ for prograde, retrograde orbits).}
\setlength{\tabcolsep}{6pt}
\renewcommand{\arraystretch}{1.6}

\begin{tabular}{@{}c|cccc@{}}
\toprule
\multirow{2}{*}{\textbf{Latitudinal Phase}} & \multicolumn{4}{c}{\textbf{Azimuthal Phase} ($\Omega \equiv \chi_H$)} \\
\cmidrule(l){2-5}
($u$) & 
\textbf{0$^\circ$} ($+\hat{\bm{x}}$) & 
\cellcolor{prune_light}\textbf{90$^\circ$} ($+\hat{\bm{y}}$) & 
\textbf{180$^\circ$} ($-\hat{\bm{x}}$) & 
\cellcolor{prune_light}\textbf{270$^\circ$} ($-\hat{\bm{y}}$) \\ 
\midrule

\cellcolor{prune_light}\textbf{0$^\circ$} (Ascending) & 
\cellcolor{prune_light}\ID{A}{+x} \enspace \xp{ascending} & 
\cellcolor{prune_dark}\ID{A}{+y} \enspace \yp{ascending} & 
\cellcolor{prune_light}\ID{A}{-x} \enspace \xm{ascending} & 
\cellcolor{prune_dark}\ID{A}{-y} \enspace \ym{ascending} \\

\textbf{90$^\circ$} (North) & 
\cellcolor{white}\ID{N}{+x} \enspace \xp{northern} & 
\cellcolor{prune_light}\ID{N}{+y} \enspace \yp{northern} & 
\cellcolor{white}\ID{N}{-x} \enspace \xm{northern} & 
\cellcolor{prune_light}\ID{N}{-y} \enspace \ym{northern} \\

\cellcolor{prune_light}\textbf{180$^\circ$} (Descending) & 
\cellcolor{prune_light}\ID{D}{+x} \enspace \xp{descending} & 
\cellcolor{prune_dark}\ID{D}{+y} \enspace \yp{descending} & 
\cellcolor{prune_light}\ID{D}{-x} \enspace \xm{descending} & 
\cellcolor{prune_dark}\ID{D}{-y} \enspace \ym{descending} \\

\textbf{270$^\circ$} (South) & 
\cellcolor{white}\ID{S}{+x} \enspace \xp{southern} & 
\cellcolor{prune_light}\ID{S}{+y} \enspace \yp{southern} & 
\cellcolor{white}\ID{S}{-x} \enspace \xm{southern} & 
\cellcolor{prune_light}\ID{S}{-y} \enspace \ym{southern} \\

\bottomrule
\end{tabular}
\end{table}

\begin{figure}[htbp]
    \centering
    \begin{subfigure}[b]{0.49\textwidth}
        \centering
        \includegraphics[width = 1.0\textwidth]{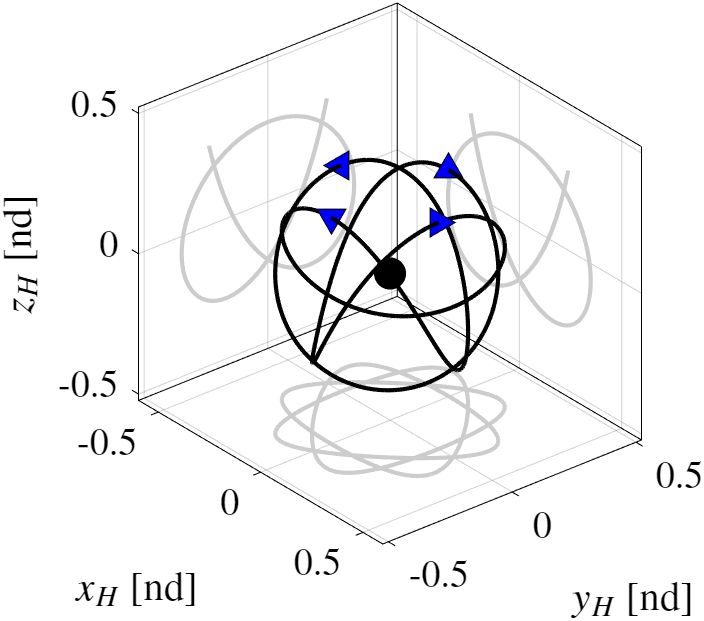}
        \caption{\label{fig:dadm_4_1_circular_N}Geometry 1: $\xp{northern}\to\ym{northern}\to\xm{northern}\to\yp{northern}$.}
    \end{subfigure}
    \hfill
    \begin{subfigure}[b]{0.49\textwidth}
        \centering
        \includegraphics[width = 1.0\textwidth]{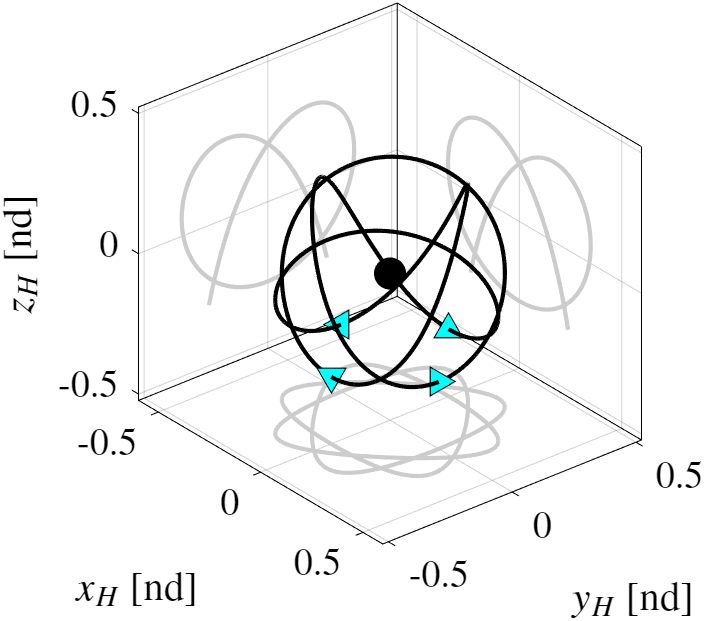}
        \caption{\label{fig:dadm_4_1_circular_S}Geometry 2: $\xp{southern}\to\ym{southern}\to\xm{southern}\to\yp{southern}$.}
    \end{subfigure}
    \caption{\label{fig:dadm_4_1_circular}Doubly symmetric (XOZ/YOZ) circular PO in the \acrshort{dadm} ($p = 4$, $q = 1$, $a_H \approx 0.37$ [nd], $i = 50^\circ$).}
\end{figure}

\begin{figure}[htbp]
    \centering
    \begin{subfigure}[b]{0.49\textwidth}
        \centering
        \includegraphics[width = 1.0\textwidth]{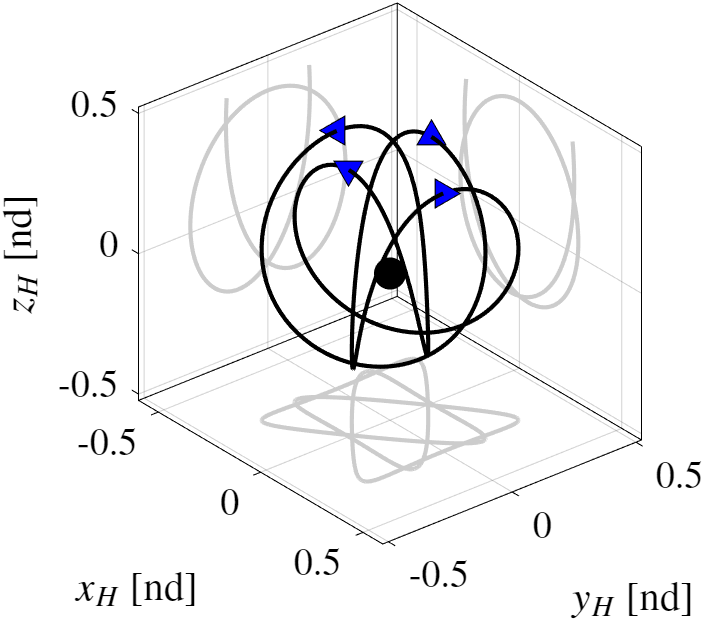}
        \caption{\label{fig:hr3bp_4_1_circular_N}Geometry 1: $\xp{northern}\to\ym{northern}\to\xm{northern}\to\yp{northern}$.}
    \end{subfigure}
    \hfill
    \begin{subfigure}[b]{0.49\textwidth}
        \centering
        \includegraphics[width = 1.0\textwidth]{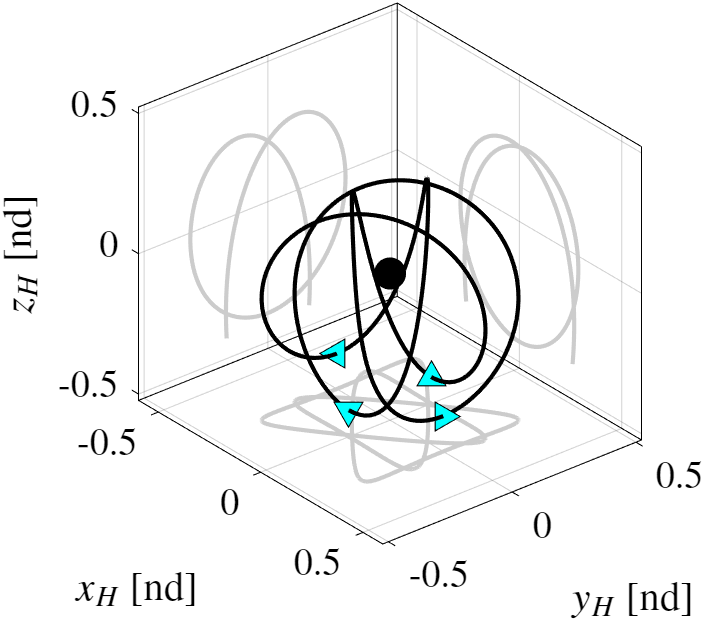}
        \caption{\label{fig:hr3bp_4_1_circular_S}Geometry 2: $\xp{southern}\to\ym{southern}\to\xm{southern}\to\yp{southern}$.}
    \end{subfigure}
    \caption{\label{fig:hr3bp_4_1_circular}Doubly symmetric (XOZ/YOZ) circular PO in the \acrshort{hr3bp} (initialized from Fig. \ref{fig:dadm_4_1_circular}). }
\end{figure}

\begin{figure}[htbp]
    \centering
    \begin{subfigure}[b]{0.49\textwidth}
        \centering
        \includegraphics[width = 1.0\textwidth]{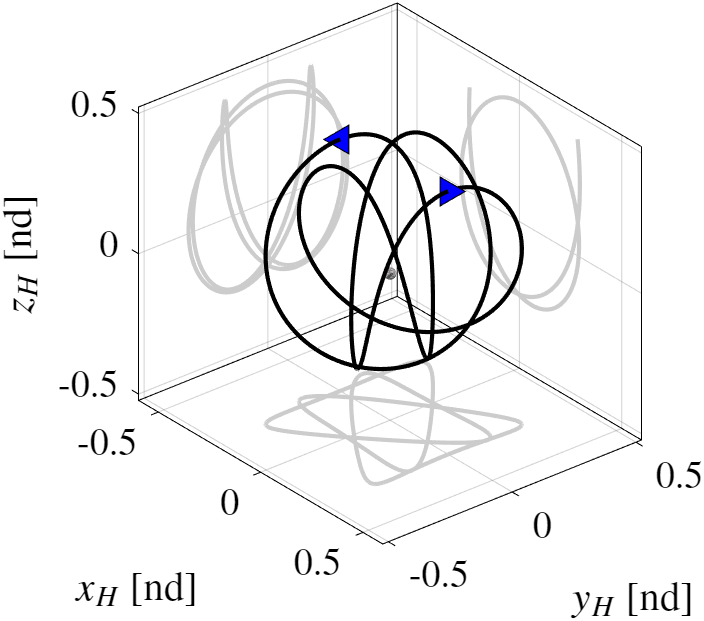}
        \caption{\label{fig:cr3bp_4_1_circular_N}Geometry 1: $\xp{northern}\to\xm{northern}$.}
    \end{subfigure}
    \hfill
    \begin{subfigure}[b]{0.49\textwidth}
        \centering
        \includegraphics[width = 1.0\textwidth]{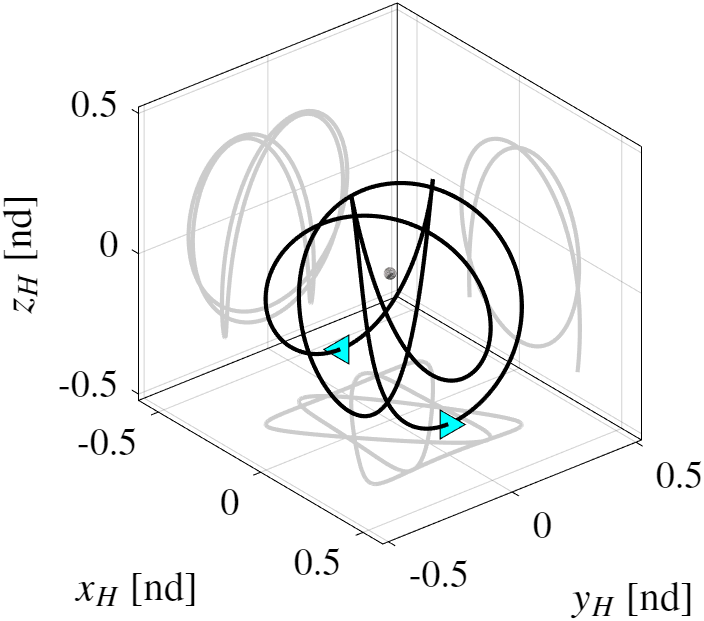}
        \caption{\label{fig:cr3bp_4_1_circular_S}Geometry 2: $\xp{southern}\to\xm{southern}$.}
    \end{subfigure}
    \caption{\label{fig:cr3bp_4_1_circular}Singly symmetric (XOZ) circular PO in the \acrshort{cr3bp} (initialized from Fig. \ref{fig:dadm_4_1_circular}). }
\end{figure}

Consider the $p = 4, q = 1$ resonance as an illustrative example for periodic orbits originating from circular equilibria. Assuming an inclination of $i = 50^\circ$, the resonance condition uniquely determines the semi-major axis as $a_H \approx 0.37$ [nd]. According to Table \ref{tab:symmetry_evolution}, this even-odd parity predicts the existence of doubly symmetric (XOZ/YOZ) configurations within the \acrshort{hr3bp}. Examination of Table \ref{tab:modular_classification} reveals that two distinct configurations exist, following the symmetric apsidal sequences as (1) $\xp{northern}\to\ym{northern}\to\xm{northern}\to\yp{northern}$ or (2)  $\xp{southern}\to\ym{southern}\to\xm{southern}\to\yp{southern}$. Based on the initialization (Table \ref{tab:initialization_matrix}), these two geometries are initialized utilizing the $(\Omega, M)$ combinations of $(0^\circ, 90^\circ)$ and $(0^\circ, 270^\circ)$, respectively. Figure \ref{fig:dadm_4_1_circular} depicts the resulting trajectories in the \acrshort{dadm}, clearly demonstrating the symmetric configurations within the \acrshort{hrf}. These configurations are then differentially corrected in the \acrshort{hr3bp} model by fixing the nd period and targeting the quarter period symmetric apse condition. The analog circular periodic orbit for the \acrshort{hr3bp} is plotted in Fig. \ref{fig:hr3bp_4_1_circular}. Here, the radius of the secondary body is arbitrarily set to $0.05$ [nd]. While the geometry deviates significantly from the \acrshort{dadm} solution, the qualitative symmetry characteristics and the number of distinct configurations are preserved. In the \acrshort{cr3bp} (Fig. \ref{fig:cr3bp_4_1_circular}), representing the Earth-Moon system, the doubly symmetric orbits degenerate into singly symmetric (XOZ) orbits due to the absence of YOZ symmetry. Nevertheless, the apsidal evolution at the half-period mark remains consistent with the predictions in Table \ref{tab:modular_classification}.

\begin{figure}[htbp]
    \centering
    \begin{subfigure}[b]{0.49\textwidth}
        \centering
        \includegraphics[width = 1.0\textwidth]{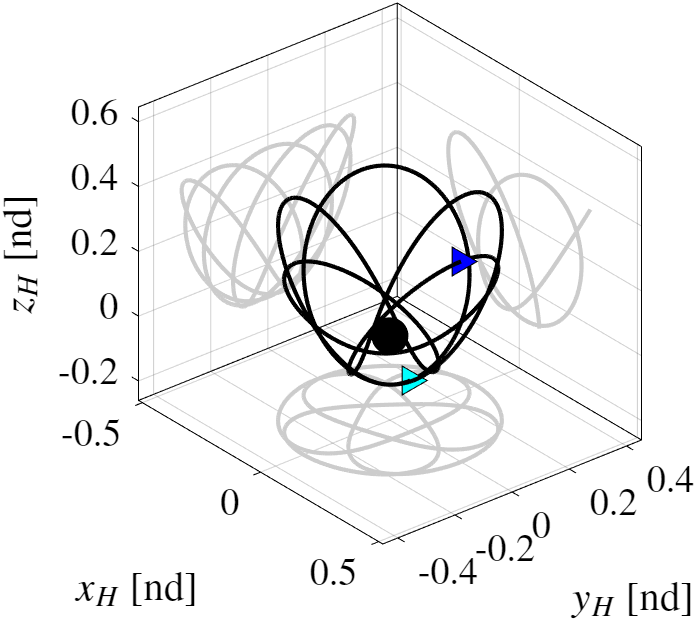}
        \caption{\label{fig:dadm_5_1_frozen_n_xp}Geometry 1: $\xp{northern}\to\xp{southern}$ with $\omega = -90^\circ$.}
    \end{subfigure}
    \hfill
    \begin{subfigure}[b]{0.49\textwidth}
        \centering
        \includegraphics[width = 1.0\textwidth]{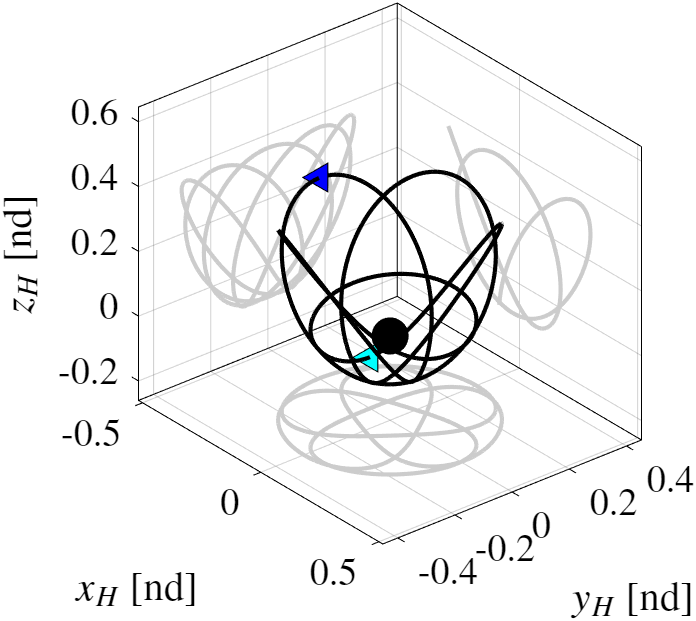}
        \caption{\label{fig:dadm_5_1_frozen_n_xm}Geometry 2: $\xm{northern}\to\xm{southern}$ with $\omega = -90^\circ$.}
    \end{subfigure}
    \begin{subfigure}[b]{0.49\textwidth}
        \centering
        \includegraphics[width = 1.0\textwidth]{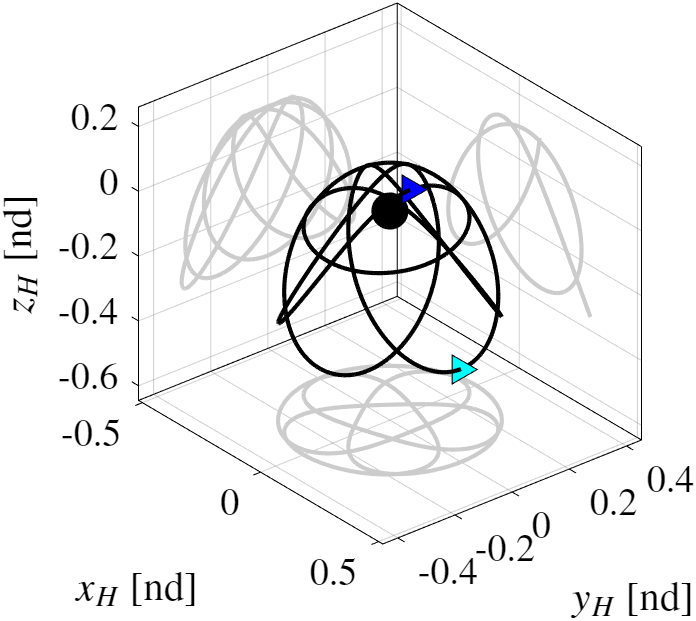}
        \caption{\label{fig:dadm_5_1_frozen_s_xp}Geometry 3: $\xp{northern}\to\xp{southern}$ with $\omega = 90^\circ$.}
    \end{subfigure}
    \hfill
    \begin{subfigure}[b]{0.49\textwidth}
        \centering
        \includegraphics[width = 1.0\textwidth]{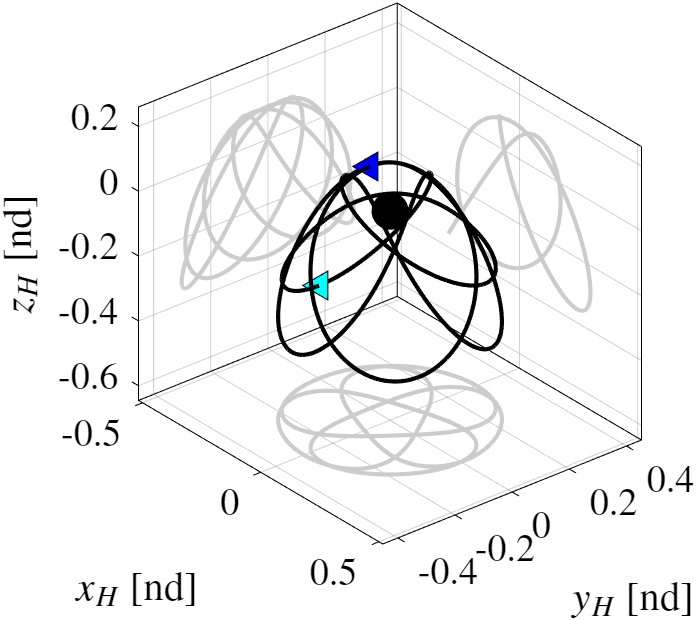}
        \caption{\label{fig:dadm_5_1_frozen_s_xm}Geometry 4: $\xm{northern}\to\xm{southern}$ with $\omega = 90^\circ$.}
    \end{subfigure}
    \caption{\label{fig:dadm_5_1_frozen}Singly symmetric (XOZ) frozen PO in the \acrshort{dadm} ($p = 5$, $q = 1$, $a_H \approx 0.30$ [nd], $i = 50^\circ$).}
\end{figure}

\begin{figure}[htbp]
    \centering
    \begin{subfigure}[b]{0.49\textwidth}
        \centering
        \includegraphics[width = 1.0\textwidth]{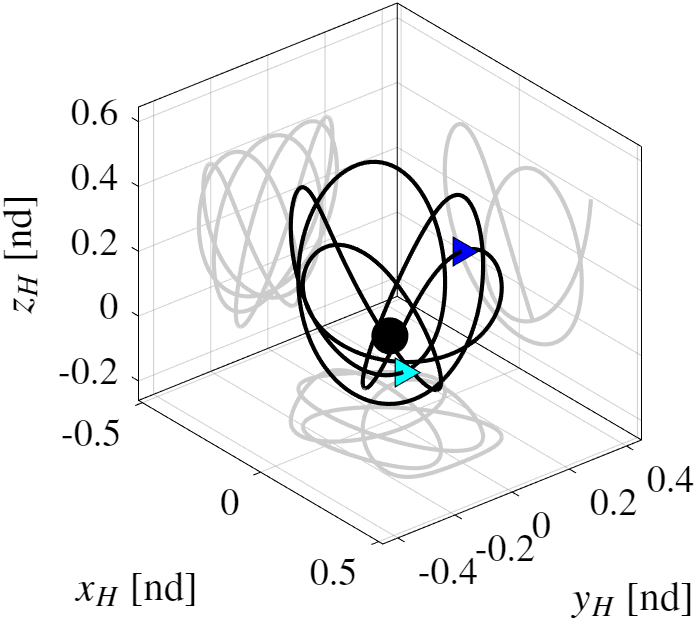}
        \caption{\label{fig:hr3bp_5_1_frozen_n_xp}Geometry 1: $\xp{northern}\to\xp{southern}$ with $\omega = -90^\circ$.}
    \end{subfigure}
    \hfill
    \begin{subfigure}[b]{0.49\textwidth}
        \centering
        \includegraphics[width = 1.0\textwidth]{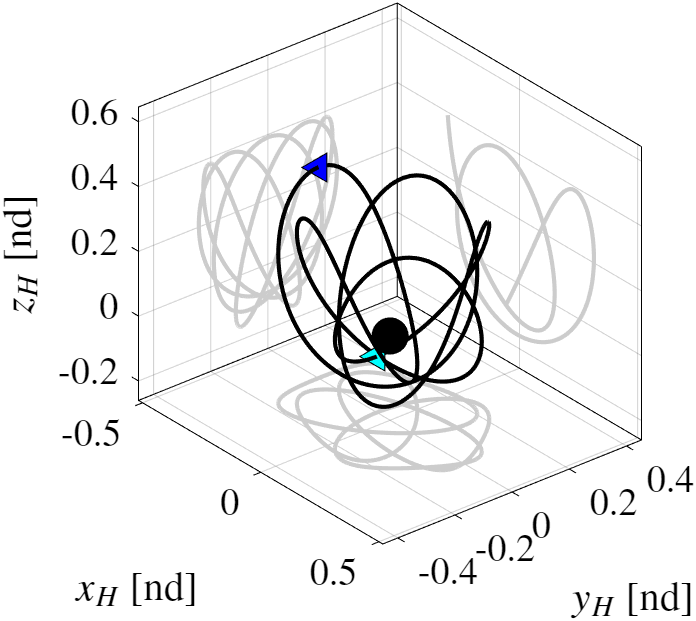}
        \caption{\label{fig:hr3bp_5_1_frozen_n_xm}Geometry 2: $\xm{northern}\to\xm{southern}$ with $\omega = -90^\circ$.}
    \end{subfigure}
    \begin{subfigure}[b]{0.49\textwidth}
        \centering
        \includegraphics[width = 1.0\textwidth]{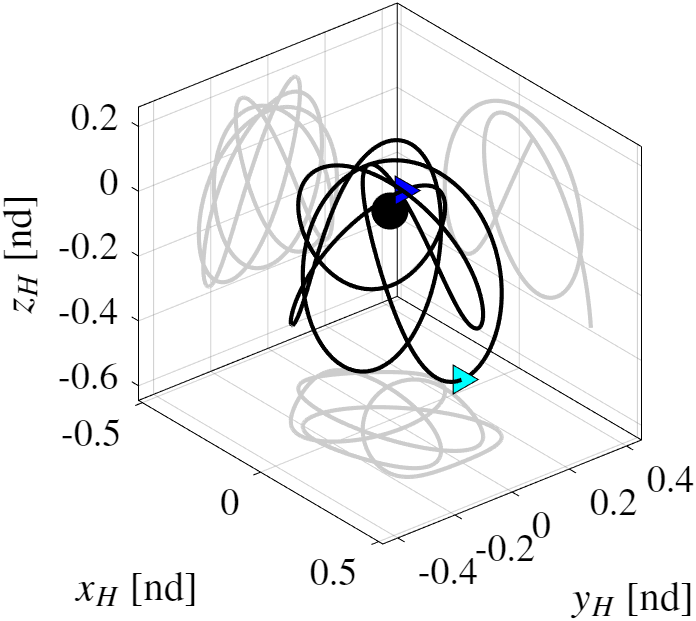}
        \caption{\label{fig:hr3bp_5_1_frozen_s_xp}Geometry 3: $\xp{northern}\to\xp{southern}$ with $\omega = 90^\circ$.}
    \end{subfigure}
    \hfill
    \begin{subfigure}[b]{0.49\textwidth}
        \centering
        \includegraphics[width = 1.0\textwidth]{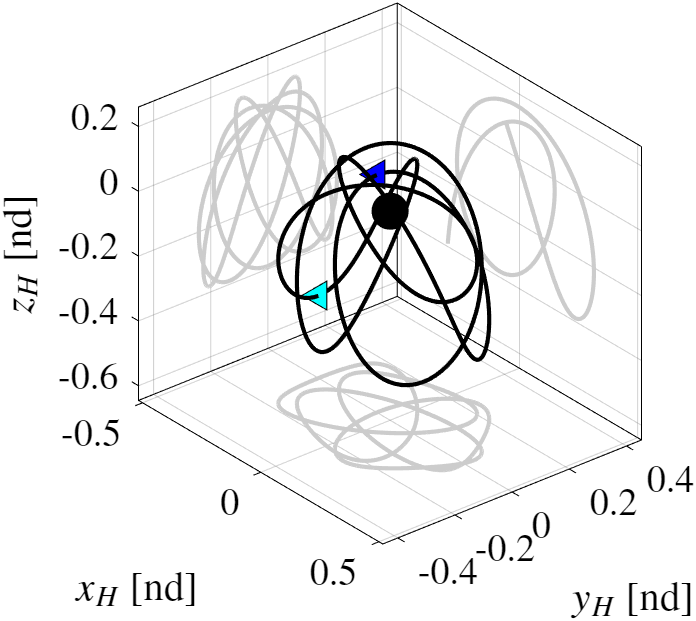}
        \caption{\label{fig:hr3bp_5_1_frozen_s_xm}Geometry 4: $\xm{northern}\to\xm{southern}$ with $\omega = 90^\circ$.}
    \end{subfigure}
    \caption{\label{fig:hr3bp_5_1_frozen}Singly symmetric (XOZ) frozen PO in the \acrshort{hr3bp} (initialized from Fig. \ref{fig:dadm_5_1_frozen}). }
\end{figure}

\begin{figure}[htbp]
    \centering
    \begin{subfigure}[b]{0.49\textwidth}
        \centering
        \includegraphics[width = 1.0\textwidth]{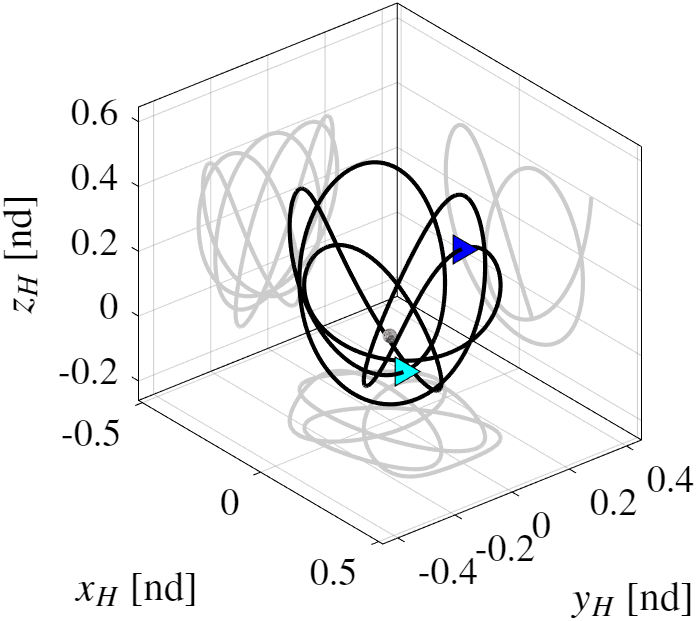}
        \caption{\label{fig:cr3bp_5_1_frozen_n_xp}Geometry 1: $\xp{northern}\to\xp{southern}$ with $\omega = -90^\circ$.}
    \end{subfigure}
    \hfill
    \begin{subfigure}[b]{0.49\textwidth}
        \centering
        \includegraphics[width = 1.0\textwidth]{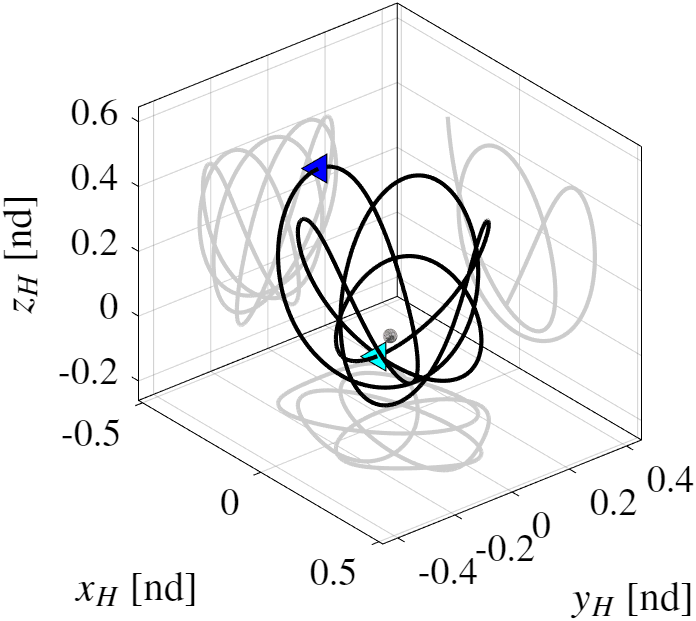}
        \caption{\label{fig:cr3bp_5_1_frozen_n_xm}Geometry 2: $\xm{northern}\to\xm{southern}$ with $\omega = -90^\circ$.}
    \end{subfigure}
    \begin{subfigure}[b]{0.49\textwidth}
        \centering
        \includegraphics[width = 1.0\textwidth]{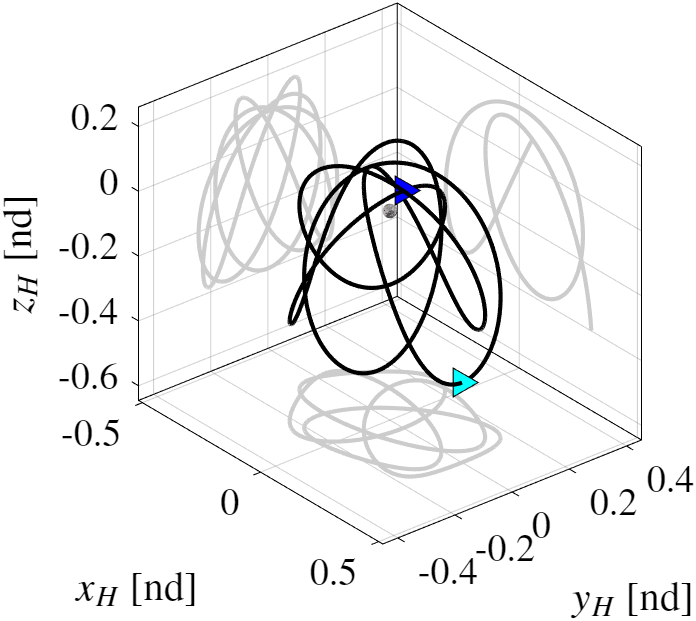}
        \caption{\label{fig:cr3bp_5_1_frozen_s_xp}Geometry 3: $\xp{northern}\to\xp{southern}$ with $\omega = 90^\circ$.}
    \end{subfigure}
    \hfill
    \begin{subfigure}[b]{0.49\textwidth}
        \centering
        \includegraphics[width = 1.0\textwidth]{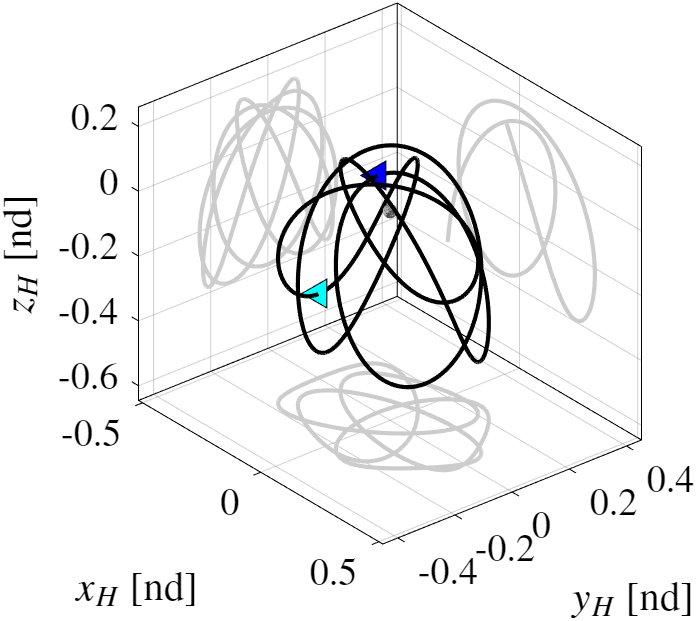}
        \caption{\label{fig:cr3bp_5_1_frozen_s_xm}Geometry 4: $\xm{northern}\to\xm{southern}$ with $\omega = 90^\circ$.}
    \end{subfigure}
    \caption{\label{fig:cr3bp_5_1_frozen}Singly symmetric (XOZ) frozen PO in the \acrshort{cr3bp} (initialized from Fig. \ref{fig:dadm_5_1_frozen}).}
\end{figure}

A second example examines the frozen equilibria with $p = 5$ and $q = 1$. Assuming an inclination of $i = 50^\circ$ leads to $a_H \approx 0.30$ [nd]. The odd-odd parity allows a singly symmetric configuration with XOZ across the \acrshort{hr3bp} and \acrshort{cr3bp} (Table \ref{tab:symmetry_evolution}). Two possible apsidal sequences are: (1) $\xp{northern}\to\xp{southern}$ and (2) $\xm{northern}\to\xm{southern}$ from the logic supplied in Table \ref{tab:modular_classification}. As indicated in Table \ref{tab:symmetry_evolution}, the number of configurations doubles to four due to the additional degree of freedom in selecting $\omega = \pm 90^\circ$. In generating an initial symmetric apse at $\xp{northern}$, $\Omega = 0, f = 0^\circ, \omega = 90^\circ$ is incorporated. Repeating the process for other geometries, the symmetric configurations from the \acrshort{dadm} are illustrated in Fig. \ref{fig:dadm_5_1_frozen}. The first and second rows correspond to northern ($\omega = -90^\circ$) and southern ($\omega = 90^\circ$) configurations\footnote{Naming convention for frozen periodic orbit follows that of the halo orbits within the \acrshort{r3bps}, dictated by the location of the apoapsis relative to the $xy$-plane.}. The two columns each represent the apsidal sequence for (1) $\xp{northern}\to\xp{southern}$ and (2) $\xm{northern}\to\xm{southern}$. The corresponding periodic orbits are converged in the \acrshort{hr3bp} and \acrshort{cr3bp}, plotted in Figs. \ref{fig:hr3bp_5_1_frozen} and \ref{fig:cr3bp_5_1_frozen}, respectively. Notably, even with significant geometric variations arising from the dynamical disparity, the symmetric apsidal configurations and the multiplicity remain consistent across multiple models as expected from Table \ref{tab:symmetry_evolution}. This consistency validates the proposed transition strategy for locating symmetric periodic orbits within the \acrshort{r3bps} originating from the \acrshort{dadm} equilibria.

\clearpage

\section{\label{sec:bifurcation}Families of Symmetric Periodic Orbits: Connections}

While the preceding section focuses on the initialization of individual symmetric periodic orbits derived from the \acrshort{dadm}, the current section expands the scope to a global characterization. Within the \acrshort{r3bps}, periodic orbits generally organize into one-parameter families \cite{henon2002generating}, interconnected via bifurcations. This section investigates the evolution of these families, explicitly linking the symmetric POs lifted from averaged equilibria to well-known dynamical structures. Particular attention is focused on the connections between circular and frozen POs, where pitchfork bifurcations in the averaged model often manifest as ``broken'' bifurcations in the unaveraged dynamics. The framework of frequency, symmetry, and resonance parity ($p:q$) established in the preceding sections serves as the foundation for this analysis.

\subsection{\label{sec:circular_bif}Families of Circular Periodic Orbits}

The families of circular periodic orbits within the \acrshort{r3bps} are expected to evolve smoothly across the inclination range. This continuous behavior is visualized in Fig. \ref{fig:a_i_eta_C_hr3bp}, where the contour lines for constant resonance ratios $\eta = p/q = \nu_s/\nu_m$ exhibit no singularities for $0^\circ < i < 180^\circ$. Such continuity is consistent with the analytical expression derived in the \acrshort{dadm} (Eq. \eqref{eq:nu_m_circular}), where the medium-period frequency $\nu_m$ varies continuously with inclination. However, a kinematic distinction arises at the planar limits ($i = 0^\circ$ and $i = 180^\circ$). Recall that the short-period frequency $\nu_s$ determines the $p$ revolutions measured in the inertial frame (\acrshort{eof}). A subtlety arises due to the combined dynamics of the rotating frame that revolves at a constant unit rate, and the secular drift of the \acrshort{raan} ($\Omega$). These factors jointly alter the geometric ``closure points'' of the circular orbits relative to the rotating frame. This relative motion is effectively quantified by the rate of the nodal line as measured within the rotating frame, i.e., $\dot{\Omega}_{\mathrm{R}} = \dot{\Omega} - 1$. Consequently, although the inertial periodicity involves $p$ revolutions, the satellite effectively completes $p \mp q$ revolutions within the rotating frame to compensate for this relative nodal drift and form a closed periodic loop. For prograde orbits ($i= 0^\circ$), the nodal precession acts to reduce the apparent revolutions, leading to the `$-$' sign ($p-q$). Conversely, for retrograde orbits ($i = 180^\circ$), the relative motion leads to the `$+$' sign ($p+q$). Thus, the effective (synodic) resonance ratios for the planar limits are consistently defined as,
\begin{align}
    \label{eq:synodic}\eta_{0^\circ} = \frac{\nu_s - \nu_m}{\nu_m} = \eta -1, \quad \eta_{180^\circ} = \frac{\nu_s + \nu_m}{\nu_m} = \eta + 1.
\end{align}
Building on these properties, the connections to well-known dynamical structures in the \acrshort{r3bps} are discussed subsequently for the prograde and retrograde limits.

\subsubsection{\label{sec:circular_bif_prograde}Connections at Planar, Prograde ($i = 0^\circ$)}

\begin{figure}[htbp]
    \centering
    \begin{subfigure}[b]{0.49\textwidth}
        \centering
        \includegraphics[width = 1.0\textwidth]{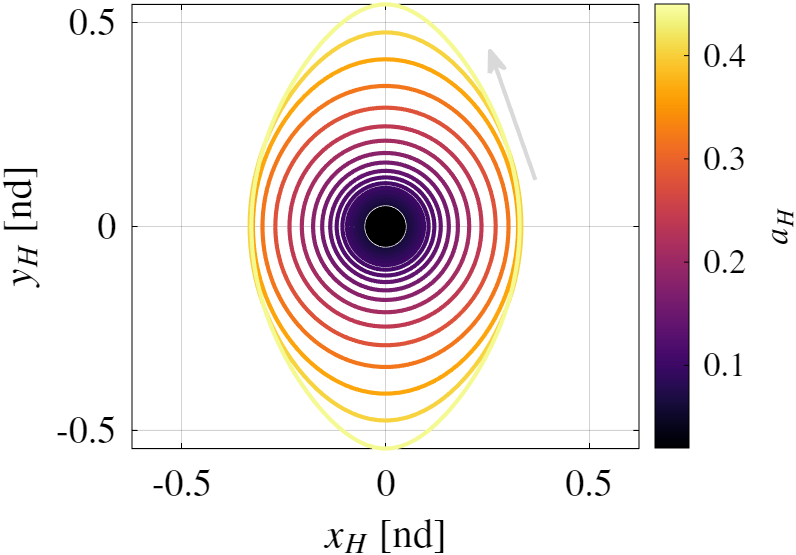}
        \caption{\label{fig:hr3bp_family_g}\acrshort{hr3bp}: family $g$.}
    \end{subfigure}
    \hfill
    \begin{subfigure}[b]{0.49\textwidth}
        \centering
        \includegraphics[width = 1.0\textwidth]{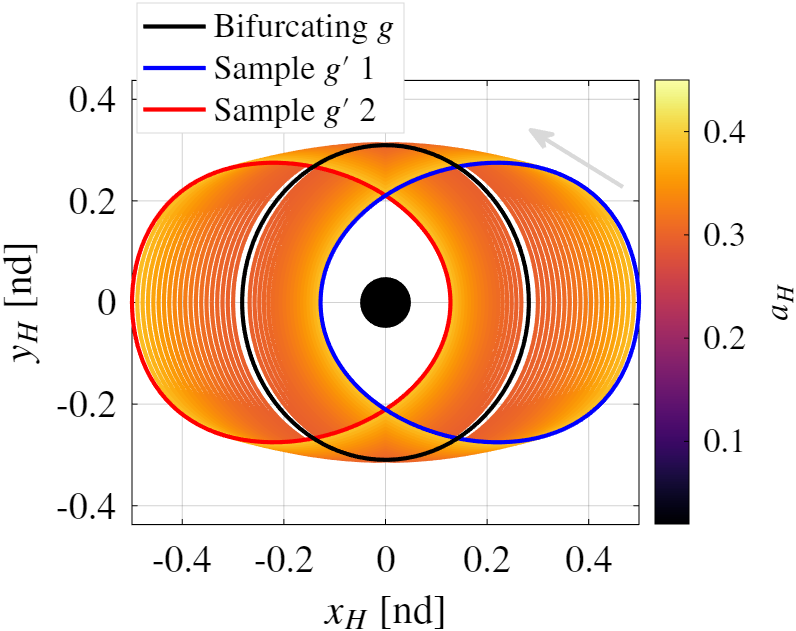}
        \caption{\label{fig:hr3bp_family_gp}\acrshort{hr3bp}: family $g'$.}
    \end{subfigure}
    \begin{subfigure}[b]{0.49\textwidth}
        \centering
        \includegraphics[width = 1.0\textwidth]{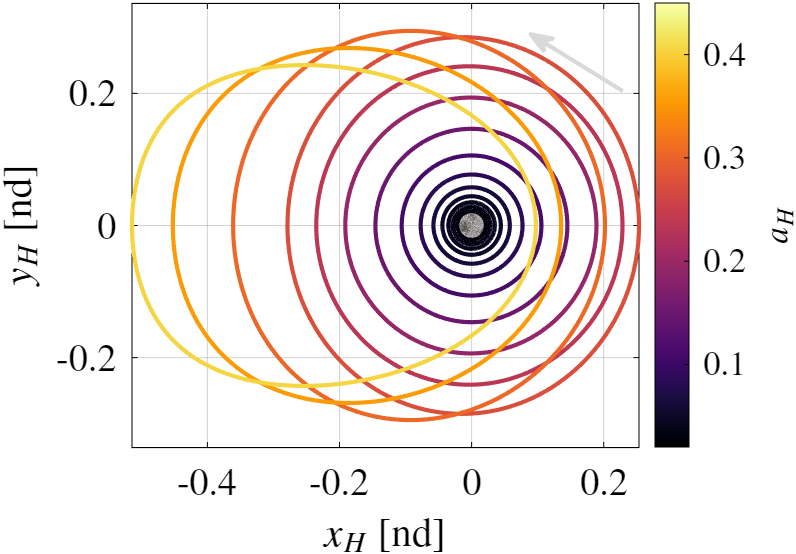}
        \caption{\label{fig:cr3bp_family_lpo}\acrshort{cr3bp}: \acrfull{lpo} family.}
    \end{subfigure}
    \hfill
    \begin{subfigure}[b]{0.49\textwidth}
        \centering
        \includegraphics[width = 1.0\textwidth]{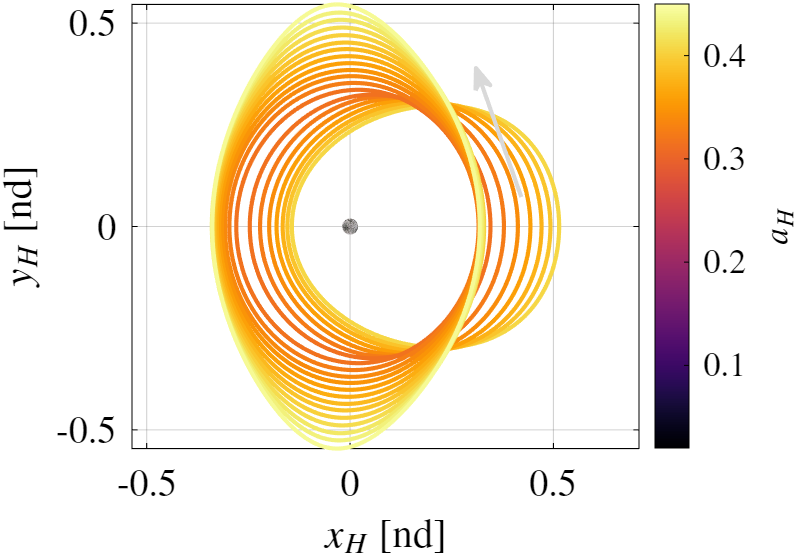}
        \caption{\label{fig:cr3bp_family_dpo}\acrshort{cr3bp}: \acrfull{dpo} family.}
    \end{subfigure}
    \begin{subfigure}[b]{0.49\textwidth}
        \centering
        \includegraphics[width = 1.0\textwidth]{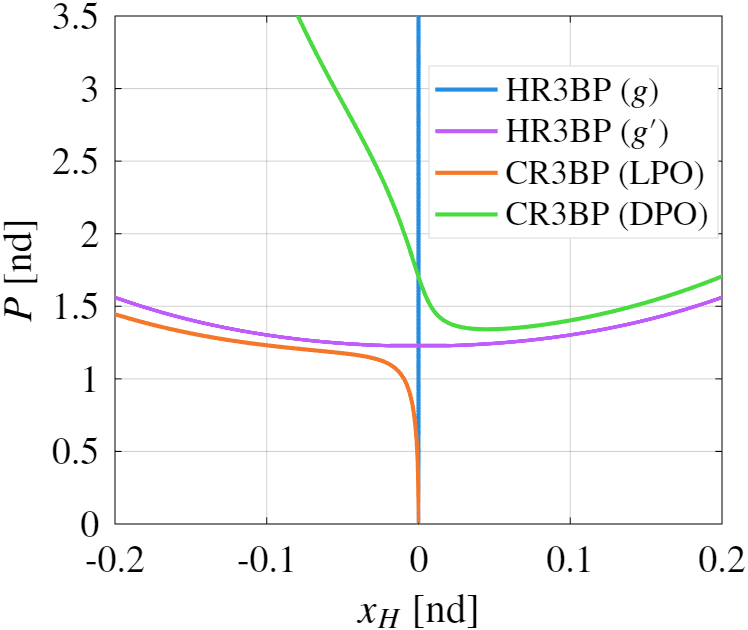}
        \caption{\label{fig:planar_broken}\acrshort{cr3bp}: broken bifurcation between \acrshort{lpo} and \acrshort{dpo} families.}
    \end{subfigure}
    \hfill
    \begin{subfigure}[b]{0.49\textwidth}
        \centering
        \includegraphics[width = 1.0\textwidth]{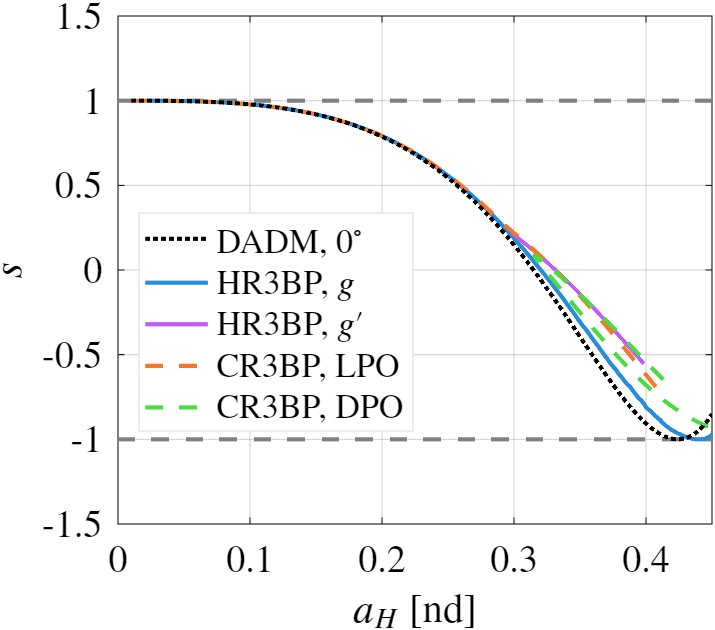}
        \caption{\label{fig:prograde_stability}Stability indices.}
    \end{subfigure}
    \caption{\label{fig:prograde_family} Analysis of circular periodic orbits at $i = 0^\circ$.}    
\end{figure}

The planar prograde limit ($i=0^\circ$) serves as a fundamental baseline, linking the averaged circular equilibria to well-documented periodic orbit families within the \acrshort{r3bps}. These structures are known as family $g$ (following Strömgren's classification \cite{stromgren1933connaissance}) or the \acrfull{lpo} family \cite{guzzetti2016rapid}. Within the \acrshort{hr3bp}, family $g$ is triply symmetric, satisfying $\rho_X, \rho_Y,$ and $\sigma$. Sample periodic orbits from this family are plotted within the \acrshort{hrf} in Fig. \ref{fig:hr3bp_family_g}, colored by their average nondimensional semi-major axis ($a_H$). Since the orbits are non-Keplerian, the averaged distance from $\mathcal{P}_2$ serves as a proxy for $a_H$. At $a_H \ll 1$, the orbit indeed demonstrates circular behavior ($e \approx 0$) and easily correlates to the circular periodic orbit derived from the \acrshort{dadm}. However, as $a_H$ grows, the orbit exhibits an ``eccentric'' behavior, elongating along the $\hat{\bm{y}}$-direction. Within the \acrshort{hr3bp}, family $g$ undergoes a pitchfork bifurcation, giving rise to family $g'$. At this bifurcation, the $\rho_Y$-symmetry is broken, creating two branches of mirror images as plotted in Fig. \ref{fig:hr3bp_family_gp}; these structures correspond to doubly symmetric ($\rho_X, \sigma$) periodic orbits. These orbits evolve with more consistent $a_H$ values but elongate into the $\pm \hat{\bm{x}}$-direction in an asymmetric manner. As family $g'$ does not demonstrate circular behavior, it is not in direct correspondence with the \acrshort{dadm} circular equilibria. However, enumerating these structures and their symmetries is crucial for characterizing the behavior within the \acrshort{cr3bp}. Since the \acrshort{cr3bp} does not possess $\rho_Y$-symmetry, the ideal pitchfork bifurcation observed within the \acrshort{hr3bp} cannot persist \cite{moreno2024bifurcation}. Rather, the dynamics demonstrate a ``broken'' bifurcation, where two families that previously intersected now deform into two separate, non-connecting branches. Specifically, the \acrshort{lpo} family within the \acrshort{cr3bp} appears to merge the members of family $g$ at low $a_H$ values with the ``western'' branch of family $g'$, as illustrated in Fig. \ref{fig:cr3bp_family_lpo}. The remaining members of families $g$ and $g'$ constitute the \acrfull{dpo} family within the \acrshort{cr3bp}, depicted in Fig. \ref{fig:cr3bp_family_dpo}. These planar prograde families are illustrated in the $x_H-P$ space in Fig. \ref{fig:planar_broken}, where $x_H$ is recorded at $\dot{y}_H = 0$. The \acrshort{hr3bp} curves clearly demonstrate a pitchfork bifurcation; the triply symmetric family $g$ appears as a straight line, while family $g'$ manifests as mirror images in the $x_H-P$ space. In contrast, the \acrshort{cr3bp} families reconnect these curves into two disconnected branches, each maintaining a doubly symmetric configuration.

Beyond the geometric resemblance of both family $g$ (\acrshort{hr3bp}) and the \acrshort{lpo} family (\acrshort{cr3bp}) to the circular equilibria at low energies, frequency information serves as a quantitative criterion for connection. Integrating the variational dynamics along the periodic orbits in the \acrshort{r3bps} for a full period yields the monodromy matrix. The non-trivial eigenvalues ($\lambda$) from this matrix contain critical information regarding stability and bifurcations. For planar periodic orbits, the out-of-plane eigenvalues are decoupled from the in-plane components, with eigenvectors directed strictly in the $\hat{\bm{z}}$-direction. Denoting the out-of-plane eigenvalue pair as $\lambda_z$ and $1/\lambda_z$, the stability index is defined as
\begin{align}
s := (\lambda_z + 1/\lambda_z)/2. \label{eq:stability_index}
\end{align} A common scenario is when $|\lambda_z| = 1$, representing a stable oscillation on the linear center manifold. In such cases, $\lambda_z = \exp{(\mathrm{i}\Psi)}$ (where $\mathrm{i}^2 = -1$), and the stability index becomes $s = \cos \Psi$. Here, $\Psi$ represents the rotation angle within the out-of-plane subspace. From the \acrshort{dadm} analysis, an analogous rotation angle is derived based on the relative frequency between the satellite and the rotating frame. Recall that the synodic resonance ratio is defined as $\eta_{0^\circ} = (\nu_s - \nu_m)/\nu_m$ from Eq. \eqref{eq:synodic}. The rotation angle $\Psi_{\text{DADM}, 0^\circ}$ corresponds to the change in the medium-period angle over one full synodic revolution of the satellite, i.e.,
\begin{align}
\Psi_{\text{DADM}, 0^\circ} = \frac{2\pi}{\nu_s - \nu_m} \nu_m = \frac{2\pi}{\eta_{0^\circ}}. \label{eq:psi_dadm}
\end{align}
Consequently, the stability index predicted by the \acrshort{dadm} is defined as,
\begin{align}
s_{\text{DADM}, 0^\circ} = \cos(\Psi_{\text{DADM}, 0^\circ}) = \cos \left( \frac{2\pi}{\eta_{0^\circ}} \right).
\end{align}
These stability indices are compared in Fig. \ref{fig:prograde_stability}, plotted against $a_H$. The geometric behaviors and frequency structures confirm the correspondence between the circular equilibria and the $g$/\acrshort{lpo} families. The correspondence is strong in the $a_H \ll 1$ regime but disparities grow as the unaveraged structures deform into complex, non-circular shapes. Notably, the \acrshort{cr3bp} \acrshort{lpo} continuously evolves into the geometry of family $g'$, while the original circular-like geometry is effectively dominated by the \acrshort{dpo} family. As such, for $a_H \gtrapprox 0.32$ [nd], the circular equilibria from the DADM are more appropriately identified as linking to the \acrshort{dpo} family within the \acrshort{cr3bp}. Note that the DPO family evolves non-monotonically in $a_H$, as evident from Fig.~\ref{fig:cr3bp_family_dpo}. As such, two green curves exist at a given $a_H$ value within Fig.~\ref{fig:prograde_stability}, also correlated with the broken bifurcation in Fig.~\ref{fig:planar_broken}.

\subsubsection{\label{sec:circular_bif_retrograde}Connections at Planar, Retrograde ($i = 180^\circ$)}

\begin{figure}[htbp]
    \centering
    \begin{subfigure}[b]{0.49\textwidth}
        \centering
        \includegraphics[width = 1.0\textwidth]{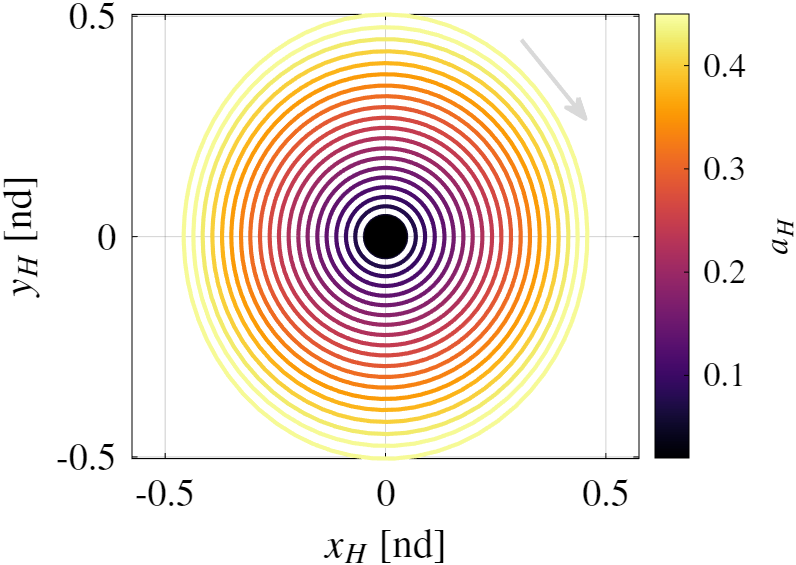}
        \caption{\label{fig:hr3bp_family_f}\acrshort{hr3bp}: family $f$.}
    \end{subfigure}
    \hfill
    \begin{subfigure}[b]{0.49\textwidth}
        \centering
        \includegraphics[width = 1.0\textwidth]{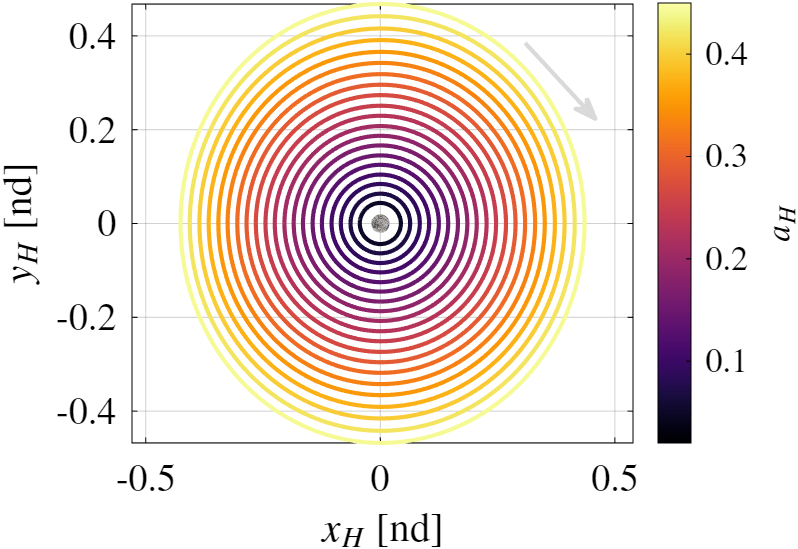}
        \caption{\label{fig:cr3bp_family_dro}\acrshort{cr3bp}: \acrfull{dro} family.}
    \end{subfigure}
    \begin{subfigure}[b]{0.49\textwidth}
        \centering
        \includegraphics[width = 1.0\textwidth]{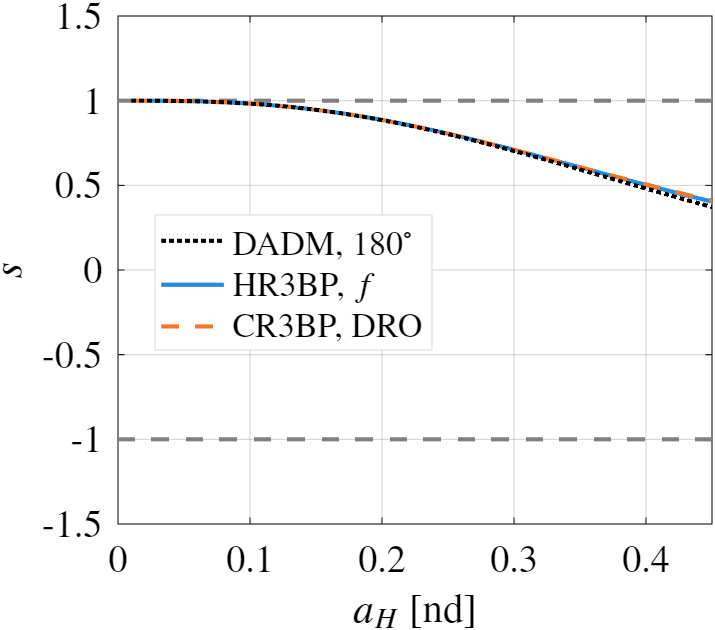}
        \caption{\label{fig:retrograde_stability}Stability indices.}
    \end{subfigure}
    \caption{\label{fig:retrograde_family}Analysis of circular periodic orbits at $i = 180^\circ$.}    
\end{figure}

The planar retrograde limit ($i = 180^\circ$) serves as the counterpart to the prograde case, linking the retrograde circular equilibria to stable retrograde families in the \acrshort{r3bps}. In this limit, the orbits correspond to Strömgren's family $f$~\cite{stromgren1933connaissance} or the \acrfull{dro} family \cite{guzzetti2016rapid}. Figures \ref{fig:hr3bp_family_f} and \ref{fig:cr3bp_family_dro} depict these families in the \acrshort{hr3bp} and \acrshort{cr3bp}, respectively, demonstrating qualitatively similar behaviors. A key distinction lies in the evolution of symmetry. While family $f$ in the \acrshort{hr3bp} retains triple symmetry ($\sigma, \rho_X, \rho_Y$), the \acrshort{dro} family within the \acrshort{cr3bp} is reduced to a doubly symmetric configuration ($\sigma, \rho_X$) due to the inherent loss of $\rho_Y$-symmetry. However, in contrast to the prograde case, family $f$ does not undergo a symmetry-breaking pitchfork bifurcation within the domain of interest. Consequently, there exists a direct, continuous correspondence between family $f$ in the \acrshort{hr3bp} and the \acrshort{dro} family in the \acrshort{cr3bp}, essentially avoiding the complexity of broken bifurcations encountered in the prograde analysis. 

The frequency analysis follows the framework established in the preceding section, with the kinematic modification that the retrograde motion increases the synodic frequency ($\nu_s + \nu_m$). Accordingly, the rotation angle $\Psi_{\text{DADM}, 180^\circ}$ is derived as,
\begin{align}
    \Psi_{\text{DADM},180^\circ} = \frac{2\pi}{\nu_s + \nu_m} \nu_m = \frac{2\pi}{\eta_{180^\circ}}.
\end{align}
The corresponding stability index is defined as,
\begin{align}
    s_{\text{DADM},180^\circ} = \cos(\Psi_{\text{DADM},180^\circ}) =  \cos \left( \frac{2\pi}{\eta_{180^\circ}} \right).
\end{align}
These analytical stability indices are compared with the numerical results in Fig. \ref{fig:retrograde_stability}. In contrast to the prograde case, where geometric deformations lead to a gradual mismatch as $a_H$ increases, the retrograde families exhibit remarkable agreement between the \acrshort{dadm} predictions and the unaveraged dynamics. 

\subsubsection{Period-Multiplying Bifurcations and Number of Branches}

Building on the planar-limit connections in Sections~\ref{sec:circular_bif_prograde}--\ref{sec:circular_bif_retrograde}, this subsection characterizes how the \emph{spatial} circular PO families branch from the planar $g$/LPO and $f$/DRO families. Owing to the circular geometry ($e=0$), these spatial families emerge through period-multiplying bifurcations; the multiplying factor at each planar limit and the resulting number of branches are established below. It is noted that, due to the circular geometry ($e=0$), the $p \mp q$ loops within the rotating frame are geometrically identical, reducing the fundamental period of the orbit by a factor of $p \mp q$. This modulation signifies a period-$(p \mp q)$-multiplying bifurcation from the respective planar PO families within the \acrshort{r3bps}. From these bifurcations, Vertical Self-Resonant (VSR) orbits are defined by \citet{robin1980numerical}, further investigated by \citet{aydin2025studying}, \citet{aydin2025exploration}, as well as \citet{peng2025analog}. Focusing on the case where $q=1$, \citet{aydin2025exploration} demonstrates that families emerging from period-$(p+1)$-multiplying bifurcations from the DRO family generally connect to period-$(p-1)$-multiplying bifurcations within the LPO families. The current work supplies a complementary viewpoint by explicitly considering the DADM dynamics to arrive at a similar conclusion; VSR orbits are spatial symmetric POs that originate from the \acrshort{dadm} circular equilibria. 

From the current analysis, a general statement is synthesized as follows: ``nominally, symmetric periodic orbits originating from the circular equilibria with a resonance ratio $p:q$ connect to $g$/LPO families via a period-$(p-q)$-multiplying bifurcation at $i=0^{\circ}$, and connect to $f$/DRO families via a period-$(p+q)$-multiplying bifurcation at $i=180^{\circ}$,'' a principle that extends to cases where $q \neq 1$. The qualifier ``nominally'' is deliberately employed to indicate that this relation represents the baseline expectation derived from the integrable \acrshort{dadm} rather than the absolute ground truth of the non-integrable \acrshort{r3bps}. As visualized in Fig. \ref{fig:prograde_family}, the unaveraged structures inevitably deviate from the averaged approximations as the effective $a_H$ increases. Furthermore, the global topology is susceptible to instabilities arising from nearby hyperbolic invariant manifolds within the unaveraged dynamics, altering the bifurcation structures. Consequently, the universality of this statement is subject to the characteristics of the specific resonance ratio and the governing dynamics represented by the mass ratio parameter $\mu \geq 0$. Still, the establishment of such an \textit{a priori} criterion remains invaluable; it provides a systematic roadmap to predict and organize the complex web of periodic families within the \acrshort{r3bps}.

The number of branches and the order of symmetries originating from the out-of-plane period-multiplying bifurcations are established in \citet{aydin2025studying}. The solution multiplicity laid out in Table \ref{tab:symmetry_evolution} is consistent with their findings and signifies the number of family branches. For instance, the HR3BP renders two configurations for each double symmetry pairing in Table \ref{tab:symmetry_evolution}. These distinct geometries connect smoothly across the period-multiplying bifurcation and may be evaluated as the same branch. For example, two OX/XOZ and OY/YOZ doubly symmetric branches emerge from a period-$(p+q)$-multiplying bifurcation along family $f$ when $p:q$ are in odd-odd parity. The two OX/XOZ branches exhibit a mirror configuration across the $yz$-plane, coinciding at $i = 0^\circ$ and $180^\circ$. In the current investigation, they are differentiated as distinct branches as a distinct geometry appears (in terms of the apsidal sequence) at the same level of inclination. Thus, the solution multiplicity in Table \ref{tab:symmetry_evolution} also depicts the number of solution branches from the period-multiplying bifurcations at the planar limits. The branches evolve as a function of parity in $p:q$ resonances as well as the dynamics, i.e., \acrshort{hr3bp} and \acrshort{cr3bp}.

\subsection{\label{sec:frozen_bif}Families of Frozen Periodic Orbits}

While the circular orbit families generally evolve continuously across the inclination range, the \acrshort{dadm} analysis predicts a singularity as the inclination approaches the polar region ($i \to 90^\circ$) for a fixed resonance ratio $\eta = p/q$ (see Fig. \ref{fig:a_i_eta_F_hr3bp}). Consequently, the domain of existence for frozen orbits is restricted to specific bands, approximately bounded by $39.2^\circ \lessapprox i < 90^\circ$ and $90^\circ < i \lessapprox 140.8^\circ$. The subsequent sections investigate the dynamical behaviors at these limiting boundaries: the asymptotic approach toward the pole ($i \rightarrow 90^\circ$) and the bifurcation points at the critical inclination ($i_{\mathrm{crit}} \approx 39.2^\circ, 140.8^\circ$).

\subsubsection{Asymptotic behaviors at $ i = 90^\circ$}

\begin{figure}[htbp]
    \centering
    \begin{subfigure}[b]{0.49\textwidth}
        \centering
        \includegraphics[width = 1.0\textwidth]{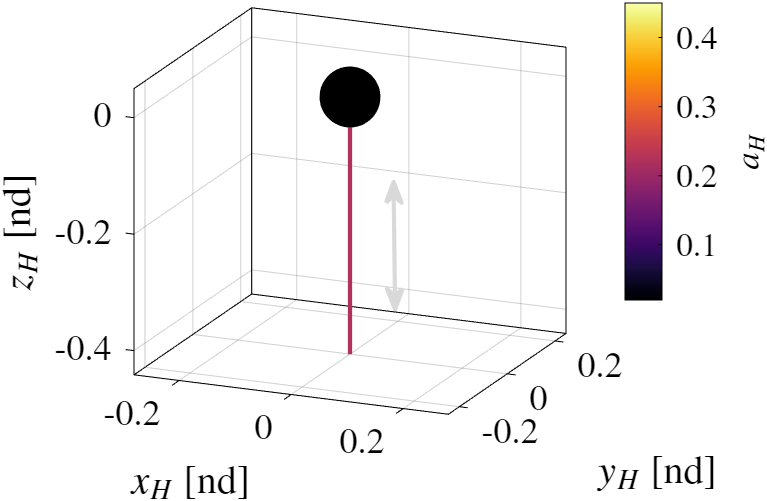}
        \caption{\label{fig:hr3bp_family_collision}\acrshort{hr3bp}: a southern collision orbit at $a_H = 0.2$ [nd].}
    \end{subfigure}
    \hfill
    \begin{subfigure}[b]{0.49\textwidth}
        \centering
        \includegraphics[width = 1.0\textwidth]{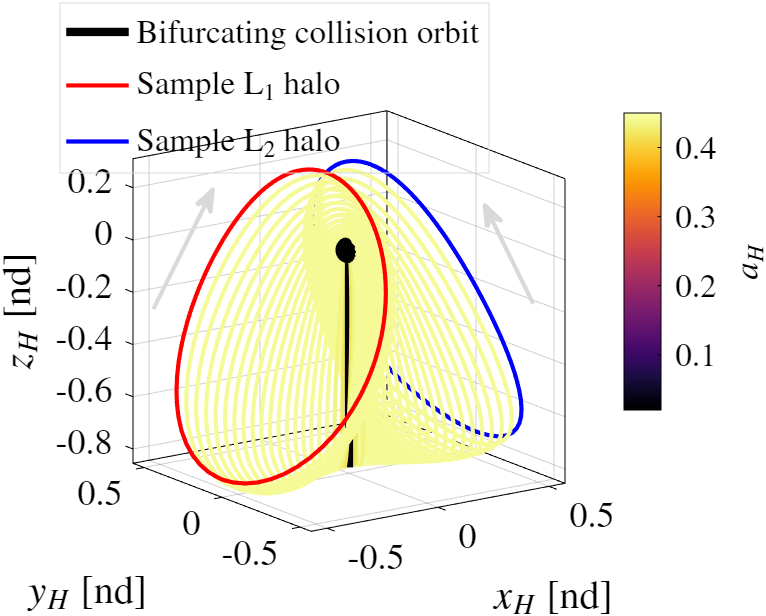}
        \caption{\label{fig:hr3bp_family_halo}\acrshort{hr3bp}: halo families.}
    \end{subfigure}
    \begin{subfigure}[b]{0.49\textwidth}
        \centering
        \includegraphics[width = 1.0\textwidth]{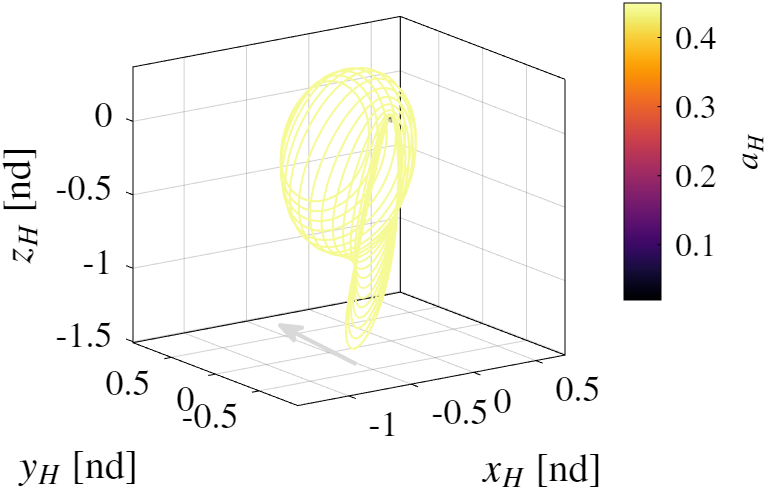}
        \caption{\label{fig:cr3bp_family_l1}\acrshort{cr3bp}: L$_1$ halo family.}
    \end{subfigure}
    \hfill
    \begin{subfigure}[b]{0.49\textwidth}
        \centering
        \includegraphics[width = 1.0\textwidth]{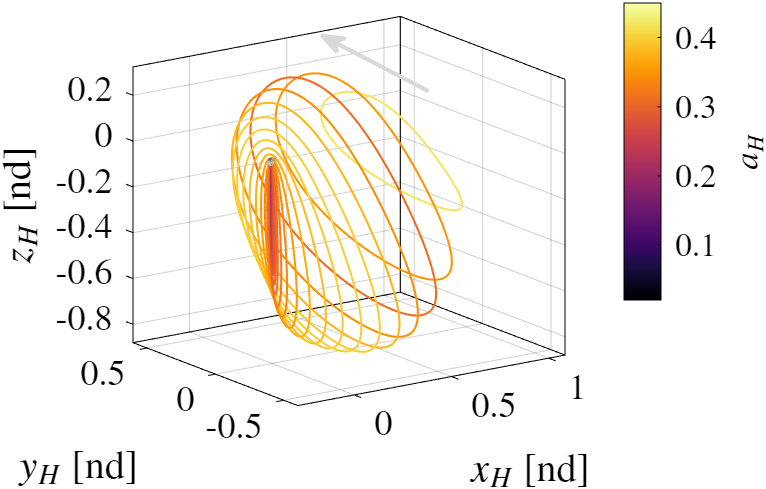}
        \caption{\label{fig:cr3bp_family_l2}\acrshort{cr3bp}: L$_2$ halo family.}
    \end{subfigure}
    \begin{subfigure}[b]{0.49\textwidth}
        \centering
        \includegraphics[width = 1.0\textwidth]{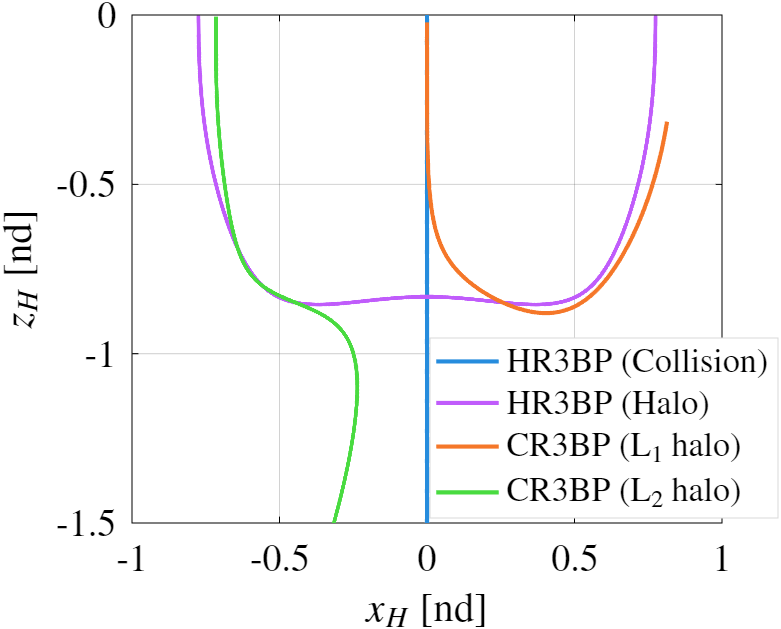}
        \caption{\label{fig:halos_broken}\acrshort{cr3bp}: broken bifurcation between L$_1$ and L$_2$ halo families.}
    \end{subfigure}
    \hfill
    \begin{subfigure}[b]{0.49\textwidth}
        \centering
        \includegraphics[width = 1.0\textwidth]{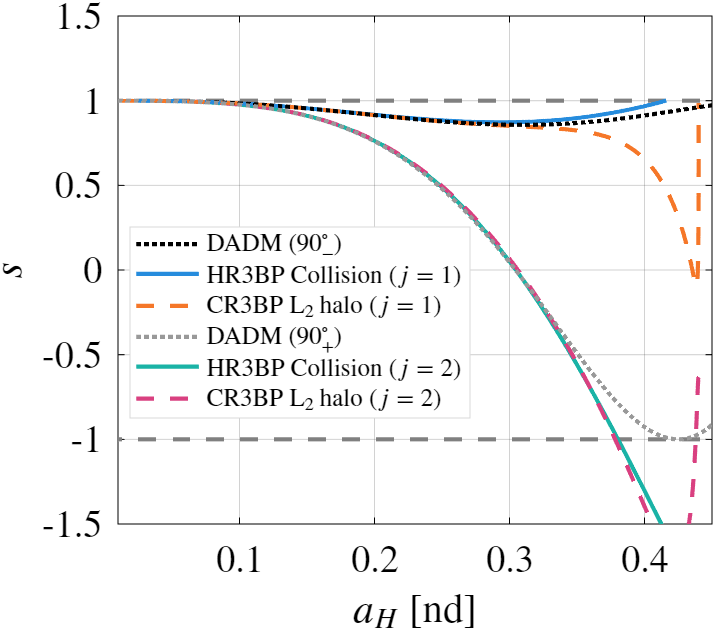}
        \caption{\label{fig:frozen_stability}Stability indices for the frozen periodic orbits near $i = 90^\circ$.}
    \end{subfigure}
    \caption{\label{fig:frozen_family}Analysis of frozen periodic orbits near $i = 90^\circ$.}    
\end{figure}

The analytical expression for the medium-period frequency (Eq. \eqref{eq:nu_m_frozen}) exhibits a discontinuity due to the $\text{sgn}(\cos i)$ term. This results in distinct asymptotic limits for the evolution of the prograde and retrograde families as
\begin{align}
    \nu_{m, 90^\circ_-} := \lim_{i \to 90^\circ_-} \nu_m &= 1 + \frac{15}{4\nu_s}(1-\mu)\sqrt{\frac{3}{5}}\label{eq:num_limit_pro} \\
    \nu_{m, 90^\circ_+} := \lim_{i \to 90^\circ_+} \nu_m &= 1 - \frac{15}{4\nu_s}(1-\mu)\sqrt{\frac{3}{5}}. \label{eq:num_limit_retro}
\end{align}
This discontinuity implies that frozen periodic orbits do not evolve continuously through the pole. Rather, the prograde and retrograde branches converge toward distinct frequency limits as they approach a rectilinear geometry ($e \to 1$). Such rectilinear trajectories correspond to ``collision'' orbits within the \acrshort{hr3bp}. At the initial epoch ($t=0$), these orbits depart from $z_H \neq 0$ with all other state components at zero; at the half-period ($t = P/2$), the state coincides with the center of $\mathcal{P}_2$. To construct these orbits, the Kustaanheimo-Stiefel (KS) regularization scheme is employed, adapted from \citet{howell1984almost}. These orbits are doubly symmetric, possessing symmetries across the $xz$- and $yz$-planes. While they do not satisfy the $\sigma$-symmetry, a mirrored family exists across the $xy$-plane, distinguishing the northern and southern geometries. A representative southern collision PO is depicted in Fig. \ref{fig:hr3bp_family_collision}, initiated with $z_H = -0.4$ [nd]. This family persists continuously across a range of initial $z_H$ values, where the effective semi-major axis is defined as $a_H = |z_H|/2$ at the apoapsis. Within the \acrshort{hr3bp}, this collision family evolves into halo orbits around the L$_1$ and L$_2$ Lagrange points via a pitchfork bifurcation, as plotted in Fig. \ref{fig:hr3bp_family_halo}. However, since the \acrshort{cr3bp} lacks $yz$-plane symmetry, this pitchfork bifurcation cannot exist as a true intersection. Instead, the structure degenerates into two disconnected branches corresponding to the L$_1$ and L$_2$ halo orbits, as depicted in Figs. \ref{fig:cr3bp_family_l1} and \ref{fig:cr3bp_family_l2}. In Fig. \ref{fig:halos_broken}, these families are projected onto the $x_H - z_H$ space, recorded at $\dot{z}_H = 0$. Similar to the planar case (Fig. \ref{fig:planar_broken}), the \acrshort{hr3bp} exhibits a clear pitchfork bifurcation connecting doubly and singly symmetric PO families, whereas the \acrshort{cr3bp} reveals two separated singly symmetric branches. For a detailed analysis of the collision-halo evolution across varying $\mu$ values, refer to \citet{joung2025bifurcations}. It is noted that the \acrshort{hr3bp} collision orbits with smaller apoapsis radii transition specifically to the L$_2$ halo orbits in the \acrshort{cr3bp}. For the halo orbits, the semi-major axis is approximated by the mean of the minimum and maximum $\hat{\bm{z}}$-excursions along the orbit.

To quantify the correspondence between the models for the frozen \acrshort{pos}, the stability of the spatial structures in the \acrshort{r3bps} is examined. Stability indices are defined as,
\begin{align}
    s_{j} = (\lambda_j + 1/\lambda_j)/2,
\end{align}
analogous to Eq. \eqref{eq:stability_index}, where $j = 1,2$ denotes the two distinct non-trivial eigenvalue pairs. Corresponding quantities for the \acrshort{dadm} limits are defined as:
\begin{align}
    s_{\text{DADM}, 90^\circ_{\pm}} = \cos \left( \frac{2\pi}{\eta_{90^\circ_{\pm}}}\right), \quad \eta_{90^\circ_{\pm}} = \frac{\nu_s}{\nu_{m, 90^\circ_{\pm}}},
\end{align}
utilizing the limits from Eqs. \eqref{eq:num_limit_pro}--\eqref{eq:num_limit_retro}. Figure \ref{fig:frozen_stability} juxtaposes these \acrshort{dadm} stability indices at $i\rightarrow 90^\circ$ with those of the collision orbits (\acrshort{hr3bp}) and L$_2$ halo orbits (\acrshort{cr3bp}). Eigenvalues for the collision orbits are constructed via the KS-regularized monodromy matrix, adapted from \citet{howell1984almost}. The comparison reveals a notable correlation: each stability index $s_{j}$ ($j = 1, 2$) for the \acrshort{r3bp} structures aligns with the distinct \acrshort{dadm} limits from the prograde and retrograde sides. For the \acrshort{r3bp} structures, the plot focuses on the range where $s_2 < 1$. For the retrograde limit, the \acrshort{hr3bp} demonstrates a stronger correspondence than the \acrshort{cr3bp}. While the \acrshort{cr3bp} L$_2$ halos resemble the \acrshort{hr3bp} collision orbits at $a_H \ll 1$, they gradually deform into the specific L$_2$ halo geometry, causing the disparity to grow rapidly. Additionally, $s_1 = -1$ marks a period-doubling bifurcation at $a_H \approx 0.38$ [nd], beyond which the unaveraged structures lose stability. Overall, the frequency-based stability analysis successfully bridges the \acrshort{dadm} frozen equilibria limit ($i \rightarrow 90^\circ$) with the collision and halo structures in the unaveraged dynamics.

Similarly to the circular POs, the frozen POs connect to the known \acrshort{r3bp} structures via period-multiplying bifurcations. Specifically, a general statement is synthesized as follows: ``nominally, symmetric periodic orbits originating from the frozen equilibria with a resonance ratio $p:q$ connect to collision (\acrshort{hr3bp}) and L$_2$ halo orbits (\acrshort{cr3bp}) via a period-$p$-multiplying bifurcation as $i \rightarrow 90^\circ$. Two distinct asymptotic behaviors are observed for prograde ($i\rightarrow 90^\circ_{-}$) and retrograde ($i\rightarrow 90^\circ_{+}$) cases that connect to the distinct eigenvalue pairs.'' While symmetric periodic orbits resembling frozen equilibria have been detected emerging from these \acrshort{r3bp} structures via period-multiplying bifurcations \cite{koblick2025novel, aydin2025studying, zhang2025time}, prior works primarily focus on the prograde scenario and do not explicitly articulate the connection to the frozen equilibria from the \acrshort{dadm}. The current work is extended to the retrograde orbits and establishes an explicit correspondence to the \acrshort{dadm} structures; consequently, the same geometries may be identified via period-multiplying bifurcations or, alternatively, through the initialization process detailed in Section \ref{sec:initialization_targeting}, supplying an adaptive framework. The entries in Table \ref{tab:symmetry_evolution} reflect the number of distinct geometries for each symmetry type for the frozen orbits at a given inclination. This count coincides with the number of family branches emerging after the period-multiplying bifurcation from the \acrshort{r3bp} structures, specifically the collision (\acrshort{hr3bp}, detailed in \citet{aydin2025studying}) and L$_2$ halo (\acrshort{cr3bp}) orbits. Note that the southern and northern frozen POs bifurcate from the southern and northern collision/L$_2$ halo orbits, respectively. 

In evaluating these general connections, the retrograde regime requires careful examination. Specifically, the asymptotic connections to the polar limit ($i \rightarrow 90^\circ_{+}$) are realized only when the resonance ratio satisfies $\eta = p/q \geq 12$, as qualitatively indicated by the contours in Fig. \ref{fig:a_i_eta_F_hr3bp}. Below this threshold, the resonant contours do not reach $i = 90^\circ$ and are expected to terminate at locations distinct from the collision and L$_2$ halo orbits. While a comprehensive analysis remains a subject for future studies, preliminary numerical investigations suggest that these families often connect to other branches of \emph{prograde} frozen periodic orbits via period-doubling bifurcations, e.g., retrograde $10:1$ to prograde $5:1$, rather than reaching the asymptotic limit at $i = 90^\circ$.

\subsubsection{\label{sec:frozen_bif_broken}Bifurcation from Circular Periodic Orbits at $i = i_{\mathrm{crit}}$}

The bifurcation of frozen periodic orbits from circular orbits, predicted at $i = i_{\mathrm{crit}}$ within the \acrshort{dadm}, warrants careful examination within the context of unaveraged dynamics. In the integrable \acrshort{dadm}, frozen equilibria branch out from circular equilibria via a perfect pitchfork bifurcation. This ideal bifurcation is supported by a specific time-reversal symmetry inherent to the reduced averaged space ($e, i, \omega$), defined as,
\begin{align}
    \rho_{h}: (k, h) \to (k, -h).
\end{align}
As visualized in the Lidov diagram (Fig. \ref{fig:2d_surface_constants_averaged_stable_unstable}), circular equilibria satisfy this symmetry ($h=0$), whereas frozen equilibria break it, creating two distinct branches corresponding to the southern ($\omega = 90^\circ$) and northern ($\omega = -90^\circ$) families.

In the \acrshort{r3bps}, however, a direct analog to the $\rho_h$-symmetry does not generally exist. Rather, the specific symmetry that frozen POs may break depends on the resonance parity ($p:q$), as detailed in Table \ref{tab:symmetry_evolution}. Within the \acrshort{hr3bp}, for odd-odd and odd-even parities, the frozen POs break an axial symmetry (OX or OY) that exists in the circular POs, preserving the pitchfork bifurcation structure. However, the even-odd parity presents a degenerate case where the set of symmetries is identical between the circular (XOZ/YOZ) and frozen (XOZ/YOZ) POs. This degeneracy is exacerbated in the \acrshort{cr3bp} due to the reduced symmetry group. As apparent in Table \ref{tab:symmetry_evolution}, both odd-even and even-odd parities result in an identical symmetry order for circular and frozen POs within the \acrshort{cr3bp}. In such cases, the mechanism for a symmetry-breaking bifurcation is absent. Consequently, the $\rho_h$-symmetry-driven pitchfork bifurcation of the \acrshort{dadm} manifests as a "broken" bifurcation in the unaveraged dynamics, analogous to the behaviors observed in the planar (Fig. \ref{fig:planar_broken}) and halo (Fig. \ref{fig:halos_broken}) families.

\begin{figure}[htbp]
    \centering
    \begin{subfigure}[b]{0.49\textwidth}
        \centering
        \includegraphics[width = 1.0\textwidth]{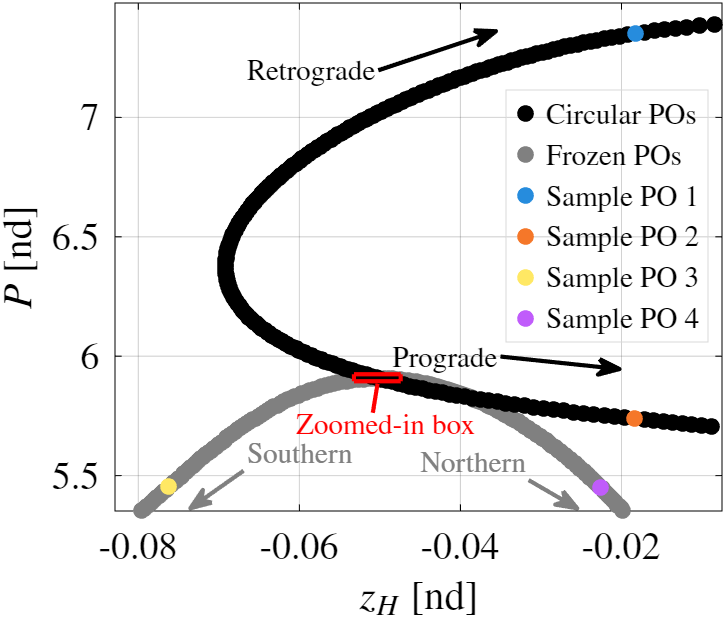}
        \caption{\label{fig:circular_frozen_broken_zoom_out}Families of circular and frozen POs (states recorded at $\xp{southern}$, the southern apex at $x > 0$).}
    \end{subfigure}
    \hfill
    \begin{subfigure}[b]{0.49\textwidth}
        \centering
        \includegraphics[width = 1.0\textwidth]{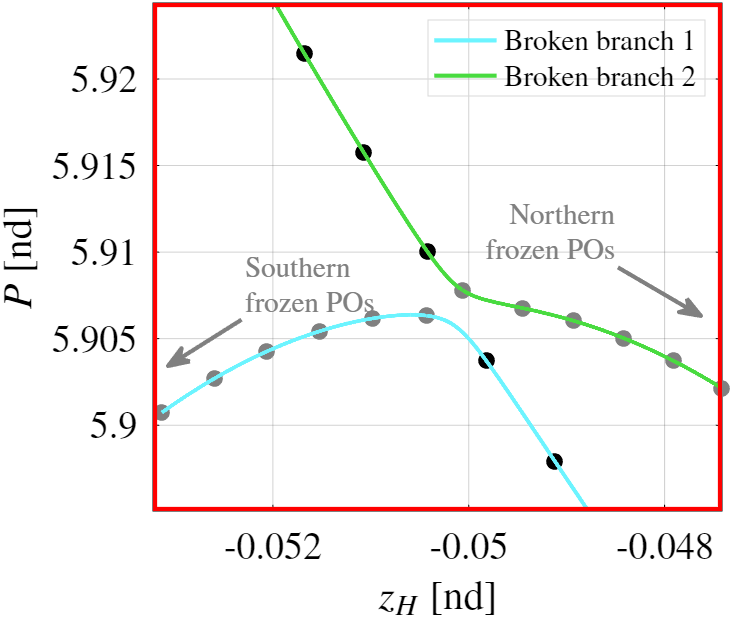}
        \caption{\label{fig:circular_frozen_broken_zoom_in}Zoom-in view of the red box in Fig. \ref{fig:circular_frozen_broken_zoom_out}.}
    \end{subfigure}
    \begin{subfigure}[b]{0.49\textwidth}
        \centering
        \includegraphics[width = 1.0\textwidth]{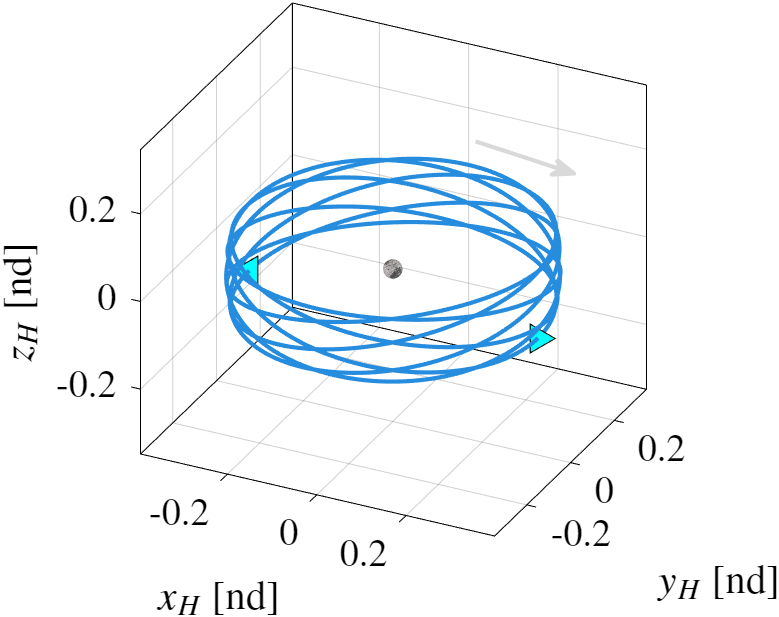}
        \caption{\label{fig:circular_frozen_broken_sample1}Sample PO 1 (circular) from Fig. \ref{fig:circular_frozen_broken_zoom_out}.}
    \end{subfigure}
    \hfill
    \begin{subfigure}[b]{0.49\textwidth}
        \centering
        \includegraphics[width = 1.0\textwidth]{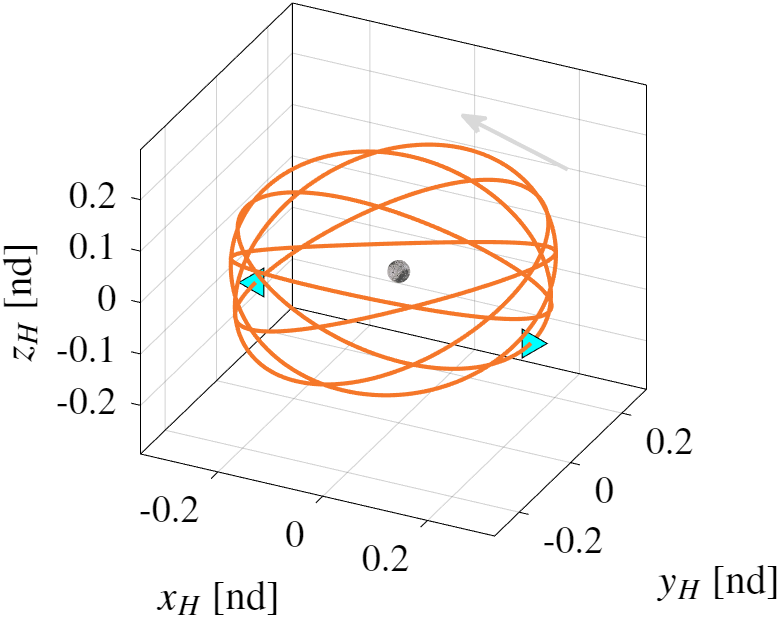}
        \caption{\label{fig:circular_frozen_broken_sample2}Sample PO 2 (circular) from Fig. \ref{fig:circular_frozen_broken_zoom_out}.}
    \end{subfigure}
    \begin{subfigure}[b]{0.49\textwidth}
        \centering
        \includegraphics[width = 1.0\textwidth]{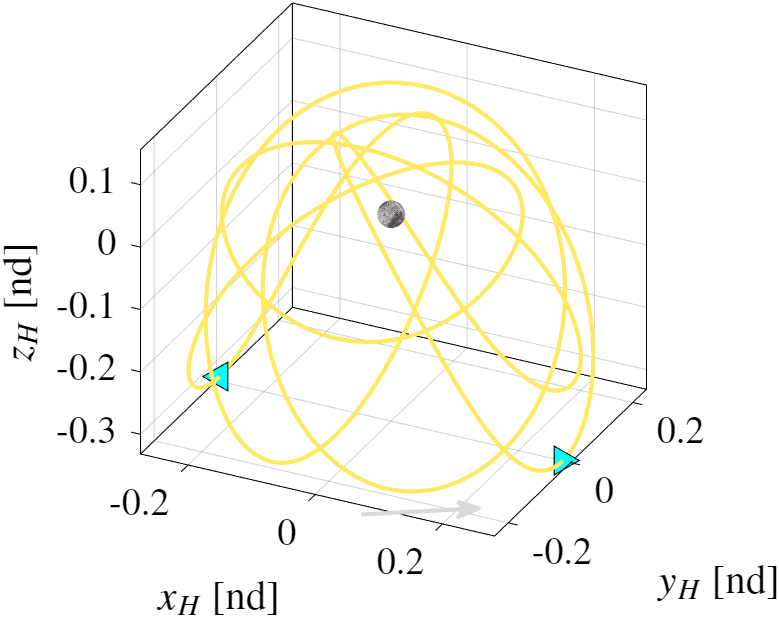}
        \caption{\label{fig:circular_frozen_broken_sample3}Sample PO 3 (frozen) from Fig. \ref{fig:circular_frozen_broken_zoom_out}.}
    \end{subfigure}
    \hfill
    \begin{subfigure}[b]{0.49\textwidth}
        \centering
        \includegraphics[width = 1.0\textwidth]{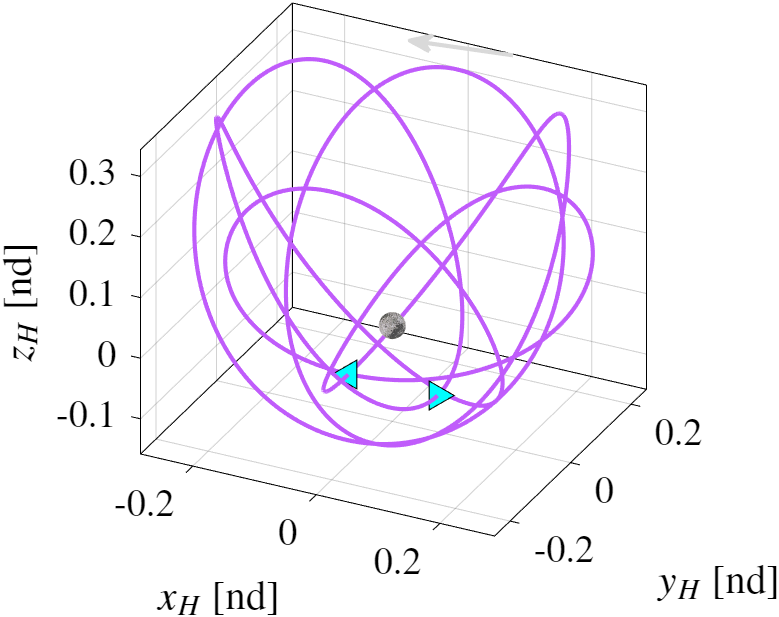}
        \caption{\label{fig:circular_frozen_broken_sample4}Sample PO 4 (frozen) from Fig. \ref{fig:circular_frozen_broken_zoom_out}.}
    \end{subfigure}
    \caption{\label{fig:circular_frozen_broken}Broken bifurcation between circular and frozen POs for $p = 6$, $q = 1$, \acrshort{cr3bp}.}    
\end{figure}

To illustrate the existence of broken bifurcations within the \acrshort{r3bps}, consider the $p = 6, q = 1$ resonance within the \acrshort{cr3bp}. This case is predicted to exhibit a broken bifurcation as the symmetry order for both circular and frozen POs is consistent at 1 (XOZ). In Figure \ref{fig:circular_frozen_broken_zoom_out}, the families of circular POs (black markers) and frozen POs (grey markers) in the $z_H - P$ space appear. Among the multiple XOZ-symmetric configurations, the plot focuses on geometries characterized by the half-period apsidal sequence $(\ID{A}{+x}, \ID{A}{-x})$, or equivalently, $(\xp{southern}, \xm{southern})$. Then, the $z_H$ coordinate at $\ID{A}{+x} (\xp{southern})$ is recorded. The focus is specifically on the transition region from the prograde circular PO to the prograde frozen PO boxed in red within Fig. \ref{fig:circular_frozen_broken_zoom_out}. Upon zooming in (Fig. \ref{fig:circular_frozen_broken_zoom_in}) and densely sampling members along each branch, it is evident that two distinct, non-intersecting branches exist. The southern and northern frozen POs each connect as separate branches to the prograde circular family, confirming the broken nature of the bifurcation. Although not supplied here, a similar broken bifurcation occurs on the retrograde side as well. 

These broken bifurcations present a unique challenge in mapping the web of symmetric POs; following a continuous branch may lead to a qualitative shift in orbital behavior; in this case, transitioning between frozen and circular characteristics. Standard stability analysis often fails to predict the existence of the disconnected branch, as typically one branch exhibits a fold (saddle-node) bifurcation while the other exhibits no stability change. Without \textit{a priori} knowledge of the family connections, detecting the secondary branch is extremely challenging\footnote{See Fig. 25 in \citet{aydin2025studying}, where the $6:1$ circular periodic orbits within the \acrshort{hr3bp} connect to the collision orbits; while fold bifurcations are recorded in their work, the current analysis predicts that broken bifurcations occur between the circular and frozen POs with the same number of symmetries (XOZ/YOZ).}. The current analysis bridges this gap by explicitly linking the averaged dynamics to the \acrshort{r3bps}. Through symmetry analysis, the specific $p:q$ parities susceptible to broken bifurcations are identified across the \acrshort{hr3bp} and \acrshort{cr3bp}, enabling the systematic tracing of distinct, disconnected branches utilizing global insights derived from the \acrshort{dadm} equilibria.

\section{\label{sec:archetype}Archetypical Bifurcation Diagrams within Unaveraged Dynamics}

\begin{figure}[htbp]
    \centering
    \begin{subfigure}[b]{0.49\textwidth}
        \centering
        \includegraphics[width = 0.8\textwidth]{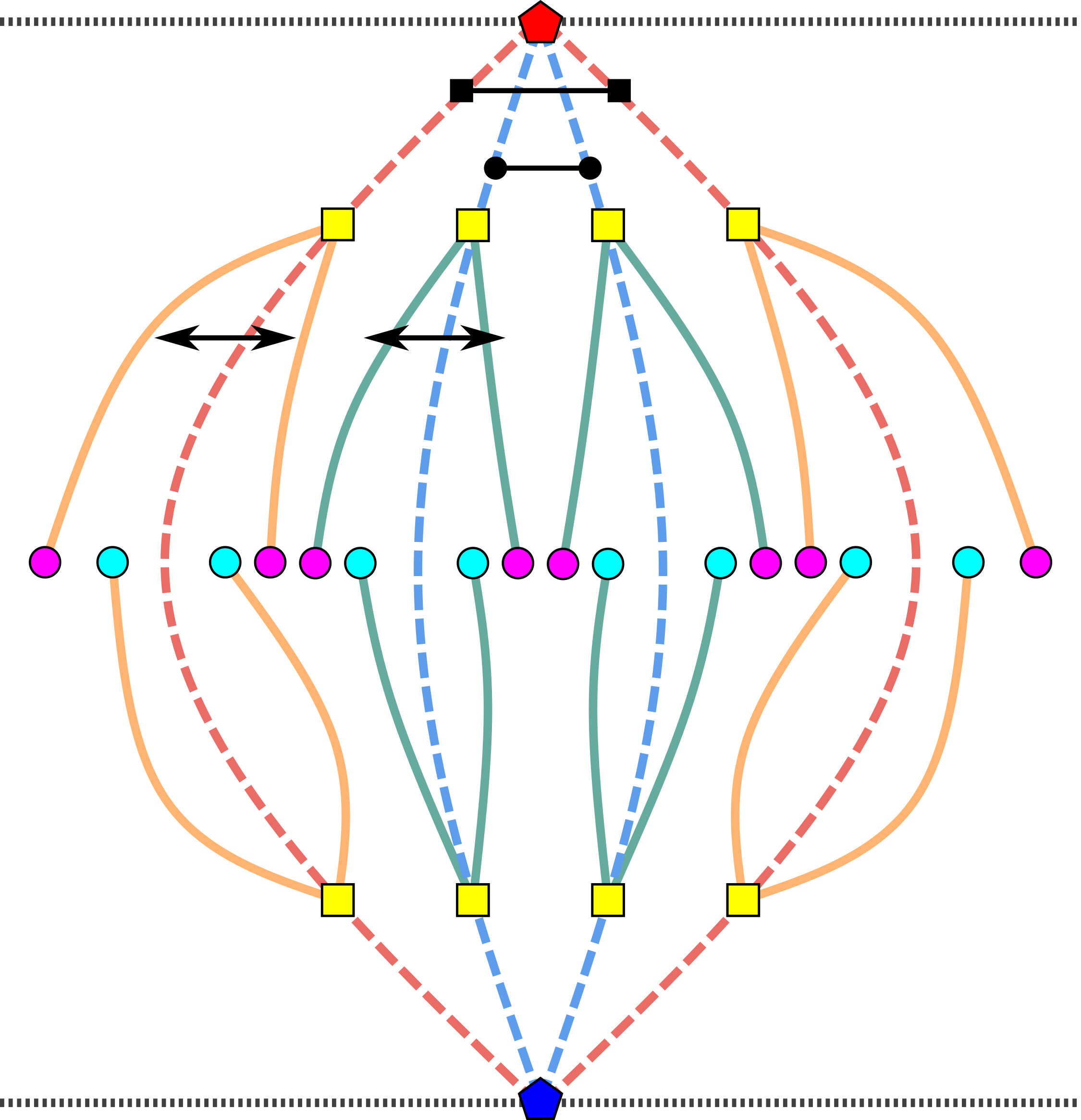}
        \caption{\label{fig:hr3bp_odd_odd}\acrshort{hr3bp}: odd-odd $p:q$.}
    \end{subfigure}
    \hfill
    \begin{subfigure}[b]{0.49\textwidth}
        \centering
        \includegraphics[width = 0.8\textwidth]{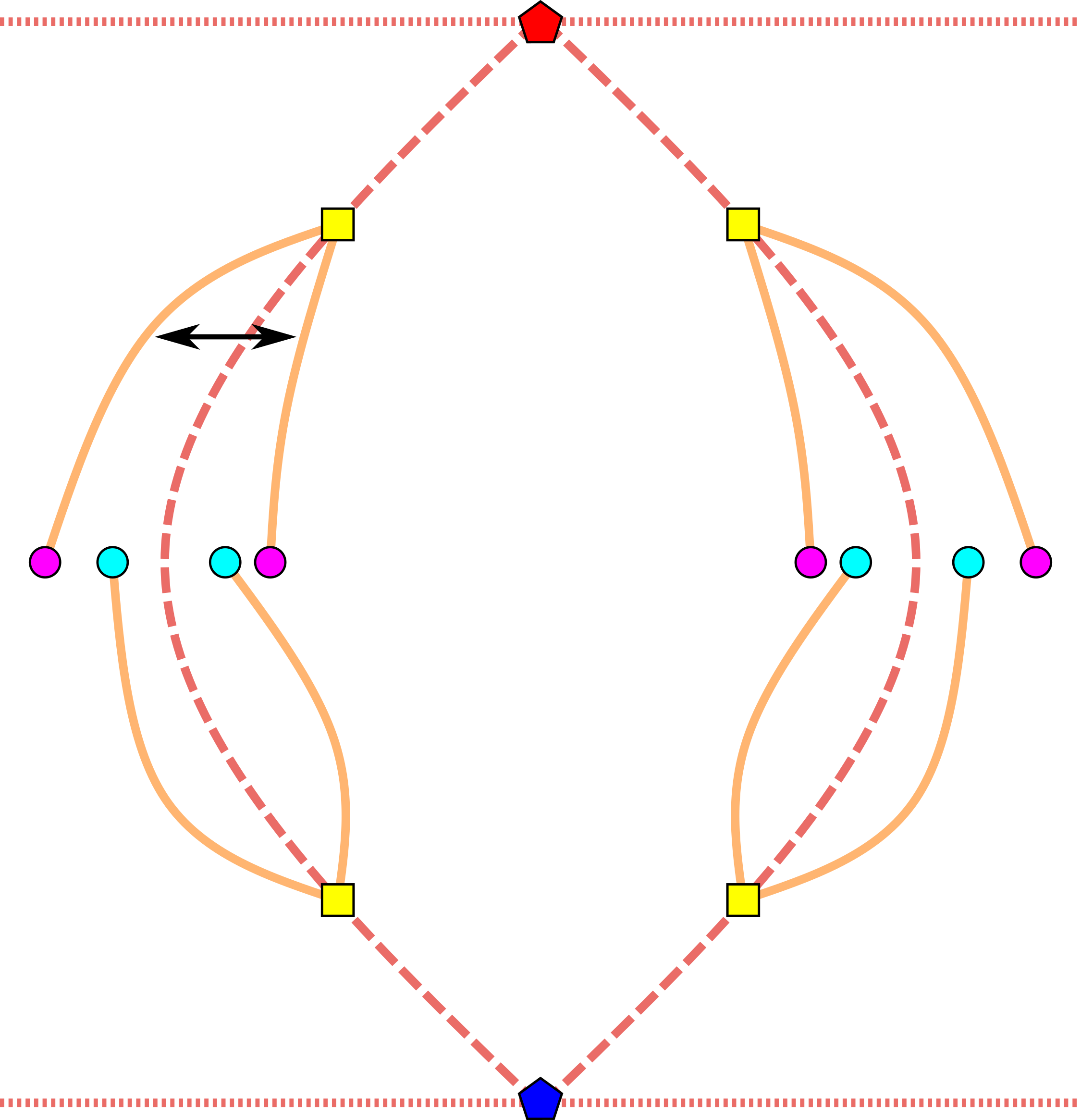}
        \caption{\label{fig:cr3bp_odd_odd}\acrshort{cr3bp}: odd-odd $p:q$.}
    \end{subfigure}
    \begin{subfigure}[b]{0.49\textwidth}
        \centering
        \includegraphics[width = 0.8\textwidth]{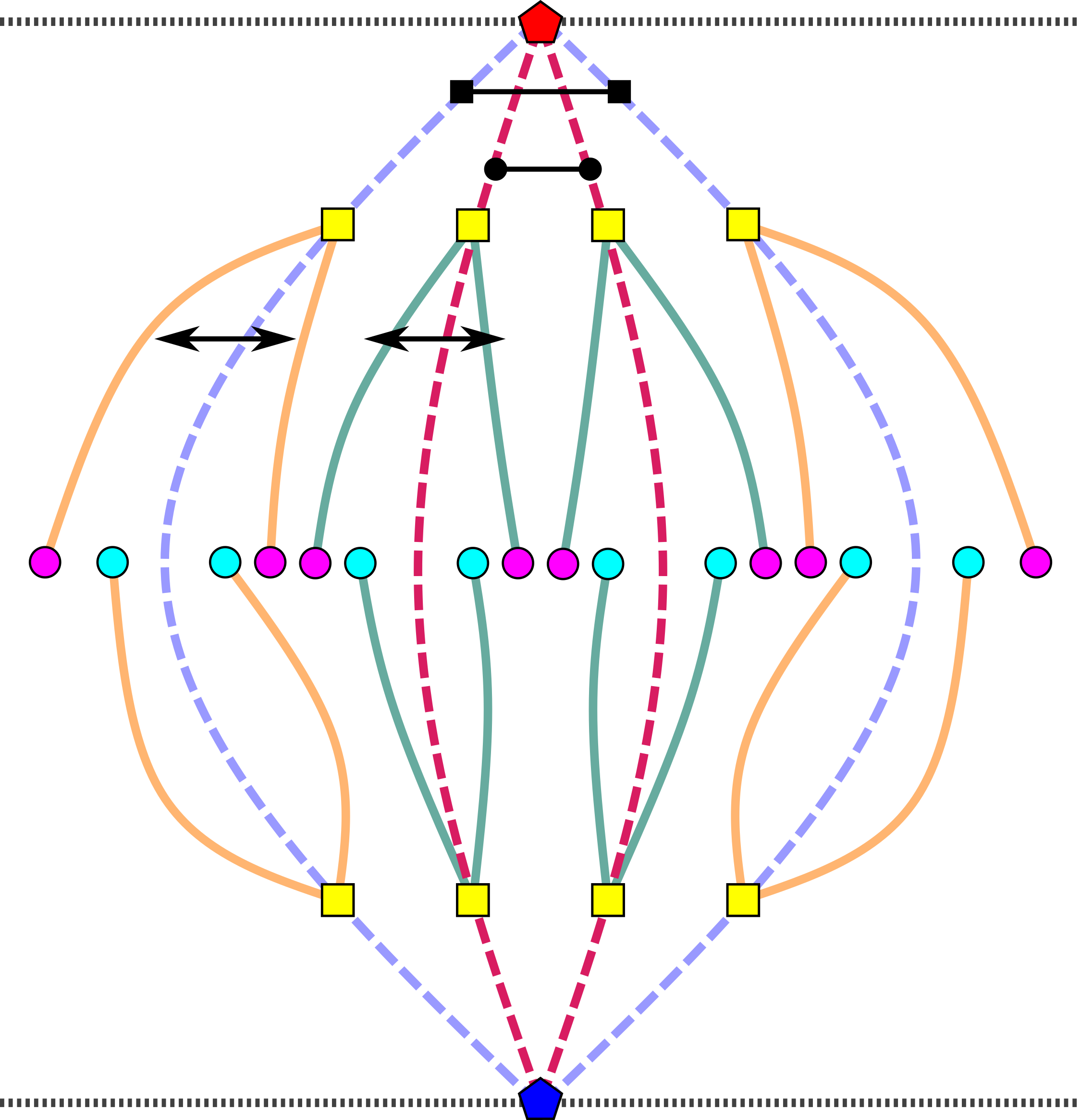}
        \caption{\label{fig:hr3bp_odd_even}\acrshort{hr3bp}: odd-even $p:q$.}
    \end{subfigure}
    \hfill
    \begin{subfigure}[b]{0.49\textwidth}
        \centering
        \includegraphics[width = 0.8\textwidth]{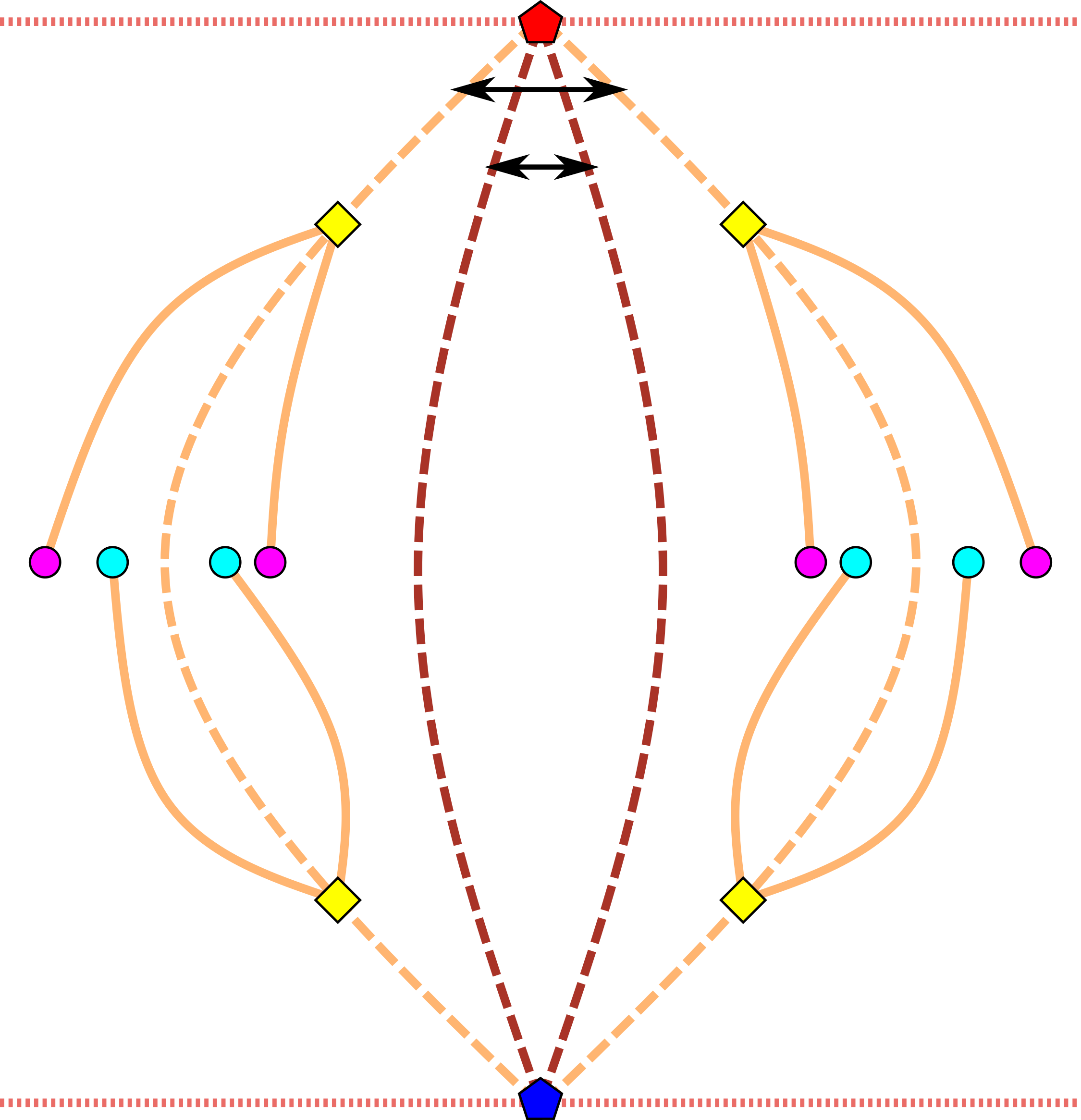}
        \caption{\label{fig:cr3bp_odd_even}\acrshort{cr3bp}: odd-even or even-odd $p:q$.}
    \end{subfigure}
    \begin{subfigure}[b]{0.49\textwidth}
        \centering
        \includegraphics[width = 0.8\textwidth]{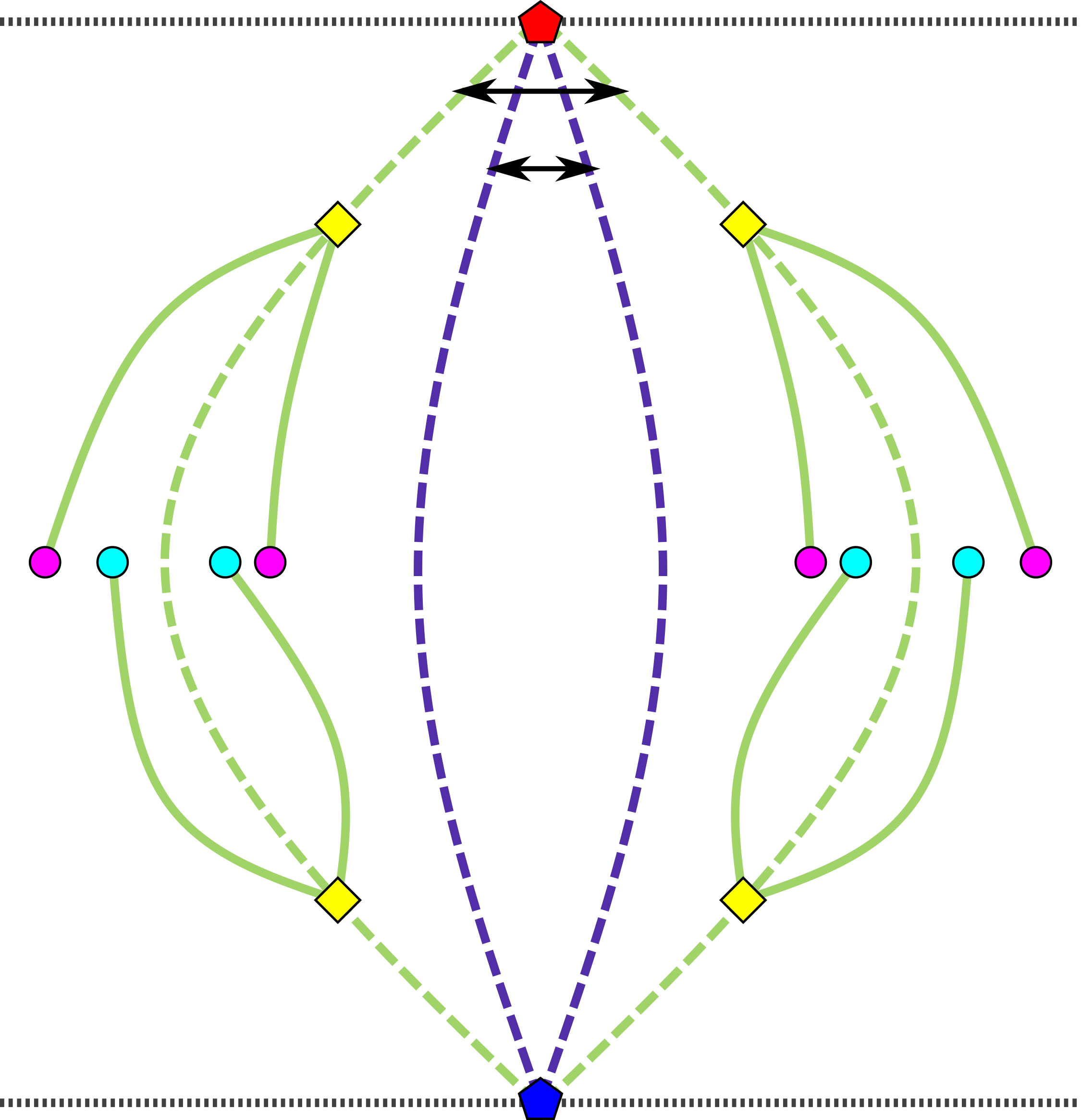}
        \caption{\label{fig:hr3bp_even_odd}\acrshort{hr3bp}: even-odd $p:q$.}
    \end{subfigure}
    \hfill
    \begin{subfigure}[b]{0.49\textwidth}
        \centering
        \includegraphics[width = 0.8\textwidth]{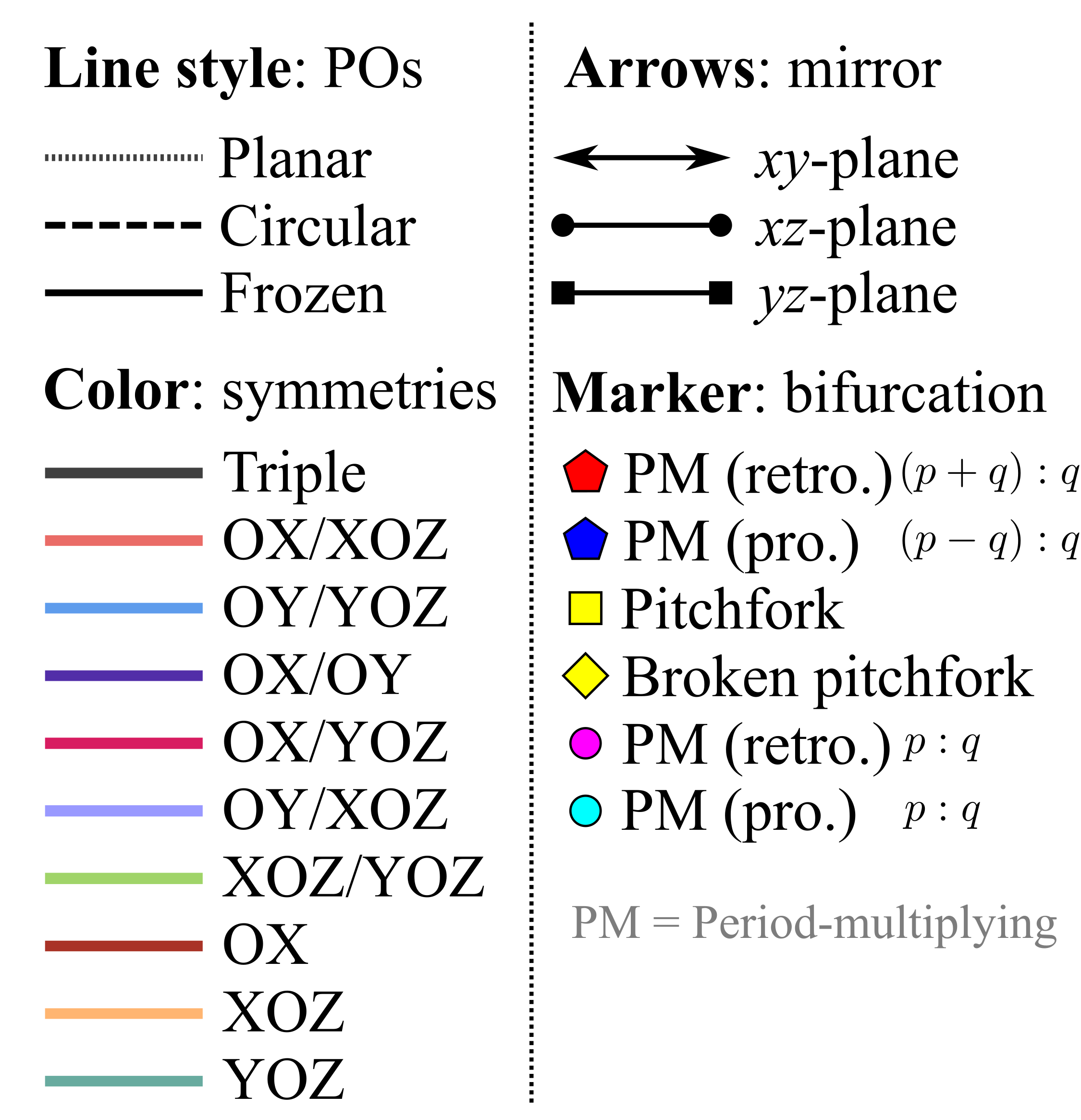}
        \caption{\label{fig:archetype_legend}Legend.}
    \end{subfigure}
    \caption{\label{fig:archetype}Archetypical bifurcation diagrams for unaveraged dynamics with different parities.}    
\end{figure}

Synthesizing the findings from the preceding sections, archetypical bifurcation diagrams are constructed to illustrate the structural connectivity within the unaveraged dynamics. Figure \ref{fig:archetype} presents these diagrams, categorized by the resonance parity and the dynamical model (\acrshort{hr3bp} vs. \acrshort{cr3bp}), with the legend provided in Fig. \ref{fig:archetype_legend}. The line styles distinguish the periodic orbit types: dotted lines represent planar PO structures, dashed lines denote circular POs, and solid lines indicate frozen POs. In each schematic, the vertical axis generally corresponds to increasing effective inclination ($i$). Consequently, the bottom and top horizontal dotted lines represent the prograde and retrograde planar \acrshort{r3bp} families, respectively, while the dashed circular POs bridge these extremal limits. The solid lines denote the frozen POs. The color coding identifies the order and specific type of symmetries maintained by each family branch. Using the odd-odd $p:q$ parity within the \acrshort{hr3bp} (Fig. \ref{fig:hr3bp_odd_odd}) as a representative example, the planar POs are triply symmetric. As predicted, four branches of doubly symmetric circular POs emerge, evenly divided into the OX/XOZ and OY/YOZ pairings (see Table \ref{tab:symmetry_evolution}). Subsequently, the frozen POs exhibit singly symmetric configurations corresponding to XOZ and YOZ types. Distinct markers highlight the bifurcation points connecting these structures. At the planar limits ($i = 0^\circ$ and $180^\circ$), the circular POs originate from the planar families via period-($p\pm q$)-multiplying bifurcations, denoted by pentagons. The bifurcation of frozen POs from circular POs at $i \approx i_{\mathrm{crit}}$ reveals a key topological distinction governed by parity and model symmetry: true pitchfork bifurcations are marked by yellow squares, whereas broken bifurcations are indicated by yellow diamonds. Approaching the polar region ($i \approx 90^\circ$), the frozen POs terminate at the collision/L$_2$ halo orbits via period-$p$-multiplying bifurcations, represented by filled circles. Note that these diagrams assume $\eta \geq 12$ to allow nominal asymptotic connections of the frozen POs to the polar limit ($i \rightarrow 90^\circ$). Finally, since certain spatial structures emerge via symmetry breaking, they exist as mirrored pairs across specific reflectional planes; these mirrored relationships are indicated by the diverging arrows.

While the specific topology in the non-integrable regime inevitably varies with the resonance ratio $\eta$ and mass parameter $\mu$ (e.g., \cite{aydin2025exploration,aydin2025studying}), these nominal diagrams provide a critical \textit{a priori} baseline. This theoretical skeleton facilitates the interpretation of the deformation of the ideal connections established in the averaged model by specific perturbations within the unaveraged dynamics. These diagrams also illustrate the linkage between the \acrshort{dadm} and the \acrshort{r3bps}, i.e., the \acrshort{hr3bp} and \acrshort{cr3bp}. 

\section{\label{sec:conclusion}Concluding Remarks} 

The current investigation establishes a unified framework linking the integrable doubly averaged dynamics to the non-integrable restricted three-body problems. By bridging the gap between model fidelity and global insight, the proposed methodology systematically maps the averaged equilibria to the symmetric periodic orbit families within the \acrfull{hr3bp} and \acrfull{cr3bp}. Central to this connection is the frequency-based resonance analysis and the identification of admissible symmetric configurations dictated by the parity of the resonance ratio. This approach enables the \textit{a priori} prediction of solution multiplicity and symmetry types, effectively converting the algebraic structure of the averaged model into a robust initialization strategy for the complex unaveraged dynamics. Through the construction of archetypical bifurcation diagrams, the global evolution of these families is elucidated. The analysis reveals that the idealized bifurcations within the averaged model often manifest as broken bifurcations or distinct period-multiplying connections in the unaveraged regime, particularly due to the symmetry-breaking nature of the Earth-Moon system. Beyond exposing the topological origins of these complex dynamical structures, the resulting atlas of symmetric periodic orbits also serves as a practical foundation for trajectory design near smaller primaries. Future efforts may extend this symmetry-aware framework to quasi-periodic orbits or higher-fidelity ephemeris models, further enhancing the applicability of dynamical insights to realistic mission scenarios.

\section*{Acknowledgments}
Dr. Beom Park acknowledges financial support from the School of Aeronautics and Astronautics at Purdue University through the Apollo 11 Postdoctoral Fellowship program. 

\bibliographystyle{unsrtnat}  
\bibliography{references}  

\clearpage

\end{document}